\documentclass[preprint,11pt]{elsarticle}
\usepackage[bookmarks=true,colorlinks=true,linkcolor=blue]{hyperref}

\pdfoutput=1

\biboptions{sort&compress}

\usepackage{amssymb,amsmath}
\usepackage{amsthm}

\newtheorem{lemma}{Lemma}
\newtheorem{proposition}{Proposition}

\usepackage{calc}
\usepackage[margin=1.in]{geometry}

\usepackage[usenames,dvipsnames,svgnames,table]{xcolor}
\newcommand{\citeCount}[1]{}

\newcommand{\dt}{\Delta t}

\newcommand{\bogus}[1]{{}}

\newcommand{\ev}{\mathbf{ e}}
\newcommand{\fv}{\mathbf{ f}}
\newcommand{\gv}{\mathbf{ g}}

\newcommand{\nv}{\mathbf{ n}}

\newcommand{\uv}{\mathbf{ u}}
\newcommand{\vv}{\mathbf{ v}}

\newcommand{\xv}{\mathbf{ x}}

\newcommand{\Fv}{\mathbf{ F}}

\newcommand{\grad}{\nabla}

\clearpage

\newcommand{\Fig}{Figure} 

\usepackage{tikz}

\newlength{\tfwidth}
\newlength{\tfheight}
\newlength{\tfxa}
\newlength{\tfxb}
\newlength{\tfya}
\newlength{\tfyb}
\newcommand{\trimFigWithBox}[6]{%
\setlength\fboxsep{0pt}%
\setlength\fboxrule{1.0pt}
\fbox{\includegraphics[width=#2, clip, trim=#3 #4 #5 #6]{#1}}%
}
\newcommand{\trimFigNoBox}[6]{%
\setlength\fboxsep{1pt}
\setlength\fboxrule{0.0pt}
\fbox{\includegraphics[width=#2, clip, trim=#3 #4 #5 #6]{#1}}%
}
\newcommand{\trimFigHeightWithBox}[6]{%
\setlength\fboxsep{0pt}%
\setlength\fboxrule{1.0pt}
\fbox{\includegraphics[height=#2, clip, trim=#3 #4 #5 #6]{#1}}%
}
\newcommand{\trimFigHeightNoBox}[6]{%
\setlength\fboxsep{1pt}
\setlength\fboxrule{0.0pt}
\fbox{\includegraphics[height=#2, clip, trim=#3 #4 #5 #6]{#1}}%
}

\newsavebox\figBox

\newcommand{\trimw}[6]{%
\sbox\figBox{\includegraphics{#1}}
\setlength{\tfwidth}{\the\wd\figBox}
\setlength{\tfheight}{\the\ht\figBox}
\setlength{\tfxa}{\tfwidth*\real{#3}}%
\setlength{\tfxb}{\tfwidth*\real{#4}}%
\setlength{\tfya}{\tfheight*\real{#5}}%
\setlength{\tfyb}{\tfheight*\real{#6}}%
\trimFigNoBox{#1}{#2}{\tfxa}{\tfya}{\tfxb}{\tfyb}%
}

\newcommand{\trimwb}[6]{%

\sbox\figBox{\includegraphics{#1}}
\setlength{\tfwidth}{\the\wd\figBox}
\setlength{\tfheight}{\the\ht\figBox}
\setlength{\tfxa}{\tfwidth*\real{#3}}%
\setlength{\tfxb}{\tfwidth*\real{#4}}%
\setlength{\tfya}{\tfheight*\real{#5}}%
\setlength{\tfyb}{\tfheight*\real{#6}}%
\trimFigWithBox{#1}{#2}{\tfxa}{\tfya}{\tfxb}{\tfyb}%
}

\newcommand{\trimh}[6]{%
\sbox\figBox{\includegraphics{#1}}
\setlength{\tfwidth}{\the\wd\figBox}
\setlength{\tfheight}{\the\ht\figBox}
\setlength{\tfxa}{\tfwidth*\real{#3}}%
\setlength{\tfxb}{\tfwidth*\real{#4}}%
\setlength{\tfya}{\tfheight*\real{#5}}%
\setlength{\tfyb}{\tfheight*\real{#6}}%
\trimFigHeightNoBox{#1}{#2}{\tfxa}{\tfya}{\tfxb}{\tfyb}%
}

\newcommand{\trimhb}[6]{%

\sbox\figBox{\includegraphics{#1}}
\setlength{\tfwidth}{\the\wd\figBox}
\setlength{\tfheight}{\the\ht\figBox}
\setlength{\tfxa}{\tfwidth*\real{#3}}%
\setlength{\tfxb}{\tfwidth*\real{#4}}%
\setlength{\tfya}{\tfheight*\real{#5}}%
\setlength{\tfyb}{\tfheight*\real{#6}}%
\trimFigHeightWithBox{#1}{#2}{\tfxa}{\tfya}{\tfxb}{\tfyb}%
}
\usepackage{algorithm} 
\usepackage{algpseudocode} 

\usepackage[subnum]{cases}

\newcommand{\pd}[2]{\frac{\partial #1}{\partial #2}}
\newcommand{\pdn}[3]{\frac{\partial^#3 #1}{\partial #2^#3}}
\newcommand{\curlcurl}{\nabla\times\nabla\times}
\newcommand{\vbase}{\mathbf{\Phi}}
\newcommand{\pbase}{\varphi}
\newcommand{\lb}{\left(}   
\newcommand{\rb}{\right)} 

\newcommand{\Pe}{\mathbb{Pe}} 

\newcommand{\oldBC}{\mathrm{TN}}
\newcommand{\newBC}{\mathrm{WABE}}
\newcommand{\norm}[1]{\|#1\|}

\newcommand{\uhat}{\hat{u}}
\newcommand{\vhat}{\hat{v}}
\newcommand{\phat}{\hat{p}}
\newcommand{\fhat}{\hat{f}}

\newcommand{\utilde}{\tilde{U}}
\newcommand{\vtilde}{\tilde{V}}
\newcommand{\ptilde}{\tilde{P}}

\newcommand{\gtilde}{\tilde{g}}

\newcommand{\bb}[1]{\mathbb{#1}} 
\newcommand{\order}[1]{\mathcal{O}\left(#1\right)} 

\begin{document}

\small

\begin{frontmatter}
\title{A split-step finite-element method for incompressible Navier-Stokes equations with high-order accuracy up-to the boundary}

\author[ul]{Longfei~Li\corref{cor1}\fnref{LongfeiThanks,RCSThanks}}
\ead{longfei.li@louisiana.edu}

%
%
%

\address[ul]{Department of Mathematics, University of Louisiana at Lafayette, Lafayette, LA  70504, USA.}


\cortext[cor1]{Department of Mathematics, University of Louisiana at Lafayette, Lafayette, LA 70504, USA}


\fntext[LongfeiThanks]{Research supported by the Margaret A. Darrin Postdoctoral Fellowship of Rensselaer Polytechnic Institute (RPI).}

\fntext[RCSThanks]{Research supported by the  Louisiana Board of Regents Support Fund under contract No. LEQSF(2018-21)-RD-A-23.}


%



\begin{abstract}
An  efficient and accurate finite-element algorithm    is described for the numerical  solution of the incompressible Navier-Stokes (INS)  
equations.  The new  algorithm that solves the INS equations in  a  velocity-pressure reformulation is  based on a split-step scheme in conjunction with the standard finite-element method. 
The split-step scheme  employed for the temporal discretization of our algorithm completely separates the pressure updates from the solution of velocity variables. When the  pressure equation is formed explicitly, the algorithm avoids solving a saddle-point problem;  therefore, our algorithm has more flexibility in choosing finite-element spaces. In contrast, popular mixed finite-element methods that  solve the INS equations in the primitive variables (or velocity-divergence formulation) lead to discrete saddle-point problems whose solution depends on  the choice of finite-element spaces for velocity and pressure that is subject to the well-known   Ladyzenskaja-Babu\v ska-Brezzi (LBB) (or $\inf$-$\sup$)  condition.
For efficiency and  robustness,  Lagrange (piecewise-polynomial) finite elements of equal order for both velocity and pressure are used. Accurate numerical boundary condition for the pressure equation is also investigated. Motivated by  a post-processing technique that calculates  derivatives of a finite element solution with super-convergent error estimates,   an alternative numerical  boundary condition is proposed for the pressure equation at the discrete level. The new numerical pressure boundary condition that can be regarded as  a better implementation of the compatibility boundary  condition  improves the boundary-layer errors of the pressure solution.  A normal-mode analysis is performed using a simplified model problem on a uniform mesh to demonstrate the numerical properties of our methods.    Convergence studies using $\Pe_1$  elements   support  the analytical results and demonstrates that our algorithm with the new   numerical boundary condition  achieves the optimal second-order  accuracy for both velocity and pressure up-to the boundary. Benchmark problems are also computed and carefully compared with existing studies. Finally, as an example to illustrate that  our approach can be easily adapted for higher-order finite elements, we solve the classical flow-past-a-cylinder problem using $\Pe_n$ finite elements with $n\geq 1$.
\end{abstract}


\end{frontmatter}

\clearpage
\tableofcontents

\section{Introduction}

The development and analysis of  numerical schemes for the fast and accurate solution of  incompressible Navier-Stokes (INS) equations have long been a very active area of research, see for example \cite{Chorin68,GiraultRaviart86,Orszag86,BCNS,ELiu95,splitStep2003,JohnstonLiu04} and the references therein. Popular existing  numerical algorithms for solving  INS equations include but are not limited to  the following pioneering methods and their follow-up variants: (i) the MAC method that uses staggered grids   for discretization \cite{Harlow65}; (ii)  projection methods \cite{Chorin68} and their extension to an implicit fractional-step method \cite{Kim85}; (iii) the method of artificial compressibility \cite{Chorin67}; (iv)  split-step methods  that solve  an equivalent velocity-pressure reformulation of the INS equations \cite{BCNS,splitStep2003,JohnstonLiu04}.  There   are also numerous other approaches based on common discretizations such as finite difference, finite element, finite volume,  spectral element, and discontinuous Galerkin method; to name just a few, see  \cite{Abdallah87a,Karniadakis,Wright93,KorczakPatera,Strikwerda,Heinrichs98,CockburnEtal05}.

The focus here  is on solving the INS equations using finite element methods (FEM).  The popular ones include methods  based on the weak formulation of  the INS equations in the primitive variables (also referred to as the velocity-divergence formulation in the literature), which often employ the $H^1(\Omega)^d\times L^2(\Omega)$ conforming elements for spatial discretization. Here $d=2$ or $3$ denotes the  spatial dimension and $\Omega$ represents the fluid domain. Methods using the  $H^1(\text{div};\Omega)\times L^2(\Omega)$ conforming finite elements  of the Raviart-Thomas type are also proposed in \cite{WangYe07} that better satisfy the divergence free condition of the fluid velocity. However, solving the INS equations in the primitive variables  leads to a discrete saddle point problem; and, in order for this problem to have a solution, the choice of finite-element spaces for velocity and pressure must satisfy the well-known   $\inf$-$\sup$ or Ladyzenskaja-Babu\v ska-Brezzi (LBB)  condition; these methods are often called mixed finite-element methods.  For example,  $\Pe_2$/$\Pe_1$ is a mixed velocity/pressure finite-element pair satisfying the LBB condition. The additional complexity posed by the LBB condition makes it hard to utilize these schemes { for more complicated multi-physics  applications  such as Fluid-Structure interaction (FSI) problems. More discussion about  popular  finite element methods for solving the INS equations can be found in the classical book by  Girault  and Raviart \cite{GiraultRaviart86}.}

Some projection methods based on the pioneering work of Chorin \cite{Chorin68} avoid solving a saddle point problem by  performing separate pressure updates through a splitting strategy. Therefore,  the solution of  these discrete problems is not subject to the LBB condition and this suggests  more flexibility in choosing finite-element spaces for spatial discretization.  However, the projection methods have their own drawback; that is, the  pressure solution near boundaries  is  found to be  less accurate with the presence of  numerical boundary-layer errors  \cite{Orszag86,ELiu95,ELiu96}. Numerous techniques  were developed attempting to reduce  these numerical boundary layers  in  modern second order projection methods \cite{ Kim85,Orszag86,vanKan86,Bell89}.  Especially,  Brown {\it et al}. \cite{brownCortezMinion:2001}  proposed a second order implicit projection method that is able to  achieve full  accuracy for both the velocity and pressure in rectangular domains;  however, within a finite-element setting, it is not straightforward that their method is applicable for more general  geometries   due to the high order spatial derivatives required in their pressure update. 
More recently, works of Liu {\it et al.} have led to the development of a series of INS algorithms using either a velocity-pressure reformulation of the original INS equations or  a variation of projection methods  based on the Laplacian and Leray projection operators; these methods have been demonstrated to achieve promising results using standard finite-element discretization  \cite{JohnstonLiu04, LiuEtal07,LiuLiuPego2010b}.

This paper  concerns  the development of  a new FEM based INS algorithm that is  suitable for FSI applications.  A long-term motivator for the work is to extend the  recently developed Added-Mass Partitioned (AMP) schemes for FSI simulations to the finite element framework. The AMP schemes were developed based on an interface condition derived at the continuous level by matching the time derivative of the kinematic interface condition. The AMP condition, which is a non-standard Robin-type boundary condition involving the fluid stress tensor, requires no adjustable parameters and in principle is applicable at the discrete level to couple the fluid and structure solvers of any accuracy and of any approximation methods (finite difference, finite element, finite volume, spectral element methods, etc). Within the finite-difference framework, the AMP algorithms have been developed and implemented to solve FSI problems involving the interaction of incompressible flows with a wide range of structures, such as elastic beams/shells \cite{fib2014,beamins2016}, bulk solids \cite{fis2014,fibrmparXiv,fibrarXiv} and rigid bodies \cite{rbinsmp2017,rbins2017,rbins3D2018}. It has been shown in these works that the AMP schemes are second-order accurate and stable without sub-time-step iterations, even for very light structures when added-mass effects are strong. In contrast, the state-of-the-art finite-element based loosely-coupled partitioned FSI algorithms have yet achieved second-order accuracy for all the solution components of the fluid-structure system \cite{CausinGerbeauNobile2005,FernandezReview2011,Degroote2013,BukacCanicMuha2015,GuidoboniGlowinskiCavalliniCanic2009,CanicMuhaBukac2012,Fernandez2011,FernandezMullaertVidrascu2013}. Given the promising results achieved using the AMP schemes within the finite-difference framework, the goal of developing  FEM based AMP schemes has the potential  to improve the accuracy and efficiency of the state-of-the-art FSI simulations within the finite-element framework. However, the extension to FEM is nontrivial, and it poses specific requirements for the underlying INS solver. With  this long-term goal in mind, the algorithm presented  in this paper is designed to deal with a number of fundamental  issues:
\begin{itemize}
\item the scheme should be able to address the checker-board instability (or correspondingly the LBB stability condition in finite elements) so that the pressure solution is free of spurious oscillations;
\item appropriate boundary conditions should be prescribed for an intermediate velocity field for projection type  methods or for the pressure if using   split-step type   methods based on a  velocity-pressure reformulation of the INS equations. These  boundary conditions are essential to keep the boundary-layer errors small, and their choices are non-trivial;
 \item the scheme should be efficient and accurate, and it is often useful to decouple the solution of the velocity from the solution of the pressure;
 \item the scheme should be able to keep the discrete divergence small if it is  not strictly enforced to be zero;
   \item the scheme should have the flexibility to work with non-standard  boundary conditions.
\end{itemize}

Motivated by  the success with the split-step  finite difference algorithm \cite{splitStep2003},   an  FEM based  split-step  algorithm  is developed   for solving a pressure-velocity reformulation of the INS equations.  The split-step strategy  that separates the  solution  of pressure from that of  the velocity variables  enables more flexibility of choosing finite-element space for spatial discretization in the same manner as projection methods by avoiding  a saddle point problem.  Therefore,    Lagrange (piecewise-polynomial) finite-elements of equal order for  both velocity and pressure can be used by our method for efficiency. It is important to note that  these standard elements   can not be used   in a straightforward way when solving the INS equations in the velocity-divergence formulation  because they   fail to satisfy the LBB  condition.  Therefore,  our  algorithm  has the potential to be more efficient than many existing FEM based fluid solvers.

Special attention has been paid to investigate   accurate numerical boundary conditions  for the pressure equation. The curl-curl boundary condition appeared  in the velocity-pressure reformulation is  a compatibility boundary condition, which is derived by applying the normal component of the momentum equations on the boundary. The correct implementation of this compatibility boundary condition as a numerical boundary condition for the pressure is found to be crucial for the stability and accuracy of similar finite difference schemes  \cite{BCNS,Petersson00,splitStep2003}.  
To this end, two approaches of incorporating the curl-curl boundary condition within the finite-element context  are considered that lead to two numerical boundary conditions of different accuracies.  The most straightforward way to incorporate the curl-curl  boundary condition, a Neumann boundary condition for the pressure, is to  implement it as a natural boundary condition within the weak formulation of the pressure-Poisson equation.  This approach, though simple and straightforward,  is found to be less accurate since a slight degradation of the pressure accuracy near the boundary is observed in numerical experiments. Motivated by a  post-processing technique that produces a super-convergent flux from finite element solutions \cite{Carey19821}, we propose an alternative way to implement the compatibility boundary condition as a more accurate numerical boundary condition for the pressure that  alleviates the boundary-layer errors  in the pressure solution.

The remainder of the paper is organized as follows. In Section~\ref{sec:VPReformulation}, we discuss  a  pressure-velocity reformulation of the INS equations. 
A split-step scheme consisting of a second-order accurate predictor-corrector algorithm  is described in Section~\ref{sec:temporalDiscretization} for the temporal discretization of the problem, and the discussion of the spatial discretization  using the standard finite-element method  follows  in Section~\ref{sec:spatialDiscretization}. The complete discrete algorithm is summarized in Section \ref{sec:algorithm}. In Section~\ref{sec:boundaryConditionsForPPE},  novel numerical boundary conditions  that  help the pressure solution achieve better accuracy near the boundary  are presented.   The numerical properties of the scheme and the two pressure boundary conditions  are analyzed for a simplified model problem discretized on a uniform mesh in Section~\ref{sec:analysis}.
Careful numerical validations are conducted  in Section~\ref{sec:results}. Finally, concluding remarks are made  in Section~\ref{sec:conclusion}.

\section{Navier-Stokes equations in velocity-pressure form}\label{sec:VPReformulation}

{An equivalent form, referred to as the ``velocity-pressure''  reformulation of the INS equations, is considered. Let $\Omega\subset \mathbb{R}^d$ ($d=2$ or $3$) denote  a bounded open domain and $\partial\Omega$ be the boundary of $\Omega$. The   velocity-pressure form are  given by
\begin{numcases}{ \label{eq:NSEvp}}
\rho\left(\pd{\uv}{t}+\uv\cdot\nabla\uv \right)= -\nabla p +\mu \Delta \uv + \Fv & $\text{for}~ \xv\in \Omega, t>0$ \label{eq:NSEvpMomentum}\\
\Delta p = - \rho\nabla\uv :\left(\nabla\uv\right)^T+\nabla\cdot \Fv + \alpha(\xv)\nabla\cdot\uv & $\text{for}~  \xv \in \Omega, t>0$\label{eq:NSEvpPPE}\\
\mathcal{B}(\uv,p) = \gv(\xv,t) & $ \text{for}~ \xv\in\partial\Omega, ~t>0$\label{eq:NSEvpBCv}\\
\nabla\cdot\uv=0 &$ \text{for}~ \xv\in\partial\Omega, ~t>0$\label{eq:NSEvpBCp}\\
\uv(\xv,0) = \fv(\xv) &  $ \text{for}~ \xv\in\Omega,~t=0$\label{eq:NSEvpIC}
\end{numcases}
where
\begin{align*}
\grad\uv:(\grad\uv)^T 
      ~ \equiv \sum_{i=1}^{d}\sum_{j=1}^{d} \frac{\partial u_i}{\partial x_j}\frac{\partial u_j}{\partial x_i}.
\end{align*}
Here $\uv=(u_1,\dots,u_d)$  and $\xv=(x_1,\dots,x_d)$  are the velocity and position in $d$-dimensional space; $p$ is the pressure; $\rho$ is the fluid density; $\mu$ is the coefficient of viscosity; and $\Fv(\xv,t) = \lb F_1(\xv,t),\dots, F_d(\xv,t) \rb$ is the external force. $ \mathcal{B}(\uv,p) = \gv(\xv,t)$ represents appropriate boundary conditions with $ \gv(\xv,t) = \lb g_1(\xv,t),\dotsm,g_d(\xv,t)\rb$ a  given vector valued function. The given initial conditions, which are assumed to be divergence free,   are represented by    $\fv(\xv)=\lb f_1(\xv),\dots,f_d(\xv)\rb$. 

Equation \eqref{eq:NSEvpPPE}
is the  pressure-Poisson equation (PPE), which is obtained by taking the divergence of the momentum equation together with the divergence free condition ($\nabla\cdot\uv=0$).  Note that a linear damping term $ \alpha(\xv)\nabla\cdot\uv$ is  included  in the PPE for numerical purposes.
The damping term, referred to as the divergence damping,   has no effect at the continuous level since $\nabla\cdot\uv=0$; however, it helps to suppress  the divergence at the discrete level, since  the numerical solution  is not exactly divergence-free due to properties of the numerical approximation. As is seen later in Section~\ref{sec:algorithm}, our numerical algorithm does not enforce the   divergence-free condition for the interior of the domain, instead this condition is only  implicitly guaranteed by the curl-curl pressure boundary condition given below in equation  \eqref{eq:curlcurlBC}.  Discretization errors could result  in  the growth of the divergence of the velocity in the numerical solution. Therefore, it is important  to include the divergence damping in the PPE to keep $\nabla\cdot\uv$  small for the whole computation. Alternatively, one could perform an extra projection after every time step to  map the velocity solutions into a divergence free space at the expense of a significant amount of additional computations.

The  velocity-pressure formulation of  the INS equations requires an extra boundary condition, and    an appropriate  choice  is to set the divergence of the velocity to be zero on the boundary, or its normal derivative \cite{splitStep2003}. {The split-step method considered in this paper  separates the updates for  velocity and   pressure; as a consequence, a Poisson problem for the pressure  needs to be solved explicitly  at each time step.
The divergence boundary condition   given in \eqref{eq:NSEvpBCp} is,  however, not convenient to be used as a  boundary condition for the PPE \eqref{eq:NSEvpPPE}  and hence an alternative condition is utilized.} Following the studies in \cite{Karniadakis,Petersson00,splitStep2003},  the curl-curl boundary condition,
\begin{align}
\pd{p}{n} =\nv\cdot\left(- \rho\pd{\uv}{t}-\rho\uv\cdot\nabla\uv  -\mu\curlcurl \uv + \Fv\right), &\quad \xv \in \partial\Omega,\label{eq:curlcurlBC}
\end{align}
has been  used in place of \eqref{eq:NSEvpBCp}  during the pressure update stage of our algorithm; noting that this curl-curl boundary condition 
is   obtained by  using  the normal component of the momentum equations \eqref{eq:NSEvpMomentum} as a compatibility boundary condition,
$$
\pd{p}{n} =\nv\cdot\left(- \rho\pd{\uv}{t}-\rho\uv\cdot\nabla\uv  +\mu \Delta \uv + \Fv\right),
$$
and then replacing the diffusion term with the following vector identity,
\begin{align*}
 \Delta\uv  = \nabla(\nabla\cdot\uv) -\curlcurl \uv. 
\end{align*}
It is important to point out that the curl-curl boundary condition  \eqref{eq:curlcurlBC}    has   the divergence free condition implicitly  implemented. 

We remark that  the  INS equations in the  velocity-pressure form were  considered  in \cite{BCNS, ICNS,splitStep2003}, in which second- and forth-order accurate finite-difference based schemes have been  developed and analyzed. In addition, the stability  analysis of the curl-curl boundary condition in the context of a centered finite difference discretization is available in \cite{Petersson00}. A similar velocity-pressure reformulation, without the divergence damping, was also investigated by Johnston and Liu \cite{JohnstonLiu04}. It is interesting to note that, without the divergence damping, the INS equations in velocity-pressure form \eqref{eq:NSEvp} are not equivalent to those in the primitive variables (also known as the  velocity-divergence form) in the case of steady-state flows because solutions with constant $\nabla\cdot \uv$ are possible for the velocity-pressure form; however, $\nabla\cdot \uv$ is forced to be identically zero if the divergence damping is present. And we will see in the analysis and results sections that including the divergence damping is essential for our scheme to achieve optimal order of accuracy.

For simplicity,   no-slip walls are assumed throughout this paper,  in which case the velocity boundary conditions \eqref{eq:NSEvpBCv} are specifically given by
$$
\uv|_{\partial\Omega}=\gv(\xv,t).
$$

}


\section{Temporal discretization}\label{sec:temporalDiscretization}

  Following \cite{splitStep2003}, we use a split-step  strategy to separate the updates of  pressure   and  velocity components in  the velocity-pressure reformulation \eqref{eq:NSEvp}. The  split-step  method is  an explicit predictor-corrector method that consists of a second-order Adam-Bashforth (AB2) predictor and a modified second-order Adam-Moulton (AM2) corrector. We note that this AB2-AM2 time-stepping method has been successfully employed to solve \eqref{eq:NSEvp} within a finite-difference framework in \cite{splitStep2003}. In this paper, we are interested in extending it to the finite-element framework. { The spatial discretization using a standard finite element method will be discussed in section \ref{sec:spatialDiscretization}, and  the summary of the full discrete algorithm follows in section \ref{sec:algorithm}.}

To be specific, the velocity-pressure formulation \eqref{eq:NSEvp} of the INS equations are advanced in time using the following time-stepping scheme. For simplicity,  the algorithm is written for a fixed time-step, $\Delta t$, so that $t_n = n\Delta t$. We note that the algorithm can be extended to  a variable $\Delta t$ with some  reasonable strategies to dynamically determine the step size and the update frequency, but this case is not considered for the scope of this paper.  Given solutions  $(\uv^{n-1},p^{n-1})$ and $(\uv^n,p^n)$ at time levels $t_{n-1}$ and $t_n$, we first predict the velocity  using the AB2 method,
  $$
  \rho \frac{\uv^p-\uv^n}{\dt} = \frac{3}{2} \left(L\uv^{n}+\Fv^{n}\right) -\frac{1}{2} \left(L\uv^{n-1}+\Fv^{n-1}\right)
~\text{with}~
  L\uv = -\rho\uv\cdot\nabla\uv -\nabla p +\mu \Delta \uv.
  $$
The pressure prediction is followed by  solving the PPE with the predicted velocity solutions,
  $$
 \Delta p^{p} = - \rho\nabla\uv^{p} :\left(\nabla\uv^{p}\right)^T+\nabla\cdot \Fv^{n+1} + \alpha(\xv)\nabla\cdot\uv^{p}.
  $$
 The velocity is  then  corrected  using the following modified AM2 method
 $$
   \rho \frac{\uv^{n+1}-\uv^n}{\dt} = \frac{1}{2} \left(L\uv^{n}+\Fv^{n}\right) +\frac{1}{2} \left(L\uv^{p}+\Fv^{n+1}\right). 
   $$
 Note that the modified AM2 method is  explicit since the predicted velocity $\uv^{p}$ is used on the right hand side of the above equation. The pressure correction follows,
   $$
    \Delta p^{n+1} = - \rho\nabla\uv^{n+1} :\left(\nabla\uv^{n+1}\right)^T+\nabla\cdot \Fv^{n+1} + \alpha(\xv)\nabla\cdot\uv^{n+1}.
   $$

   We emphasize that this algorithm is stable without the use of the corrector step. Typically,  the corrector step is included since  the scheme has a larger stability region than the predictor  alone, and the stability region includes the imaginary axis so that the scheme can be used for inviscid problems ($\mu = 0$). The time step $\Delta t$ is determined  by a diffusive stability constraint ($\Delta t \sim h^2$) for the explicit AB2-AM2 method, where $h$ is the grid spacing. Nevertheless,  the time-step restriction can be alleviated  if we treat the viscous term of the momentum equation implicitly using a Crank-Nicholson method; the time step for the  semi-implicit scheme is determined by a convective stability constraint ($\Delta t \sim h$). One can refer to \cite{Petersson00} for a proof of the time-step constraints for finite-difference based schemes. 

\section{Spatial discretization}\label{sec:spatialDiscretization}


It is important to note that the temporal discretization introduced in section \ref{sec:temporalDiscretization} fully   decouples the update of pressure from that of the velocity.
{When the pressure equation is formed explicitly,  Lagrange finite elements of equal order can be used to discretize the velocity and pressure equations in a stable manner. This is in contrast to discretizing the velocity-divergence formulation in which case the standard Lagrange basis leads to an unstable scheme that does not satisfy the LBB condition.}

In this paper,  the standard Lagrange finite elements of equal order for both velocity and pressure   are employed to discretize the  above time-difference scheme in space; that is, we look for finite element   solutions in  the finite dimensional space  $V_{h}^d \times Q_h =\bb{P}_{n}^d \times \bb{P}_n$. Here
$\bb{P}_n$ is the piecewise polynomial finite-element space of degree $n$ that is defined on a triangulation of the domain $\mathcal{T}_h(\Omega)$.
If 
 $\{\pbase_j\}_{j=1}^N$  denote the basis functions of $\bb{P}_n$ with  $N$ being the number of degrees of freedom, then  a set of basis functions for    $\bb{P}_n^d$ can be conveniently formed as $
 {\vbase_k}_j = \pbase_j \ev_k
$ ( $k=1,\dots,d$ and $j = 1,\dots,N$), where $\ev_i$'s are the standard bases of $\bb{R}^d$.  Thus, the finite-element approximations to velocity and pressure solutions at time $t_n$ can be represented as
\begin{equation}\label{eq:spatialRep}
\uv^n_h = \sum_{k=1}^d\sum_{j=1}^{N}u^n_{k_{j}}{\vbase_k}_j,~~p^n_h = \sum_{j=1}^{N}p^n_{j}\pbase_j.
\end{equation}

\newcommand{\indentI}{\hangindent\parindent\hangafter=0\noindent}

\section{The complete  numerical  algorithm}\label{sec:algorithm}

To simplify discussion,  the following notations are introduced for inner products defined over the domain $\Omega$ and its boundary $\partial\Omega$,  
\begin{alignat*}{3}
\forall\, f,g \in L_2(\Omega):       &\quad\left(f,g\right) = \int_{\Omega} fg \,dX,              &&\quad\left<f,g\right> = \int_{\partial\Omega} fg\,dS,\\
\forall\, \fv,\gv \in  L_2(\Omega)^d: &\quad\left(\fv,\gv\right) = \int_{\Omega} \fv\cdot\gv \,dX, &&\quad\left(\nabla\fv,\nabla\gv\right) = \int_{\Omega} \nabla\fv:\nabla\gv \,dX.
\end{alignat*}
And we denote   $V_{h\mathbf{0}}^d$  the subspace of  $V_{h}^d$ that vanishes on the boundary.

Given the finite element solutions $(\uv^n_h,p_h^n) \in V_{h}^d \times Q_h$ at the current time $t_n$, and $(\uv^{n-1}_h,p_h^{n-1}) \in V_{h}^d \times Q_h$ at one previous time level $t_{n-1}$, the goal of the  algorithm is to determine the solution at time $t_{n+1}$.  The complete discrete scheme using the above predictor-corrector time stepping method and the FEM spatial discretization is  as follows.

\medskip
\noindent{\bf Begin predictor.}

\medskip
\indentI {\bf Stage I - velocity prediction}: 
We predict the velocity solution by solving for   $\uv^p_h\in V_{h}^d$ such that
 \begin{equation}
   \begin{cases}
          \medskip\displaystyle
     \uv_h^p(\xv_b) =\gv^{n+1}(\xv_b),& \forall \xv_b \in \partial \Omega, \\
     \medskip\displaystyle
  \frac{\rho}{\dt} \left(\uv_h^p-\uv_h^n,\vv_h\right) = \frac{3}{2}\left[ \left(L\uv_h^{n},\vv_h\right)+(\Fv^{n},\vv_h\right)] -\frac{1}{2} \left[(L\uv_h^{n-1},\vv_h)+(\Fv^{n-1},\vv_h)\right], &  \forall \vv_h \in V_{h\mathbf{0}}^d.
      \end{cases}
 \end{equation}
 Here
 $
 (L\uv_h,\vv_h)  =-\rho( \uv_h\cdot\nabla\uv_h,\vv_h) -(\nabla p_h,\vv_h) -\mu (\nabla \uv_h,\nabla\vv_h)$. The superscripts  over the given functions $\Fv$ and $\gv$ indicate evaluating the functions  at the corresponding time level.

\medskip
\indentI {\bf Stage II - pressure update}: 
It is important to point out that, 
with the no-slip boundary condition, the pressure is only determined up-to an additive constant in the PPE (other boundary conditions may remove this singularity). {For those boundary conditions  that imply a  singular Poisson problem}, we add an additional constraint $ \lb p_h,1\rb=0$ that  sets the mean value of pressure to zero  to the PPE  as a Lagrange multiplier \cite{splitStep2003}. Specifically, we solve for $p^p_h\in Q_h$ and $\lambda\in\mathbb{R}$ such that,  for $ \forall\, q_h\in Q_h$,
\begin{equation}\label{eq:discretePPEwithRegularization}
\begin{cases}
  \medskip\displaystyle
  -\lb\nabla p^p_h , \nabla q_h\rb+\lambda \lb1,q_h\rb= \left( - \rho\nabla\uv_h^{p} :\left(\nabla\uv_h^{p}\right)^T+\nabla\cdot \Fv^{n+1} + \alpha\nabla\cdot\uv_h^{p},q_h\right) - \left<\pd{p^p_h}{n},q_h\right>, \\
  \medskip\displaystyle
  (p^p_h,1)=0.
\end{cases}
\end{equation}
In practice, we chose $\alpha$ to be inversely proportional to the square of the mesh spacing. For a nonuniform mesh, the  minimum of the spacings, $h_{\min}$, is used; that is,
$
\alpha=C_dh_{\min}^{-2}.
$
In addition,  the boundary integral is given by
\begin{equation}
\left<\pd{p^p_h}{n},q_h\right> = \left< \nv\cdot\left(-\rho\pd{\gv^{n+1}}{t}-\rho\gv^{n+1}\cdot\nabla\uv^p_h+\Fv^{n+1}\right),q_h\right> +\mu \left<\nabla\times\uv^p_h,\nv\times\nabla q_h\right>,\label{eq:discreteWeakFormulationNBC} 
\end{equation}
which is derived by utilizing the curl-curl boundary condition \eqref{eq:curlcurlBC} and the following vector identity, 
\begin{equation}\label{eq:vectorIdentity}
  \left< \nv\cdot\curlcurl\uv,q\right> = \left(\curlcurl\uv,\nabla q \right) = -\left<\nabla\times\uv,\nv\times\nabla q\right>.
  \end{equation}
We note that the vector identity \eqref{eq:vectorIdentity} plays an important role in the algorithm since it reduces   the regularity requirement for the admissible finite-element space that makes it possible to use $\mathbb{P}_1$ finite elements.

\medskip
\indentI     {\bf Remark}:  we refer to the pressure boundary condition \eqref{eq:discreteWeakFormulationNBC}  as the traditional Neumann (TN) boundary condition, which arises naturally from  testing the PPE and integration by parts. However, a slight degradation of the pressure accuracy near the boundary is observed in numerical experiments. We think this degradation stems from the fact that we need to evaluate  $\left<\nabla\times\uv_h,\nv\times\nabla q_h\right>$ for the TN pressure boundary condition, and  direct evaluation of $\nabla\times\uv_h$ from the finite element solution for the velocity field $\uv_h$ is of sub-optimal order of accuracy. Motivated by Carey's post-processing technique that produces a super-convergent  flux from finite element solutions \cite{Carey19821}, we propose an alternative compatibility boundary condition for the pressure at the discrete level. This  new condition referred to as WABE boundary condition  improves the pressure accuracy near the boundary, and its detail will be discussed in Section~\ref{sec:boundaryConditionsForPPE}.    

\medskip
\noindent{\bf Begin corrector.}

\medskip
\indentI {\bf Stage III - velocity correction}: To correct the velocity, we solve for $\uv^{n+1}_h\in V_{h}^d$ such that
 \begin{equation}
   \begin{cases}
          \medskip\displaystyle
     \uv_h^{n+1}(\xv_b) =\gv^{n+1}(\xv_b),& \forall \xv_b \in \partial \Omega, \\
     \medskip\displaystyle
  \frac{\rho}{\dt} \left(\uv_h^{n+1}-\uv_h^n,\vv_h\right) = \frac{1}{2}\left[ \left(L\uv_h^{n},\vv_h\right)+(\Fv^{n},\vv_h\right)] +\frac{1}{2} \left[(L\uv_h^{p},\vv_h)+(\Fv^{n+1},\vv_h)\right], &  \forall \vv_h \in V_{h\mathbf{0}}^d.
      \end{cases}
 \end{equation}

\medskip
\indentI {\bf Stage IV - pressure update}: Finally, we obtain the updated pressure solution $p^{n+1}_h\in Q_h$ by solving the same equation \eqref{eq:discretePPEwithRegularization} and \eqref{eq:discreteWeakFormulationNBC} but with the corrected velocity solution $\uv_h^{n+1}$ supplied on the right hand side.

\section{Weighted average over boundary elements (WABE) boundary condition} \label{sec:boundaryConditionsForPPE}


{In this section, an alternative numerical boundary condition for the PPE is introduced, which is essential for our algorithm to achieve high-order accuracy up-to the boundary.
We derive this numerical  boundary condition from   a weighted average of the momentum equations  over boundary elements; and therefore, we refer to it as WABE boundary condition for short. The WABE boundary condition avoids direct evaluation of $\nabla\times\uv_h$ in leading order terms, and  thus prevents  degradation of  the pressure accuracy from the boundary. Unless otherwise noted, the discussion is restricted to 2D in this paper.}


Using the continuity equation in the Laplacian term of the momentum equations implies
\begin{align}
\rho\lb\pd{{u_1}}{t}+\uv\cdot \nabla {u_1}\rb = -\pd{p}{x_1}+\mu \lb \pdn{u_1}{x_2}{2}- \frac{\partial ^2u_2}{\partial x_1\partial x_2}\rb ,  \label{eq:momemtumEqn1}\\
\rho\lb\pd{{u_2}}{t}+\uv\cdot \nabla {u_2}\rb = -\pd{p}{x_2}+\mu \lb \pdn{u_2}{x_1}{2}- \frac{\partial ^2u_1}{\partial x_1\partial x_2}\rb. \label{eq:momemtumEqn2} 
\end{align}
Not that  the external forcing $\Fv$ is omitted here to save space, which can be easily included in the WABE boundary condition derived below.
On the finite dimensional space $V_h^d\times Q_h$,  testing the momentum equations \eqref{eq:momemtumEqn1}  and \eqref{eq:momemtumEqn2} with the basis function $\pbase_{i_b}$  and integrating by parts, we have
\begin{align}
\left(\pd{p_h}{x_1}, \varphi_{i_b}\right)& = -\rho\left(\pd{u_{1_h}}{t}+\uv_h\cdot\nabla u_{1_h},\varphi_{i_b}\right)-\mu\left(\pd{u_{1_h}}{x_2}-\pd{u_{2_h}}{x_1},\pd{\varphi_{i_b}}{x_2}\right)+n_2\mu\left<\pd{u_{1_h}}{x_2}-\pd{u_{2_h}}{x_1},\varphi_{i_b}\right>,  \label{eq:averageMomemtumEqn1}\\
\left(\pd{p_h}{x_2}, \varphi_{i_b}\right)& =-\rho\left(\pd{u_{2_h}}{t}+\uv_h\cdot\nabla u_{2_h},\varphi_{i_b}\right)-\mu\left(\pd{u_{2_h}}{x_1}-\pd{u_{1_h}}{x_2},\pd{\varphi_{i_b}}{x_1}\right)+n_1\mu\left<\pd{u_{2_h}}{x_1}-\pd{u_{1_h}}{x_2},\varphi_{i_b}\right>. \label{eq:averageMomemtumEqn2}
\end{align}
Here $i_b$ denotes the index of any boundary degree of freedom, i.e.,  $\forall P_{i_b}\in\partial\Omega$.  
Since $\pbase_{i_b}$ is nonzero only on the elements that contain the node $P_{i_b}$, the equations \eqref{eq:averageMomemtumEqn1} and \eqref{eq:averageMomemtumEqn2} are essentially a weighted average of the momentum equations over those boundary elements. Let $\nv^{i_b} = \left( n_1^{i_b},n_2^{i_b}\right)$ be the unit outward normal vector at the boundary node $P_{i_b}$.  To mimic the curl-curl boundary condition \eqref{eq:curlcurlBC}, the averaged momentum  equations  \eqref{eq:averageMomemtumEqn1} and  \eqref{eq:averageMomemtumEqn2} are combined  in $\nv^{i_b}$ direction, and the WABE boundary condition is obtained:
\begin{align}
\left(\nv^{i_b}\cdot\nabla p_h, \varphi_{i_b}\right) = -\rho\left(\nv^{i_b}\cdot \lb \pd{\uv_h}{t}+\uv_h\cdot\nabla {\uv}_h\rb,\varphi_{i_b}\right)+\mu\left(\nabla\times\uv_h,\nv^{i_b}\times\nabla\varphi_{i_b}\right)+I_b,
\label{eq:WABEbc}
\end{align}
where
$
I_b=\mu\left<\lb\nv\times\nv^{i_b}\rb\cdot \lb\nabla\times\uv_h\rb,\varphi_{i_b}\right> 
$ accounts for the boundary integral resulted from the integration by parts.
Here $\nv$ is the normal on the edge and $\nv^{i_b}$ is the normal on the boundary node. The less accurate boundary integral $I_b$ vanishes if the  boundary of the domain is a straight line since  $\nv\times\nv^{i_b}=0$;  however,  if the boundary is not a straight line but a smooth curve, we  have $\nv\times\nv^{i_b}=\order{h^2}$ and thus the error introduced by $\nabla\times\uv_h$ will be scaled down by an order of $h^2$. In either case,  the WABE boundary condition \eqref{eq:WABEbc} should be more accurate than the discrete TN boundary condition \eqref{eq:discreteWeakFormulationNBC}, and thus improves the boundary-layer errors that appear in  the pressure solution. We note that $I_b$ can be ignored for second order accurate methods in practice.

To solve for  pressure  with the WABE  boundary condition, we only need to  replace the discrete pressure equations on the boundary nodes with the equations given by the WABE boundary condition \eqref{eq:WABEbc}.  To be specific, we solve the following modified discrete PPE at the stages of   pressure update as described in the previous section,
\begin{equation}\label{eq:discretePPEwithRegularizationWABE}
\begin{cases}
\medskip \displaystyle
\text{For $\forall i$ such that $P_i\in\Omega\backslash \partial \Omega$  and $\forall i_b$ such that  $P_{i_b}\in\partial \Omega$}, \\
\medskip \displaystyle
-\lb\nabla p^{n+1}_h , \nabla \varphi_i\rb+\lambda \lb1,\varphi_i\rb= \left( - \rho\nabla\uv_h^{n+1} :\left(\nabla\uv_h^{n+1}\right)^T+ \alpha(\xv)\nabla\cdot\uv_h^{n+1},\varphi_i\right),\\
\medskip \displaystyle
 \left(\nv^{i_b}\cdot\nabla p^{n+1}_h, \varphi_{i_b}\right) = -\rho\left(\nv^{i_b}\cdot \lb \pd{\uv^{n+1}_h}{t}+\uv^{n+1}_h\cdot\nabla {\uv}_h^{n+1}\rb,\varphi_{i_b}\right)+\mu\left(\nabla\times\uv^{n+1}_h,\nv^{i_b}\times\nabla\varphi_{i_b}\right),\\
\medskip   \displaystyle (p^{n+1}_h,1)=0.\\
 \end{cases}
\end{equation}
It is important to remark that the implementation of the WABE boundary condition is not computationally more expensive than the TN boundary condition since there is no additional calculation needed to facilitate  the WABE condition, and  all the data needed,  except for ${\partial\uv^{n+1}_h}/{\partial t}$,   are directly available from the previous velocity updates. 
However, a sufficiently accurate value for the purpose of the boundary condition can be obtained using a forward finite difference formula in time, i.e.,
$
{\partial\uv^{n+1}_h}/{\partial t}= (\uv_h^{n+1}-\uv_h^{n})/\Delta t
$.

\section{A model problem for analysis}\label{sec:analysis}
In this section, we perform  a normal-mode  analysis on  a  model problem to reveal the numerical properties of the algorithm.  The motivation that we  analyze  the model problem for a particular geometric domain using finite difference theory as opposed to using a more traditional approximation theory and energy estimation  is as follows.  As we know, the standard  energy estimate method  concerns  the accuracy results in  $L_2$ norm, which averages the numerical errors  over the entire domain. However, one of the main focuses of this paper is to address the  boundary-layer errors for the pressure solution; therefore, it is more appropriate to analyze numerical errors in $L_\infty$ norm.  Using  a normal-mode analysis for a particular geometric domain, we essentially solve the model problem analytically. Even though  this analysis does not apply to general geometries and unstructured meshes, with the analytical solution of a model problem, we are able to identify some subtleties that are buried inside the $L_2$ averaging process, i.e., the boundary-layer errors in the pressure solution. The idea of rewriting finite element schemes  on a uniform mesh as finite difference ones for analytical
purposes was also employed by other researchers; for instance, in [47], the authors rewrote their discontinuous
Galerkin schemes as finite difference ones, and then performed a Fourier type analysis since it was not easy for them to use standard finite element techniques to prove the inconsistency and weak instability of their schemes. Utilizing a  mode analysis, they were able to analytically identify the  weak instability that was observed in their numerical experiments.

Here  Stokes equations are used as a model problem for the INS equations since the nonlinear convection terms can be regarded as lower order terms that   do not affect the stability of the  scheme. To further  simplify the discussion, the  model problem  is assumed to be  $2\pi$-periodic in $x$ direction on a semi-infinite domain $[0,2\pi]\times[0,\infty]$. Although the model problem drastically simplifies the original INS equations, the analysis performed here can shed some light on the stability and accuracy of our proposed method. Specifically, we consider the following initial-boundary value problem (IBVP)  with  $\rho=1$ and $\nu=\mu/\rho$, 
\begin{align*}
 \pd{u}{t} = -\pd{p}{x}+\nu\Delta u+ f_u, &\qquad \text{for}~ (x,y)\in [0,2\pi]\times[0,\infty],\\
  \pd{v}{t} = -\pd{p}{y}+\nu\Delta v+ f_v, &\qquad \text{for}~ (x,y)\in [0,2\pi]\times[0,\infty],\\
  \Delta p = \nabla \fv  +\alpha \nabla\cdot\uv, & \qquad \text{for}~ (x,y)\in [0,2\pi]\times[0,\infty].
\end{align*}
We impose  the no-slip and curl-curl pressure boundary conditions at the boundary $y=0$,
\begin{align*}
  u(x,0,t)=v(x,0,t)=0 \quad \text{and}\quad
  &\pd{p}{y} = -\nu \frac{\partial^2 u}{\partial x\partial y}+f_v. \label{eq:curlcurlAnalysis}
\end{align*}
The homogeneous initial conditions are  imposed to complete the statement of the model problem,
\begin{align*}
  u(x,y,0)=v(x,y,0)=p(x,y,0)=0.
\end{align*}

Utilizing the  assumption of periodicity,    the IBVP can be  Fourier transformed  in $x$ direction; thus, in the Fourier space, we have 
\begin{equation}\label{eq:transformedPDEs}
  \begin{cases}
 \medskip\displaystyle\pd{\uhat}{t} = - ik \phat-\nu k^2\uhat+\nu\pdn{\uhat}{y}{2}+ \fhat_u,\\
 \medskip\displaystyle \pd{\vhat}{t} = -\pd{\phat}{y}- \nu k^2\vhat+\nu\pdn{\vhat}{y}{2}+ \fhat_v,\\
 \medskip\displaystyle -k^2\phat+\pdn{\phat}{y}{2} =  ik\fhat_u+\pd{\fhat_v}{y}  +\alpha ik  \uhat+ \alpha\pd{\vhat}{y},\\
  \end{cases}
\end{equation}
subject to the transformed  boundary and initial conditions,
\begin{equation}\label{eq:transformedICBCs}
  \begin{cases}
     \medskip\displaystyle\uhat(k,0,t)=\vhat(k,0,t)=0,\\
  \medskip\displaystyle\pd{\phat}{y}(k,0,t) = -\nu ik\frac{\partial \uhat}{\partial y}(k,0,t)+\fhat_v(k,0,t) \\
    \medskip\displaystyle   \uhat(k,y,0)=\vhat(k,y,0)=\phat(k,y,0)=0.
      \end{cases}
\end{equation}
Here $\uhat(k,y,t)$, $\vhat(k,y,t)$, and $\phat(k,y,t)$  are Fourier transformations of $u(x,y,t)$, $v(x,y,t)$   and $p(x,y,t)$ with wave number $k \in \mathbb{Z}$.
For regularity, we also require
\begin{equation}\label{eq:infinityConditions}
\norm{\uhat} <\infty,~\norm{\vhat} <\infty,~\text{and}~\norm{\phat} <\infty,
\end{equation}
where the norm is the standard $L_2$ function norm.

To analyze our numerical scheme, we discretize the above transformed equations  using $\Pe_1$ finite elements on   a  uniform Cartesian grid,  $\mathcal{G}=\{y_j = jh~ |~ \forall j\in \mathbb{N}\} $. Here $j=0$ corresponds to the boundary node. Thus, the finite element approximations to $\uhat(k,y,t)$, $\vhat(k,y,t)$, and $\phat(k,y,t)$ can  be represented as
\begin{equation*}
\uhat_h(k,y,t) = \sum_{j=0}^{\infty} u_j(k,t) \varphi_j(y),\quad \vhat_h(k,y,t) = \sum_{j=0}^{\infty} v_j(k,t) \varphi_j(y),\quad\text{and}\quad \phat_h(k,y,t) = \sum_{j=0}^{\infty} p_j(k,t) \varphi_j(y),
\end{equation*}
with $\varphi_j(y)$ denoting the $j$th basis function for $\Pe_1$ that is defined by 
$$
\varphi_j(y)=
\begin{cases}
  \frac{1}{h}(y-y_{j-1}), & y\in [y_{j-1},y_j]\\
  -\frac{1}{h}(y-y_j), & y\in [y_j,y_{j+1}]\\
  0, & \text{otherwise}
\end{cases},
~\text{for}~j>0,
\quad\text{and}\quad
\varphi_0(y) =
\begin{cases}
  -\frac{1}{h}(y-y_0), & y\in[y_0,y_1]\\
 0,& \text{otherwise}
  \end{cases}.
$$

For any  grid function $f_j$, we introduce the following difference  operators 
$$ 
Mf_j = \frac{1}{6}f_{j-1}+\frac{2}{3}f_{j} +\frac{1}{6} f_{j+1},\quad
D_{+}f_j=\frac{f_{j+1}-f_{j}}{h},\quad 
D_{-}f_j=\frac{f_{j}-f_{j-1}}{h},\quad\text{and}\quad
D_{0}f_j=\frac{f_{j+1}-f_{j-1}}{2h},
$$
where $M$ is an average operator, and  $D_{+}$, $D_{-}$ and $D_0$ are the forward, backward and centered divided difference operators, respectively.
With these difference operators, it is easily seen that 
\begin{align*}
  (\uhat_h,\varphi_j) 
  = h Mu_j, ~
  \left(\pd{\uhat_h}{y},\varphi_j\right) 
  =h D_0u_j,~ \text{and}~
  -\left(\pd{\uhat_h}{y},\pd{\varphi_j}{y}\right)  
  =hD_+D_-u_j, & \quad\forall j>0.
\end{align*}
Similar relations can be found for $\vhat_h$ and $\phat_h$ as well. Therefore, the finite element discretization of the transformed equations can be regarded   as a finite difference scheme, and thus can be analyzed  using  finite difference techniques. 

Specifically, we rewrite our finite element scheme as the following finite difference scheme,
\begin{align}\label{eq:fdEqn}
  \text{for}~j>0:~
  \begin{cases}
 M\dot{u}_j = - ik Mp_j-\nu k^2Mu_j+\nu D_+D_-u_j+ (\fhat_u,\varphi_j)/h,\\
  M\dot{v}_j = -D_0p_j- \nu k^2Mv_j+\nu D_+D_-v_j+ (\fhat_v,\varphi_j)/h, \\
     -k^2Mp_j+D_+D_-p_j=  \alpha ik  Mu_j+ \alpha D_0v_j  +\left( ik\fhat_u+\pd{\fhat_v}{y},\varphi_j\right)/h,
    \end{cases}
\end{align}
subject to  the no-slip boundary conditions
$u_0=v_0=0$,
and either one of  the   discrete TN \eqref{eq:modelEqnTN} and  WABE \eqref{eq:modelEqnWABE} boundary conditions  that are described below.  Here  $\alpha=C_d/h^2$ is assumed where $C_d$ is a constant independent of $h$.

To derive the discrete TN boundary condition for the model problem, we  test  the pressure equation in \eqref{eq:transformedPDEs} with $\varphi_0$; that is
$$
-k^2(\phat_h,\varphi_0)-\left(\pd{\phat_h}{y},\pd{\varphi_0}{y}\right)+\left.\pd{\phat_h}{y}\varphi_0\right|_0^\infty =    \alpha ik  (\uhat_h,\varphi_0)+ \alpha\left(\pd{\vhat_h}{y},\varphi_0\right)+\left(ik\fhat_u+\pd{\fhat_v}{y},\varphi_0\right).
$$
With the regularity condition \eqref{eq:infinityConditions} and the curl-curl condition  given in \eqref{eq:transformedICBCs}, we have
$$
\left.\pd{\phat_h}{y}\varphi_0\right|_0^\infty  = -\pd{\phat_h}{y}(k,0,t)\varphi_0(0)=\nu ik\pd{\uhat_h}{y}(k,0,t) -\fhat_v(k,0,t) =\nu ik \frac{u_1-u_0}{h} -\fhat_{v_0},
$$
thus the TN boundary condition can be written as
\begin{equation}\label{eq:modelEqnTN}
  D_+p_0 +\nu ik D_+u_0 =\fhat_{v_0}+k^2\left(\frac{1}{3}p_0+\frac{1}{6}p_1\right)h+   \alpha ik \left(\frac{1}{3}u_0+\frac{1}{6}u_1\right)h+\alpha\frac{v_1-v_0}{2}+\left(ik\fhat_u+\pd{\fhat_v}{y},\varphi_0\right). 
\end{equation}

On the other hand, the WABE condition for the transformed model problem is given by 
$$
\left(\pd{\phat}{y},\varphi_0\right) =-\left(\pd{\vhat}{t},\varphi_0\right)-\nu k^2\left(\vhat,\varphi_0\right) -\nu ik\left(\pd{\uhat}{y},\varphi_0\right) +(\fhat_v,\varphi_0),
$$
which can be written as the following finite difference form,

\begin{equation}\label{eq:modelEqnWABE}
D_+p_0+\nu ikD_+u_0=\frac{2}{h}(\hat{f}_v,\varphi_0)-\left(\frac{2}{3}\dot{v}_0+\frac{1}{3}\dot{v}_1\right)-\nu k^2\left(\frac{2}{3}{v}_0+\frac{1}{3}{v}_1\right).
\end{equation}

\subsection{Consistency and Error Equations}
Assuming that the continuous problem  given in \eqref{eq:transformedPDEs} and \eqref{eq:transformedICBCs} has a smooth solution, $(\uhat,\vhat,\phat)$, we introduce it into the difference equations  \eqref{eq:fdEqn}  and obtain  that the order of the truncation error for the discretization  is $\order{h^2}$. In addition,  introducing the smooth solution into the boundary conditions \eqref{eq:modelEqnTN} and \eqref{eq:modelEqnWABE},  we  obtain the orders of the   truncation errors for the  TN and WABE boundary conditions are $\order{h}$ and $\order{h^2}$, respectively. 
Furthermore, we can readily write down the equations for the errors $U_i=u_i-\uhat_i$, $V_i=v_i-\vhat_i$ and $P_i=p_i-\phat_i$, where  $\uhat_j=\uhat(k,y_i,t)$, $\vhat_j=\vhat(k,y_i,t)$ and $\phat_j=\phat(k,y_i,t)$ are the exact solutions evaluated at  $y_j$ for wavenumber $k$. The consistency and error equations of the discretized problem are summarized in the following proposition; more detail of the analysis can be found in  \ref{sec:appendix}
\begin{proposition}\label{prop:errorEqns}
  The   numerical scheme \eqref{eq:fdEqn} subject to either the TN \eqref{eq:modelEqnTN} or the WABE \eqref{eq:modelEqnWABE} boundary conditions is consistent. The discretization is formally $\order{h^2}$ for the interior equations,   $\order{h}$ for the TN boundary condition and  $\order{h^2}$ for the WABE boundary condition. The error equations for the discretized problem are given by
\begin{align}\label{eq:errorEqnLumped}
  \text{for}~j>0:~
  \begin{cases}
 \dot{U}_j = - ik P_j-\nu k^2U_j+\nu D_+D_-U_j+h^2F_u,\\
  \dot{V}_j = -D_0P_j- \nu k^2V_j+\nu D_+D_-V_j+h^2F_v, \\
     -k^2P_j+D_+D_-P_j=  \alpha ik  U_j+ \alpha D_0V_j +h^2F_p,
  \end{cases}
\end{align}
where $F_u$, $F_v$ and $F_p$ are some functions of $\order{1}$, and the boundary errors are  given by 
\begin{equation}\label{eq:errorBCLumped}
  \begin{cases}
    U_0=V_0=0\\
    D_+P_0+\nu ikD_+U_0=h^rg_0,\\
   \end{cases}
\end{equation}
where $r=1$ for TN  condition and $r=2$ for WABE condition with some function $g_0=\order{1}$.
  
\end{proposition}

As a remark,  we look at the order of the truncation error of the divergence damping term to see why the choice of a large coefficient, $\alpha=C_d/h^2$, does not affect the order of the truncation error of the whole scheme. The damping term appears in the pressure equation as well as in the TN boundary condition.
In the pressure equation, the expansion of  the damping term leads up to
$$
\alpha ik  M\uhat_j+ \alpha D_0\vhat_j = \alpha \left[\left(ik \uhat_j+\pd{\vhat_j}{y}\right)+\frac{1}{6}\left(ik\pdn{\uhat}{y}{2}+\pdn{\vhat}{y}{3}\right)+\order{h^4}\right],~~\forall j>0.
$$
The continuity condition implies
$$
ik \uhat_j+\pd{\vhat_j}{y}=ik\pdn{\uhat_j}{y}{2}+\pdn{\vhat_j}{y}{3}=0.
$$
So the divergence damping term with the choice of  $\alpha=C_d/h^2$ contributes an $\order{h^2}$ error that is in line with the accuracy of the other difference operators in the scheme. Similarly, for the divergence damping term in the TN boundary condition, we have
$$
\alpha ik \left(\frac{1}{3}\uhat_0+\frac{1}{6}\uhat_1\right)h+\alpha\frac{\vhat_1-\vhat_0}{2}= \alpha \frac{h}{2}\left[ ik \left(\frac{2}{3}\uhat_0+\frac{1}{3}\uhat_1\right)+\frac{\vhat_1-\vhat_0}{h}\right]= \alpha \frac{h}{2}\left[ik\uhat_{\frac{1}{3}}+\pd{\vhat_{\frac{1}{3}}}{y}+\order{h^2}\right].
$$
Thus, with the continuity condition, we see that the error  contributed by the divergence damping is $\order{h}$ that is also consistent  with the accuracy of the  discrete TN boundary condition.

\subsection{Stability analysis}
It suffices to consider the homogeneous  version of the problem  \eqref{eq:errorEqnLumped} \& \eqref{eq:errorBCLumped} when analyzing the stability of the scheme.
For simplicity, we do not discretize in time and   analyze the stability properties of the semi-discrete problem directly using  Laplace transformation method and normal-mode analysis.  It is mentioned in \cite{ICNS,INSDIV,GustafssonKreissOliger95} that any dissipative time discretization can be used and the resulting fully discrete problem will be stable provided the semi-discrete problem is stable.

As is pointed out in  \cite{GustafssonKreissOliger95}, there are several possible stability definitions.  Here we show the semi-discrete problem is stable in the sense of Godunov-Ryabenkii condition;
that is, we demonstrate its stability by showing that there is no eigenvalue $s$ with  $\Re(s)>0$ for a related eigenvalue problem, which  is obtained by Laplace transforming the semi-discrete problem \eqref{eq:errorEqnLumped} in time with  $s$ denoting the dual variable.

To be specific,  after Laplace transforming the homogeneous version of \eqref{eq:errorEqnLumped} \& \eqref{eq:errorBCLumped}, we obtain the eigenvalue problem,
\begin{align}\label{eq:eigenValueProblem}
    \text{for}~j>0:~
  \begin{cases}
 s {\utilde}_j = - ik \ptilde_j-\nu k^2\utilde_j+\nu D_+D_-\utilde_j,\\
  s {\vtilde}_j = -D_0\ptilde_j- \nu k^2\vtilde_j+\nu D_+D_-\vtilde_j \\
     -k^2\ptilde_j+D_+D_-\ptilde_j=  \alpha ik \utilde_j+ \alpha D_0\vtilde_j ,
    \end{cases}
\end{align}
with the boundary conditions, 
\begin{equation}\label{eq:eigenValueProblemBC}
  \begin{cases}
    \utilde_0=\vtilde_0=0,\\
     D_+\ptilde_0+\nu ikD_+\utilde_0=0,
    \end{cases}
\end{equation}
and the regularity condition,
$$
||\utilde||_h<\infty,~||\vtilde||_h<\infty,~||\ptilde||_h<\infty.
$$
Here ($\utilde_j$, $\vtilde_j$, $\ptilde_j$)  denote the transformed solutions and   $||\cdot||_h$ represents the discrete $L_2$ norm defined on the grid $\mathcal{G}$.

Note that if $k=0$, we  have a non-trivial solution to the eigenvalue problem \eqref{eq:eigenValueProblem}  \&  \eqref{eq:eigenValueProblemBC}  , i.e.,
$$
\begin{pmatrix}
  \utilde_j\\
  \vtilde_j\\
  \ptilde_j
\end{pmatrix}
=
\begin{pmatrix}
  0\\
  0\\
  \ptilde_0
\end{pmatrix} ~\text{for any constant}   ~\ptilde_0.
$$
This solution should be excluded due to  the regularity condition $||\ptilde||_h<\infty$.
This case corresponds to the undetermined constant in the pressure which is regularized in our algorithm by enforcing the mean of the pressure to be zero as is described in \eqref{eq:discretePPEwithRegularization}.
So we proceed the stability analysis  assuming $k\neq0$.

\subsubsection{Without divergence damping ($\alpha=0$)} First, let us consider  $\alpha=0$. In this case, the pressure is decoupled from the velocity equations, and the  solution of the eigenvalue problem can be  found explicitly,  {which is described  in the following lemma.
  
\begin{lemma}\label{lemma:eigenSolutionNodd}
  If $\alpha=0$, the solution of the eigenvalue problem  \eqref{eq:eigenValueProblem} \&  \eqref{eq:eigenValueProblemBC} is found to be
  \begin{equation}\label{eq:noddAnalysisSolutions}
  \begin{cases}
  \medskip\displaystyle\utilde_j=-\frac{ ik C_p}{s}\left(e^{-\xi y_j}- e^{-\gamma y_j}\right),\\
  \medskip \displaystyle \vtilde_j=\frac{1}{hs}\sinh(\xi h)C_p\left(e^{-\xi y_j}- e^{-\gamma y_j}\right),\\
  \medskip  \displaystyle  \ptilde_j=C_p e^{-\xi y_j},
    \end{cases}
  \end{equation}
  where $\xi$ and $\gamma$ satisfy
  \begin{align*}
    \frac{4}{h^2}\sinh^2\left(\frac{\xi h}{2}\right)=k^2, & \quad   \xi>0, \\
    \frac{4}{h^2}\sinh^2\left(\frac{\gamma h}{2}\right)=\frac{s}{\nu}+ k^2, &\quad \Re(\gamma)>0.
  \end{align*}
  The remaining coefficient $C_p$ will be determined by the pressure boundary condition in \eqref{eq:eigenValueProblemBC}, which implies
\begin{equation}\label{eq:CpEqn}
C_p\frac{1}{h}\left[ \left(e^{-\xi h}-1\right)+\nu  \frac{ k^2 }{s}\left(e^{-\xi h}- e^{-\gamma h}\right)\right]=0.
\end{equation}
  \end{lemma}

If we let
$$
q_1(s)=\frac{1}{s}\left(e^{-\xi h}- e^{-\gamma h}\right)~~\text{and}~~q(s)=\left(e^{-\xi h}-1\right)+\nu k^2q_1(s),
$$
then  \eqref{eq:CpEqn} can be written as $C_pq(s)=0$.  For $\Re(s)>0,$ we have the following lemmas  that imply $q(s)\neq 0$; therefore, we conclude that $C_p=0$ and  the solution given in \eqref{eq:noddAnalysisSolutions} is trivial.

\begin{lemma}\label{lemma:noddLemma_q1} 
 If $q_1(s)$ is real, then $s$ is real.
\end{lemma}

\begin{lemma}\label{lemma:noddLemma_q} 
 $q(s)\neq0$ for  $\Re(s)>0$
\end{lemma}

For a concise presentation of the main results,  the more technical proofs for these Lemmas are shown  later  in   \ref{sec:appendix},  and  here  some examples of $q_1(s)$  and $q(s)$ for $s>0$ and various wavenumbers are plotted in Figure~\ref{fig:qsExamples} to help demonstrate the analytical results.  From the plots, we see that both $q_1(s)$ and $q(s)$ are decreasing and no root of $q(s)$ is observed for $s>0$, which are consistent with the analysis.

    {
\newcommand{\figWidth}{7cm}
\newcommand{\trimfig}[2]{\trimw{#1}{#2}{0.}{0.}{0.}{0.0}}
\begin{figure}[h]
\begin{center}
\begin{tikzpicture}[scale=1]
  \useasboundingbox (0.0,0.0) rectangle (14.,6);  

  \draw(-0.5,0.) node[anchor=south west,xshift=0pt,yshift=0pt] {\trimfig{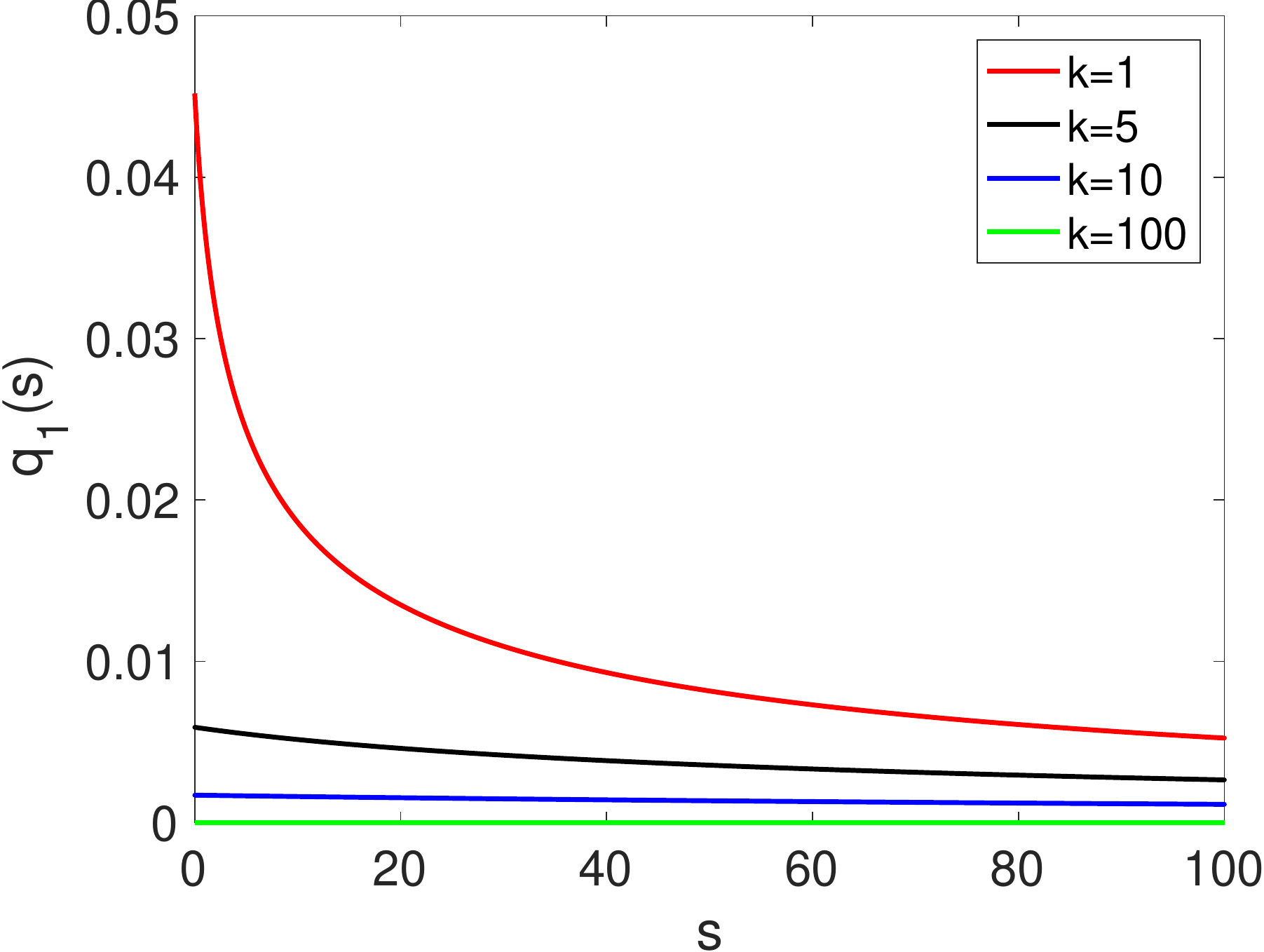}{\figWidth}};
\draw(6.8,0.) node[anchor=south west,xshift=0pt,yshift=0pt] {\trimfig{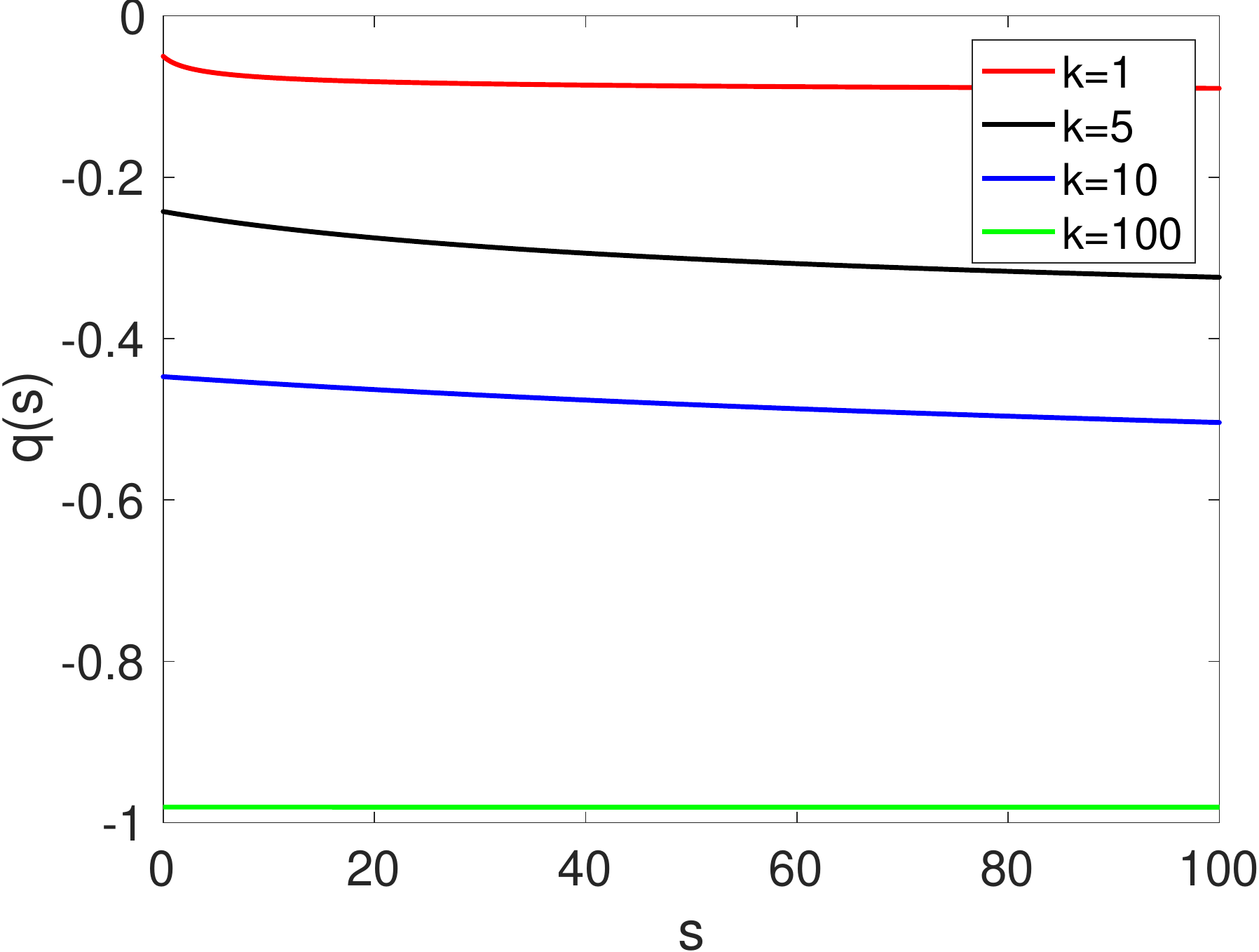}{\figWidth}};
%
\end{tikzpicture}

\end{center}
    \caption{Plots of $q_1(s)$ (left image) and $q(s)$  (right image) with $h=0.1$ and $\nu=1$ for various wavenumbers.  }\label{fig:qsExamples}
\end{figure}
    }

In conclusion, the eigenvalue problem given in  \eqref{eq:eigenValueProblem}  \&  \eqref{eq:eigenValueProblemBC}  has no eigenvalue $s$ with  $\Re(s)>0$, and hence the following proposition concerning the stability of the scheme   holds.
\begin{proposition}\label{pp:noDDStable}
If the divergence damping coefficient  $\alpha=0$, the semi-discrete problem \eqref{eq:errorEqnLumped} with boundary condition \eqref{eq:errorBCLumped} is stable in the sense of Godunov-Ryabenkii condition.
\end{proposition}

}

\subsubsection{With divergence damping ($\alpha\neq 0$)}
Now, we consider $\alpha\neq 0$.  In this case, the pressure and velocity components are coupled. To solve the difference equations, we make the ansatz
$$
\begin{pmatrix}
  \utilde_j\\
  \vtilde_j\\
  \ptilde_j
\end{pmatrix}
= \lambda^j
\begin{pmatrix}
  \utilde_0\\
  \vtilde_0\\
  \ptilde_0
\end{pmatrix}.
$$
Since we want the solutions to be bounded at $y=\infty$, we look for solutions with $|\lambda|<1$.

Inserting the ansatz into the eigenvalue problem \eqref{eq:eigenValueProblem}, we have 
\begin{equation}\label{eq:ansatzEqns}
\begin{pmatrix}
  h^2(s+\nu k^2)-\nu d_2(\lambda) & 0 & ikh^2\\
  0 &  h^2(s+\nu k^2)-\nu d_2(\lambda) & hd_1(\lambda)\\
  \alpha ikh^2 & \alpha h d_1(\lambda) & k^2h^2-d_2(\lambda)
\end{pmatrix}
\begin{pmatrix}
    \utilde_0\\
  \vtilde_0\\
  \ptilde_0
\end{pmatrix}
= \mathbf{0},
\end{equation}
where the following observations have been utilized,
\begin{align*}
  h D_0 \lambda^j = d_1(\lambda) \lambda^j, & ~~\text{where} ~~d_1(\lambda)=\frac{\lambda-\lambda^{-1}}{2},\\
  h^2D_+D_- \lambda^j =d_2(\lambda)\lambda^j,  & ~~\text{where} ~~d_2(\lambda)=\lambda-2+\lambda^{-1}.
\end{align*}
The characteristic equation of \eqref{eq:ansatzEqns}  is 
\begin{equation}\label{eq:deterinantCondition}
\left[h^2(s+\nu k^2)-\nu d_2\right]\left\{\left[h^2(s +\nu k^2+\alpha)-\nu d_2\right](k^2h^2-d_2)+\alpha h^2 (d_2-d_1^2) \right\}=0.
\end{equation}
Noticing  that  $d_1^2=d_2(d_2+4)/4$, the characteristic equation is in fact a cubic equation for $d_2$; namely,
$$
\left[h^2(s+\nu k^2)-\nu d_2\right]
\left\{
(\nu - \alpha h^2/4)d_2^2 - h^2(s+2\nu k^2 + \alpha )d_2 + h^4k^2(s+\nu k^2 + \alpha )
 \right\}=0.
 $$
 So there are three roots for $d_2$:
 \def\aa{(\nu - \alpha h^2/4)}
 \def\nbb{ h^2(s+2\nu k^2 + \alpha )} 
 \def\cc{h^4k^2(s+\nu k^2 + \alpha )}
 \def\dd{[(\alpha+s)^2+(s+\alpha+\nu k^2) \alpha h^2 k^2 ] h^4}
 \begin{equation}\label{eq:d2Solutions}
   \begin{cases}
  \medskip\displaystyle d^{(1)}_2= \frac{h^2(s+\nu k^2)}{\nu},\\
  \medskip\displaystyle d^{(2)}_2=\frac{\nbb + \sqrt{\dd}}{2\aa},\\
  \medskip\displaystyle d^{(3)}_2=\frac{\nbb - \sqrt{\dd}}{2\aa}.\\
     \end{cases}
 \end{equation}
 Note that the three roots are distinct if $\Re(s)>0$.
 For each $d_2^{(n)}$, we have  an equation for $\lambda$,
 \begin{equation}\label{eq:lambdaEquation}
 \lambda^2-\left(2+d_2^{(n)}\right)\lambda+1=0, ~~n=1,2,3.
 \end{equation}
The $\lambda$ equation has two reciprocal roots, and   the  root with magnitude less than one  is denoted as    $\lambda_{(n)}$.
The corresponding solutions for $(\utilde_0,\vtilde_0, \ptilde_0)^T$ are
 $$
 \begin{pmatrix}
\medskip\displaystyle   \utilde_0^{(1)}\\
\medskip\displaystyle   \vtilde_0^{(1)}\\
\medskip\displaystyle   \ptilde_0^{(1)}
\end{pmatrix}
=
\begin{pmatrix}
\medskip \displaystyle -\frac{d^{(1)}_1}{h}\\
\medskip\displaystyle   ik\\
\medskip\displaystyle   0
\end{pmatrix},
\quad
\begin{pmatrix}
\medskip\displaystyle   \utilde_0^{(2)}\\
\medskip\displaystyle   \vtilde_0^{(2)}\\
\medskip\displaystyle   \ptilde_0^{(2)}
\end{pmatrix}
=
\begin{pmatrix}
\medskip\displaystyle   ik\\
\medskip\displaystyle \frac{d_1^{(2)}}{h}\\
\medskip\displaystyle  -(s+\nu k^2)+\nu \frac{d_2^{(2)}}{h^2}
\end{pmatrix},
\quad
\begin{pmatrix}
\medskip\displaystyle   \utilde_0^{(3)}\\
\medskip\displaystyle   \vtilde_0^{(3)}\\
\medskip\displaystyle  \ptilde_0^{(3)}
\end{pmatrix}
=
\begin{pmatrix}
\medskip\displaystyle   \frac{ik}{\alpha}\\
\medskip\displaystyle   \frac{d_1^{(3)}}{(h\alpha)}\\
\medskip\displaystyle   [-(s+\nu k^2)+\nu \frac{d_2^{(3)}}{h^2}]\frac{1}{\alpha}
\end{pmatrix},
$$
where  $d_1^{(n)}=\left(\lambda_{(n)}-\lambda^{-1}_{(n)}\right)/2$.

Therefore, the general solution to   \eqref{eq:eigenValueProblem}  is 
\begin{equation}\label{eq:generalSolution}
\begin{pmatrix}
  \utilde_j\\
  \vtilde_j\\
  \ptilde_j
\end{pmatrix}
=
\sum_{n=1}^3\sigma_n\lambda_{(n)}^j
\begin{pmatrix}
  \utilde_0^{(n)}\\
  \vtilde_0^{(n)}\\
  \ptilde_0^{(n)}
\end{pmatrix}.
\end{equation}
After applying the boundary conditions \eqref{eq:eigenValueProblemBC}, we have a system of equations for $\mathbf{\sigma}=(\sigma_1,\sigma_2,\sigma_3)^T$,
$$
Z\mathbf{\sigma}=\mathbf{0},
$$
where
\begin{equation}\label{eq:matrixZ}
Z=
\begin{pmatrix}
 \medskip \displaystyle -\frac{d^{(1)}_1}{h}&   ik & \displaystyle \frac{ik}{\alpha} \\
 \medskip  ik            & \displaystyle  \frac{d_1^{(2)}}{h} &\displaystyle \frac{d_1^{(3)}}{h\alpha} \\
 \medskip \displaystyle -\nu ik  \frac{d^{(1)}_1(\lambda_{(1)}-1)}{h^2} &\displaystyle  \left[-(s+2\nu k^2)+\nu \frac{d_2^{(2)}}{h^2}\right]\frac{\lambda_{(2)}-1}{h} & \displaystyle   \left[-(s+2\nu k^2)+\nu \frac{d_2^{(3)}}{h^2}\right]\frac{\lambda_{(3)}-1}{h \alpha}
\end{pmatrix}.
\end{equation}
The numerical scheme given in  \eqref{eq:errorEqnLumped} \& \eqref{eq:errorBCLumped}  is stable  in the Godunov-Ryabenkii sense if we can show the determinant condition $\det(Z) \neq 0$ holds for $\Re(s)>0$. Note that the determinant condition implies  $\mathbf{\sigma}=\mathbf{0}$ and  the  solution \eqref{eq:generalSolution} is trivial; thus,  there is no eigenvalue $s$ with $\Re(s)>0$.

Here we provide some evidence that the determinant condition is  not violated by  plotting  zero contours of the real and imaginary parts of  $\det(Z(s)), \forall s\in\mathbb{C}$ for some wavenumbers ($k=1,5,10,100$)   in Figure~\ref{fig:zeroContoursDetZ}. Note that the intersections  of $\Re(\det(Z))=0$ and $\Im(\det(Z))=0$ indicate the values of $s$ that makes  $\det(Z)= 0$. As is shown in Figure~\ref{fig:zeroContoursDetZ},  no intersections are  found  on the right half of the complex plane (i.e.,$\Re(s)>0$) for all the examples considered.

{\bf Remark:} 
 it is hard  to prove the  the determinant condition for the general case.  But following \cite{INSDIV,BCNS}, we are able to establish the so-called ``local stability'' for the scheme. Further, by solving the leading order terms of the error equations explicitly, we can identify that with $\alpha\sim 1/h^2$ the scheme is second-order accurate because the error introduced on the boundary will be damped out by the divergence damping term. However, since TN boundary condition corresponds to a first-order discretization, the boundary-layer errors in the pressure solution is expected for this case.
Details of the local stability and accuracy are discussed in the following subsection.

{
  \newcommand{\figWidth}{7cm}
  \def\xa{14}
  \def\ya{12}
\newcommand{\trimfig}[2]{\trimw{#1}{#2}{0.}{0.}{0.}{0.0}}
\begin{figure}[h]
\begin{center}
\begin{tikzpicture}[scale=1]
  \useasboundingbox (0.0,0.0) rectangle (\xa,\ya);  

  \draw(-0.5,6.) node[anchor=south west,xshift=0pt,yshift=0pt] {\trimfig{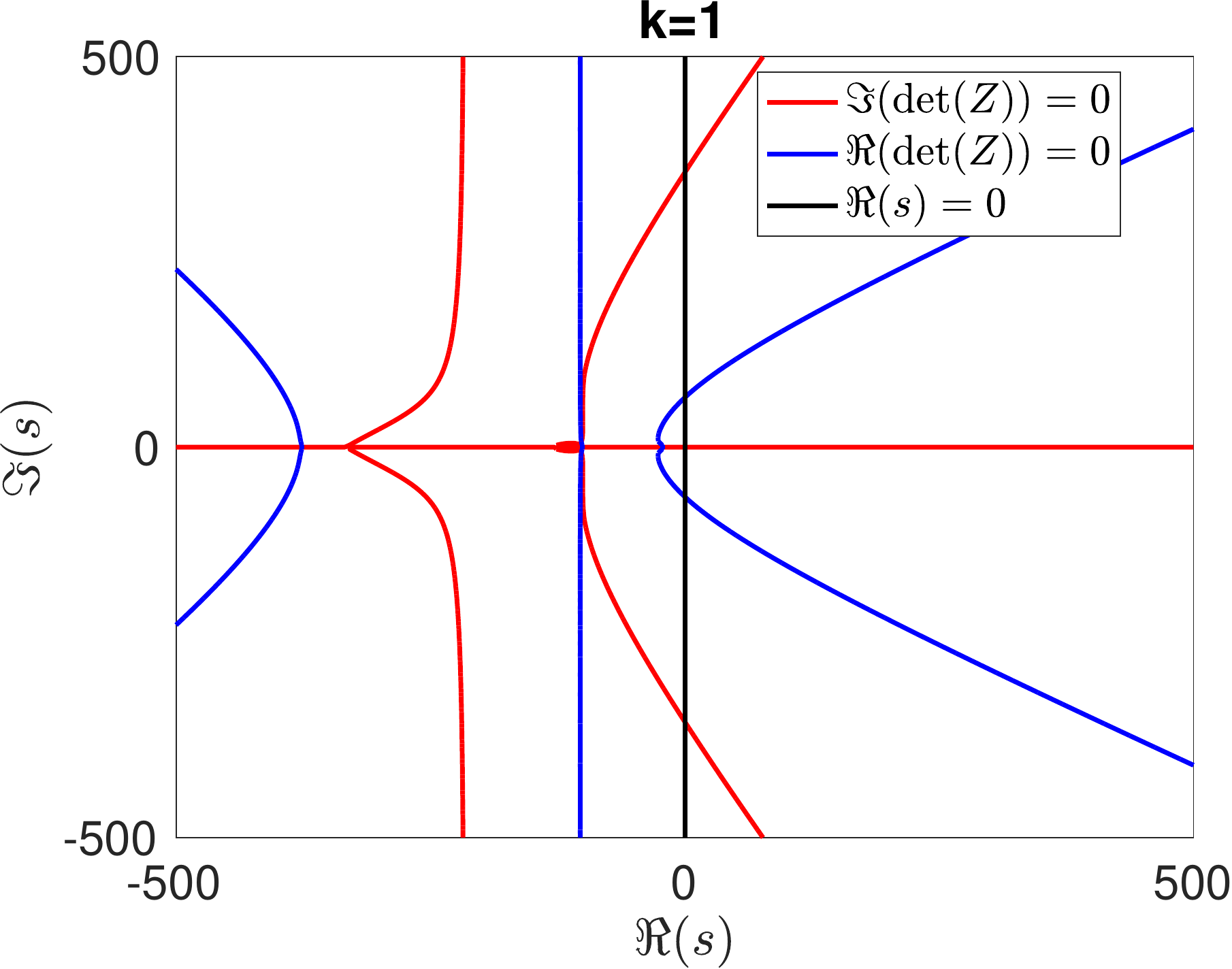}{\figWidth}};
  \draw(6.8,6.) node[anchor=south west,xshift=0pt,yshift=0pt] {\trimfig{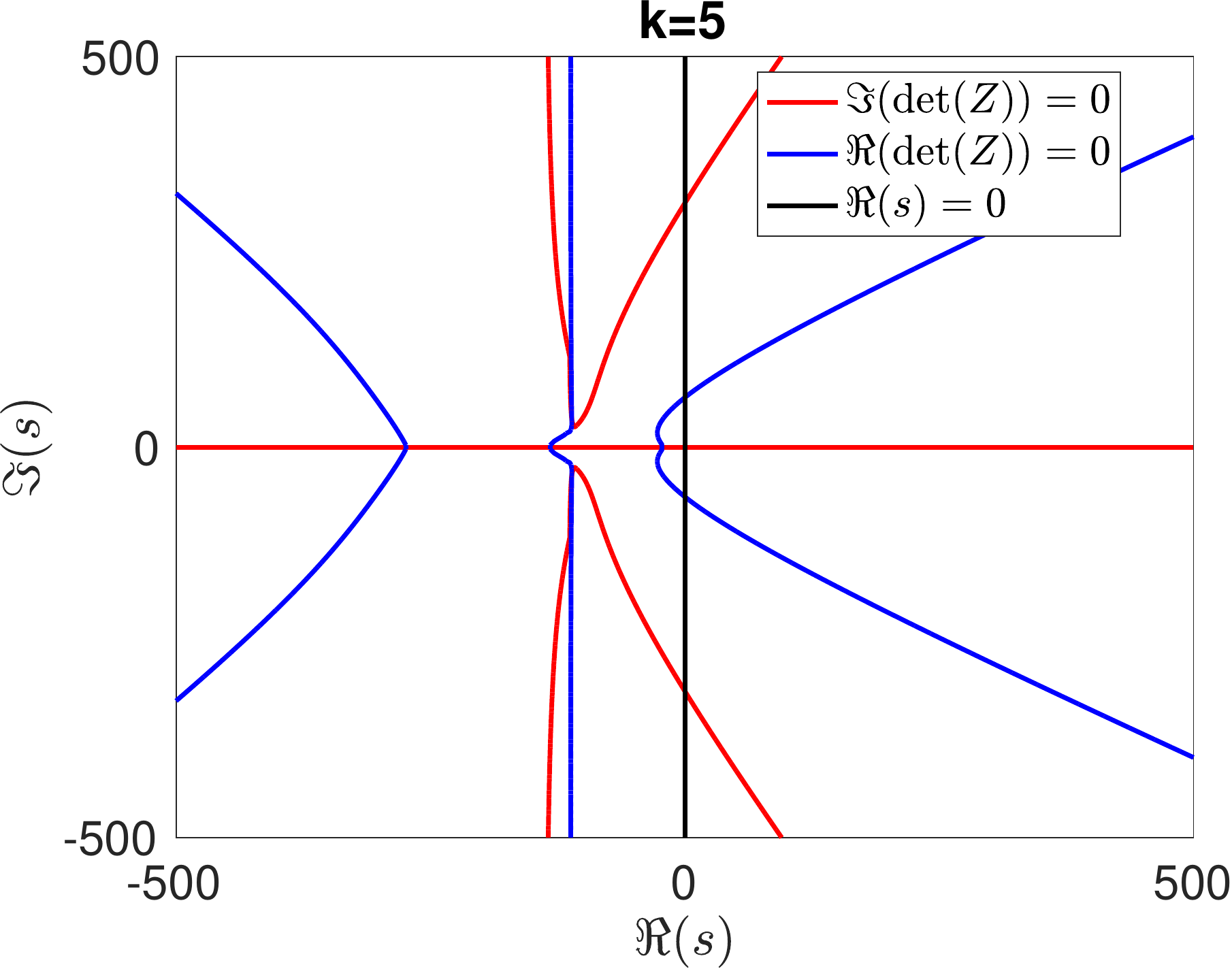}{\figWidth}};
    \draw(-0.5,0.) node[anchor=south west,xshift=0pt,yshift=0pt] {\trimfig{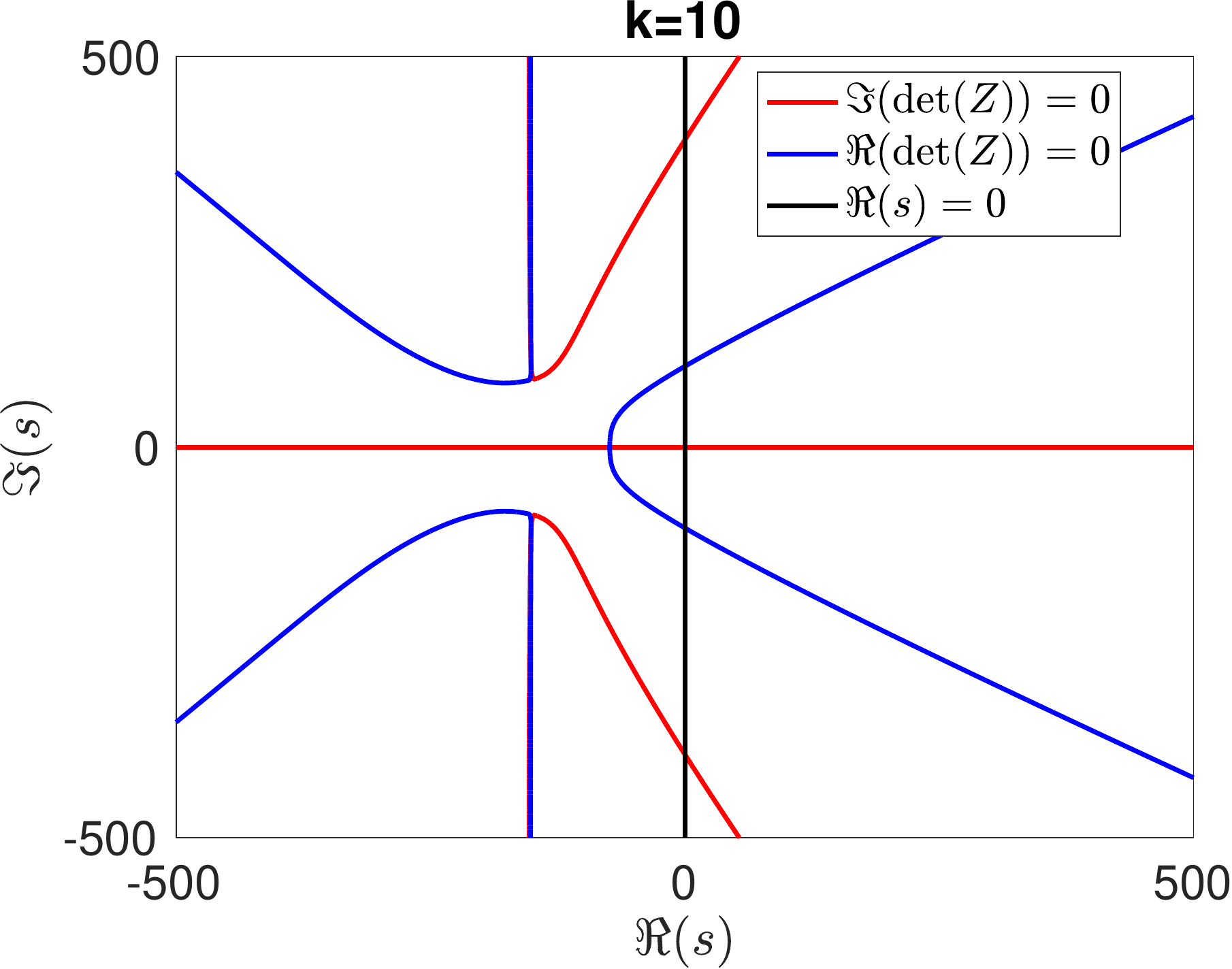}{\figWidth}};
\draw(6.8,0.) node[anchor=south west,xshift=0pt,yshift=0pt] {\trimfig{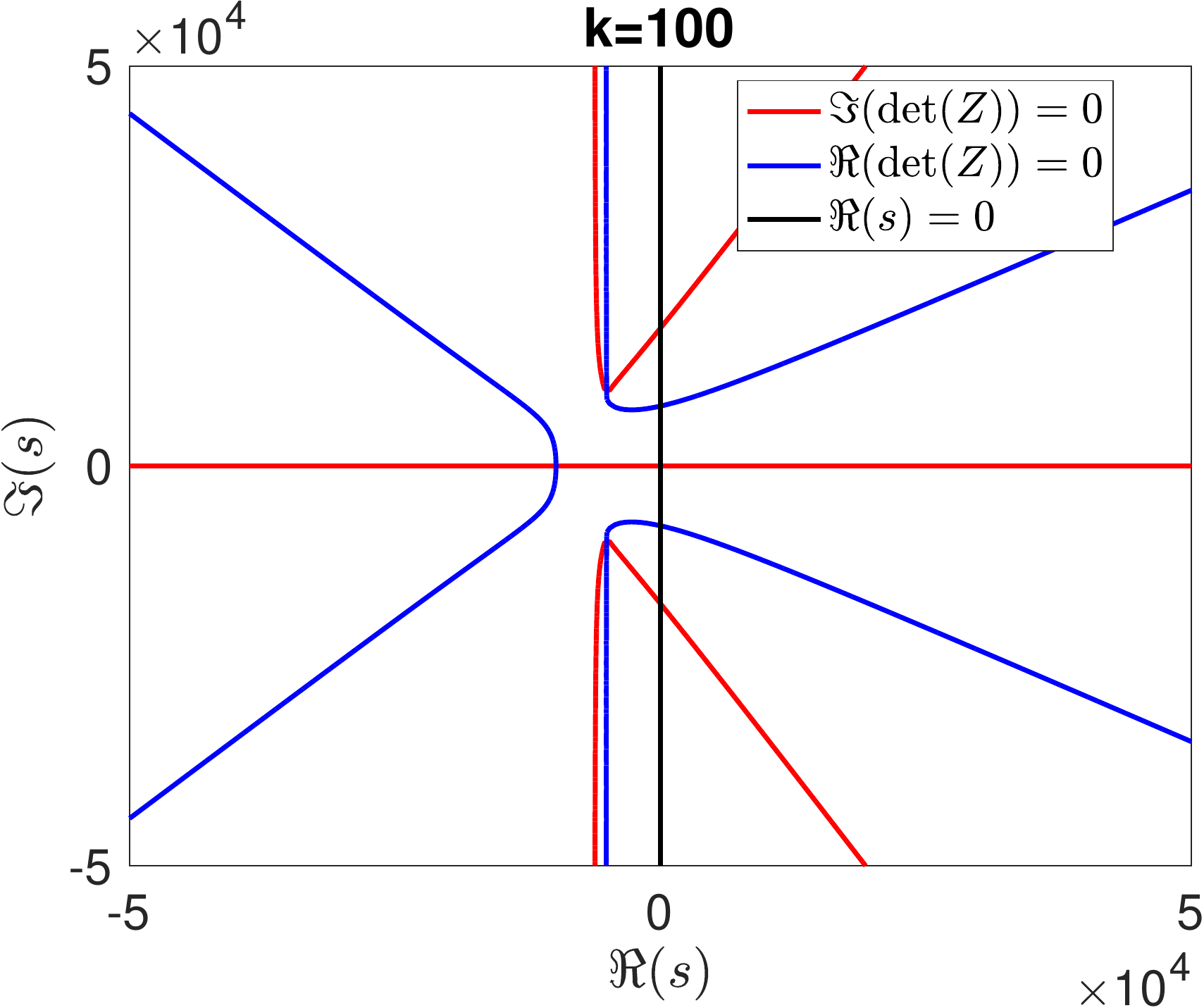}{\figWidth}};
%
\end{tikzpicture}

\end{center}
    \caption{Zero contours of $\Re(\det(Z))$ and $\Im(\det(Z))$  with $h=0.1$ and $\nu=1$ for various wavenumbers. }\label{fig:zeroContoursDetZ}
\end{figure}
}

\subsection{Local Stability and Accuracy}
In this section, we show  stability and accuracy  of the scheme assuming that $h\sqrt{s/\nu+k^2}\ll 1$. Stability results established under this assumption is referred to as local stability in the literature  \cite{INSDIV,BCNS}. A scheme that is locally stable but not stable in the global sense can be quickly  identified in computations since the unstable modes occur at high frequencies.

\begin{proposition}\label{ddLemma1} Assuming that $h\sqrt{s/\nu+k^2}\ll 1$ and  $\alpha=C_d/h_0^2$  for some fixed grid spacing $h_0$,
 there exists a constant $h_c$ such that the determinant condition $\det(Z)\neq 0$ holds for $\forall h<h_c$ for $\Re(s)>0$ and a fixed wavenumber $k$.
\end{proposition}
\begin{proof}
 By solving equation \eqref{eq:lambdaEquation}  with the assumption $\Re(s)>0$, we know the root with magnitude less one is of the following form:   
   $$
 \lambda_{(n)}=\frac{\left(2+d_2^{(n)}\right)-\sqrt{d_2^{(n)}\left(d_2^{(n)}+4\right)}}{2},~~n=1,2,3.
 $$
 From \eqref{eq:d2Solutions}, we have
 $$
     \begin{cases}
  \medskip\displaystyle \lim_{h\rightarrow 0}\frac{d^{(1)}_2}{h^2}= \frac{s+\nu k^2}{\nu},\\
  \medskip\displaystyle  \lim_{h\rightarrow 0}\frac{d^{(2)}_2}{h^2}=\frac{s+\nu k^2+\alpha}{\nu},\\
  \medskip\displaystyle \lim_{h\rightarrow 0}\frac{d^{(3)}_2}{h^2}=k^2.\\
     \end{cases}
 $$
Further, we  can show that 
     $$
    \lim_{h\rightarrow 0} \frac{\lambda_{(n)}-1}{h}=-\lim_{h\rightarrow 0} \sqrt{\frac{d_2^{(n)}}{h^2}}~~\text{and}~~ \lim_{h\rightarrow 0} \frac{d_1^{(n)}}{h}=-\lim_{h\rightarrow 0} \sqrt{\frac{d_2^{(n)}}{h^2}}.
    $$
    Therefore, we have
    \begin{equation}\label{eq:leadingOrderZ}
     \lim_{h\rightarrow 0}Z =\begin{pmatrix}
 \medskip\displaystyle \sqrt{\frac{s+\nu k^2}{\nu}}&   ik &  \displaystyle \frac{ik}{\alpha} \\
 \medskip  ik            &\displaystyle -\sqrt{\frac{s+\nu k^2+\alpha}{\nu},}& \displaystyle -\frac{|k|}{\alpha} \\
 \medskip  - ik (s+\nu k^2)  & \displaystyle  (\nu k^2-\alpha)\sqrt{\frac{s+\nu k^2+\alpha}{\nu}} & \displaystyle \frac{(s+\nu k^2)|k|}{\alpha}
     \end{pmatrix}
     \end{equation}
     and hence
     $$
     \lim_{h\rightarrow 0 }\det(Z)=-\frac{1}{\alpha}(\alpha + s) \left(|k| \sqrt{\frac{\nu k^2 + s}{\nu}} - k^2\right)\sqrt{\frac{\nu k^2 + \alpha + s}{\nu}}.
     $$
     Obviously, if $\Re(s)>0$, we have
     $$
      \lim_{h\rightarrow 0 }\det(Z)\neq 0.
      $$
      Since $\det(Z)$ continuously depends on $h$, there exists a constant $h_c$ such that
      $
      \det(Z)\neq 0,~~\forall h<h_c.
      $
This completes the proof.
\end{proof}

As an example, we plot the zero contours of the real and imaginary parts of  $\lim_{h\rightarrow 0 }\det(Z)$ with wavenumber $k=1$ and $10$ in Figure~\ref{fig:zeroContoursLimDetZ}. As expected, no intersection is observed when $\Re(s)>0$ for both cases.
{
\newcommand{\figWidth}{7cm}
\newcommand{\trimfig}[2]{\trimw{#1}{#2}{0.}{0.}{0.}{0.0}}
\begin{figure}[h]
\begin{center}
\begin{tikzpicture}[scale=1]
  \useasboundingbox (0.0,0.0) rectangle (14.,6);  

  \draw(-0.5,0.) node[anchor=south west,xshift=0pt,yshift=0pt] {\trimfig{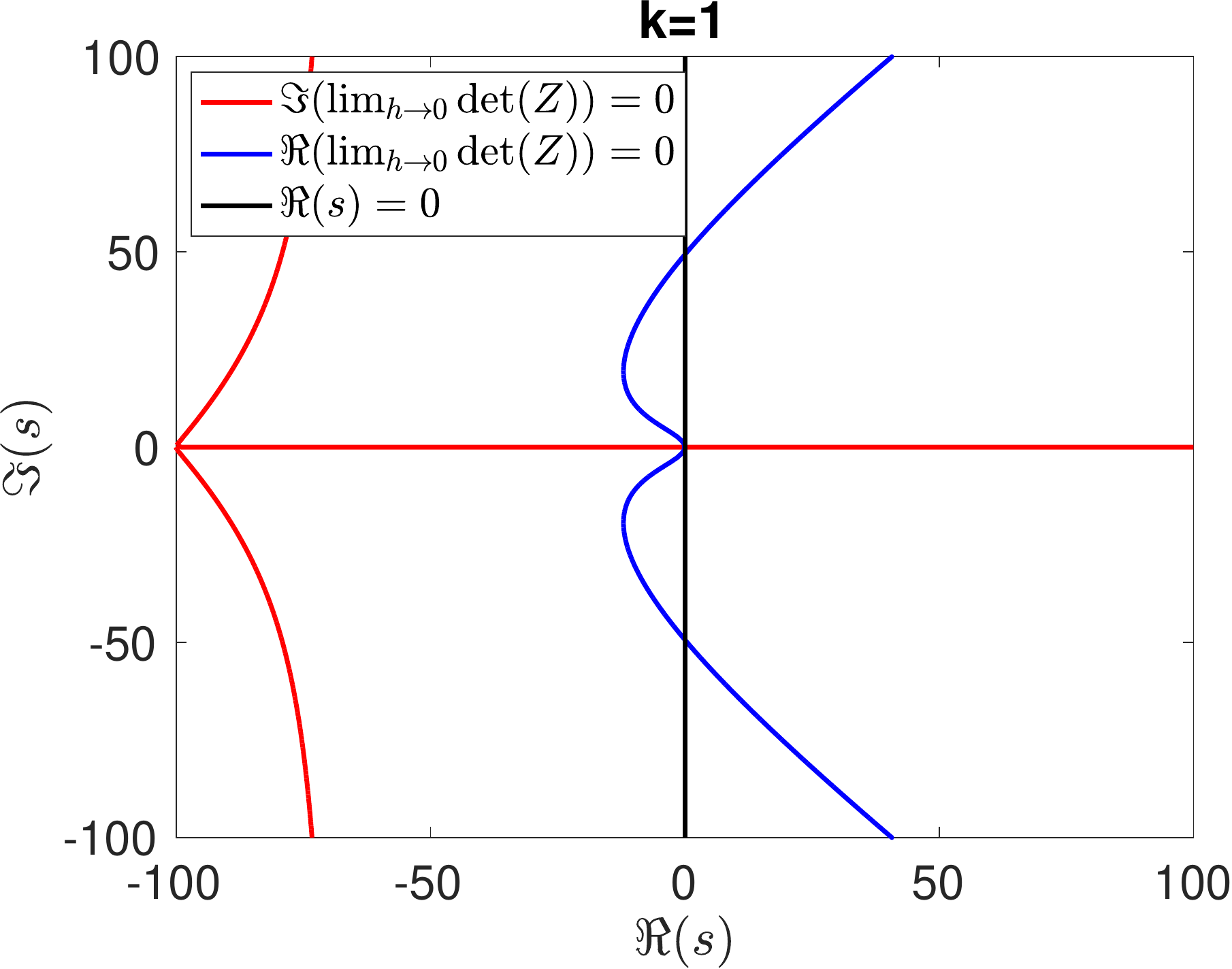}{\figWidth}};
\draw(6.8,0.) node[anchor=south west,xshift=0pt,yshift=0pt] {\trimfig{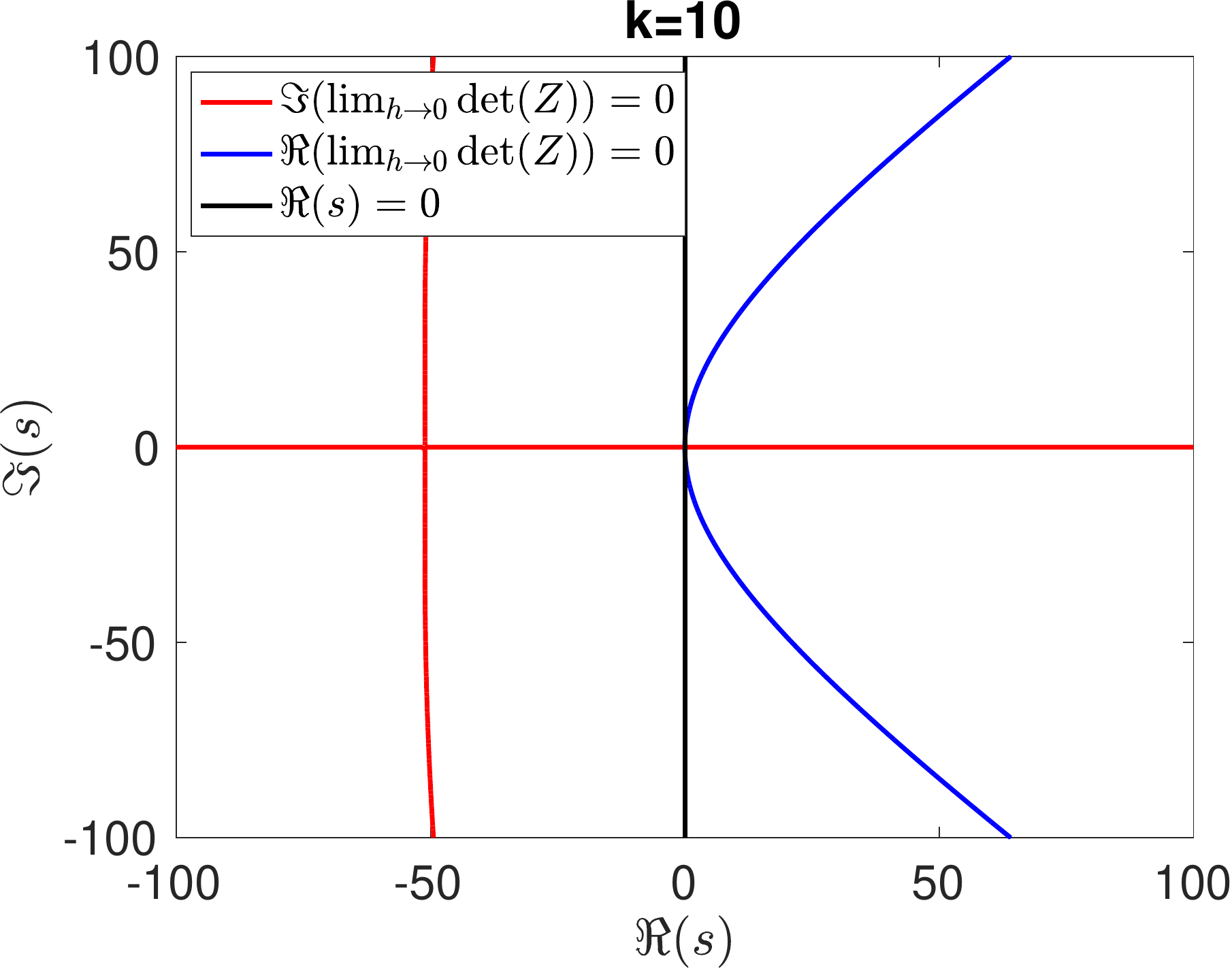}{\figWidth}};
%
\end{tikzpicture}

\end{center}
    \caption{Zero contours of the real and imaginary parts of $\lim_{h\rightarrow0}(\det(Z))$ with $\alpha=100$ and $\nu=1$ for wavenumbers $k=1$ (left figure) and $k=10$ (right figure). } \label{fig:zeroContoursLimDetZ}
\end{figure}
}

With the local stability assumption, we can 
  further show how the TN and WABE boundary conditions affect the accuracy of the numerical results by solving  the leading-order terms of the  error equation \eqref{eq:eigenValueProblem} with the non-homogeneous boundary condition that accounts for the truncation error of the corresponding boundary conditons, i.e.,
\begin{equation}\label{eq:nonhomoBCForError}
     D_+\ptilde_0+\nu ikD_+\utilde_0=h^r\gtilde_0,
\end{equation}
where $r=1$ for TN  condition and $r=2$ for WABE condition; here  $\gtilde_0=\order{1}$. Notice that this boundary condition is the Laplace transformation of \eqref{eq:errorBCLumped}. It has already been  shown that the general solution to  \eqref{eq:eigenValueProblem} is  \eqref{eq:generalSolution}. After applying the no-slip boundary conditions $\utilde=\vtilde=0$ and the boundary condition \eqref{eq:nonhomoBCForError}, we have
$$
Z\mathbf{\sigma} =
\begin{pmatrix}
  0\\
  0\\
  h^r\gtilde_0
  \end{pmatrix},
$$
where $Z$ is the same as \eqref{eq:matrixZ} and its leading order is given by equation~\eqref{eq:leadingOrderZ}. The solutions to  the leading order equations
are 
\begin{align*}
\medskip  &\sigma_1\sim   \frac{i\gtilde_0 h^r k \left(|k| - \sqrt{(\nu  k^2 + \alpha + s)/\nu}\right)}{(\alpha + s) (|k| \sqrt{(\nu k^2 + s)/\nu} - k^2)\sqrt{(\nu k^2 + \alpha + s)/\nu}},\\
\medskip  & \sigma_2\sim -\frac{\gtilde_0 h^r}{(\alpha + s)\sqrt{(\nu k^2 + \alpha + s)/\nu}},\\
\medskip  & \sigma_3\sim  \frac{\alpha \gtilde_0 h^r \left(\sqrt{(\nu k^2 + s)/\nu} \sqrt{(\nu k^2 + \alpha + s)/\nu} - k^2\right)}{(\alpha + s) (|k| \sqrt{(\nu k^2 + s)/\nu} - k^2) \sqrt{(\nu k^2 + \alpha + s)/\nu}}.
\end{align*}
Thus, we  see that if $g_0\neq 0$, the error introduced by the boundary condition can affect the interior. However, with large $\alpha$, the boundary error is rapidly damped out producing a boundary layer of order $r$ ($r=1$ for TN and $r=2$ for WABE).

{\bf Remark:} 
the error estimate of the IBVP \eqref{eq:errorEqnLumped} \&  \eqref{eq:errorBCLumped} can be obtained in two stages. First, we  obtain estimates for a pure initial-value problem on a periodic domain satisfying the forcing terms (i.e., $h^2F_u$, $h^2F_u$, $h^2F_p$); this problem is the same as the 2nd-order finite-difference scheme, and it is shown in  \cite{INSDIV} that the error estimate is $\order{h^2}$. 
And then,  after subtracting the solutions of the pure initial-value problem  from the IBVP, we have a new IBVP with zero forcing on the interior equations and inhomogeneous boundary conditions.  We have already shown that the boundary errors are quickly damped out  producing a boundary layer of order $r$. Therefore, after using the Parseval's relation, we know that the scheme is $\order{h^2}$ in $L_2$ norm, but with the existence of  numerical boundary-layer errors that is $\order{h^r}$ ($r=1$ for TN and $r=2$ for WABE). The results are supported by careful numerical mesh refinement studies shown below.

\section{Numerical results}\label{sec:results}
We now present the results for a series of simulations chosen to demonstrate the properties of our numerical approach. We  first consider the INS equations on a sequence of refined unit square meshes to study the accuracy of the scheme; cases  with and without the divergence damping are considered subject to both the  TN and WABE boundary conditions. Some  benchmark problems are also considered to further illustrate the numerical properties of our scheme and to compare with existing results. Finally, as a demonstration that our scheme can be easily extended to work with higher-order elements, we solve the classical flow-past-a-cylinder problem using  $\Pe_n$ finite elements with $n\geq 1$.

  {\bf Remark}: for all  the  test problems with  known exact solutions,  errors of the numerical solutions are measured using  both $L_\infty$ and $L_2$ norms. To be specific, given an exact solution $v$ and its  FEM approximation $v_h$ in the finite space $V_h=\mathrm{span}\{\varphi_1,\varphi_2,\dots,\varphi_n\}$,  we define  the error function as
  $$
  E(v) = |v_h-v_e|,
  $$
where $v_e$ is the projection of the exact solution $v$ onto the  finite space $V_h$, i.e.,
$
v_e = \sum_{i=1}^{n} v(\xv_i)\varphi_i
$. Here  $\xv_i$ is the coordinates of the corresponding degree of freedom. In $V_h$, the $L_\infty$ and $L_2$ norms of the error function is given by 
$$
||E(v) || _\infty = \max(E(v) )~~\text{and}~~||E(v)||_2=\left(\int_\Omega E(v)^2d\xv\right)^{1/2}.
$$
A numerical quadrature with sufficient order of accuracy is used to compute the integral; for example, for $\Pe_1$ elements, a third order accurate quadrature rule is used.

\subsection{Manufactured Solutions}
To numerically investigate  the accuracy and stability  of the finite element scheme, we perform careful mesh refinement study using the method of  manufactured solutions \cite{Roache1998}. Exact solutions of the INS  equations can be constructed by adding forcing functions to the governing equations. The forcing is specified so that a chosen function becomes an exact solution to the forced equation. Here
we use the following trigonometric functions as the exact solutions for our convergence tests,
\begin{align*}
 u_e = a\sin(f_x\pi x)\sin(f_y\pi y)\cos(f_t\pi t),\\
v_e = a\cos(f_x\pi x)\cos(f_y\pi y)\cos(f_t\pi t),\\  
p_e = a\sin(f_x\pi x)\cos(f_y\pi y)\cos(f_t\pi t).
\end{align*}
Note that the exact solutions are chosen to be divergence free.
Parameters for the exact solutions are specified as 
$
a=0.5,
f_x=2,
f_y=2,
$ and
$
f_t=2.
$

\begin{figure}[h]
  \begin{center}
    \includegraphics[width=4cm]{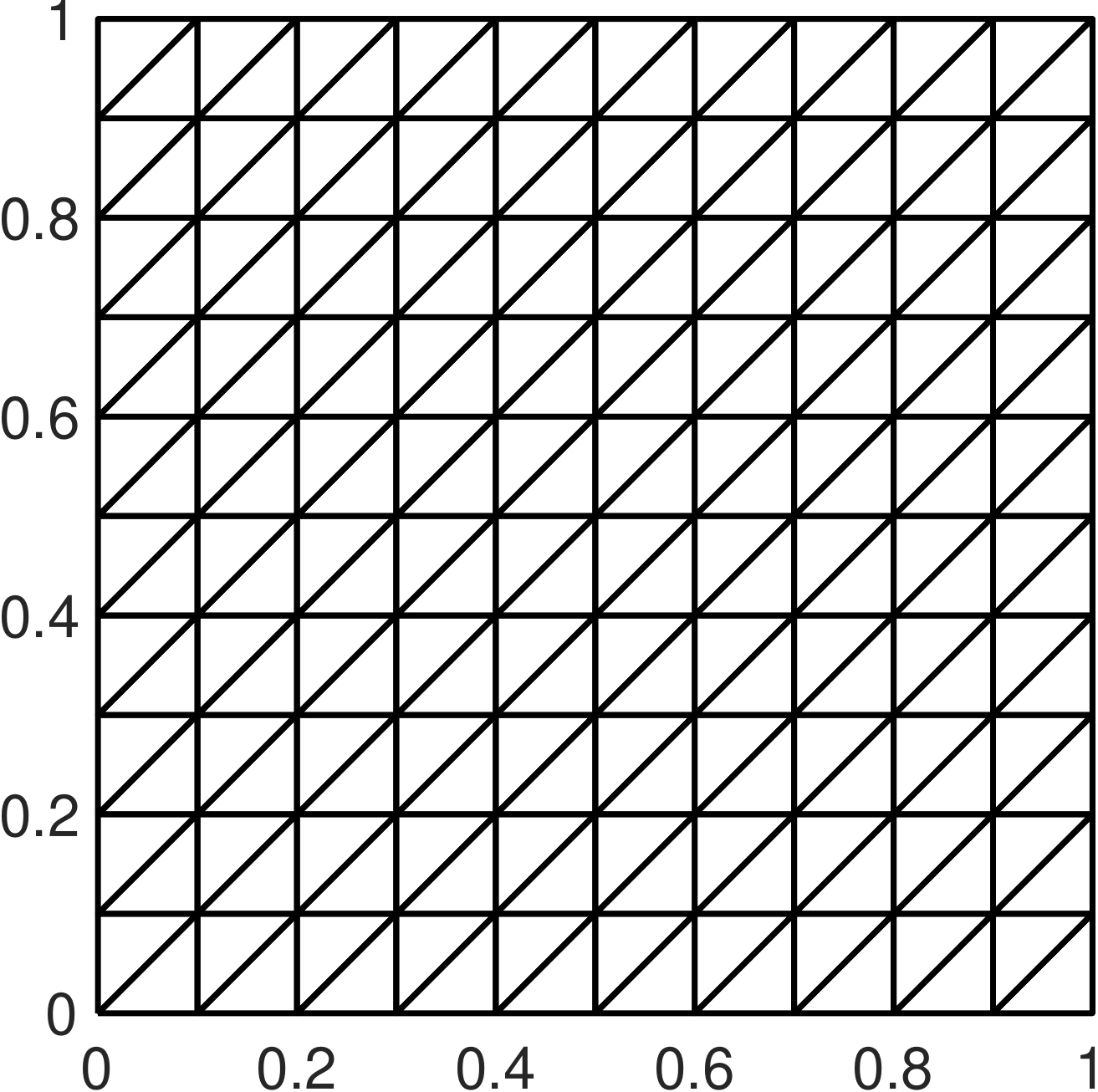}
    \caption{A unit square with uniform triangular mesh. The mesh size is $h=0.1$.}\label{fig:squareMesh}
  \end{center}
\end{figure}

For simplicity,  the test problems are solved on a unit  square  domain ($\Omega=[0,1]\times[0,1]$) with uniform triangular meshes using grid spacings $h=1/(10j), j=1,2,...$. The mesh with $h=1/10$ is  shown in  Figure~\ref{fig:squareMesh}.
The INS equations are discretized using $\Pe_1$ finite element and numerically solved using the   algorithm described in Section \ref{sec:algorithm}. We note that the discrete system obtained using $\Pe_1$ element on this particular mesh is equivalent to the 2nd-order centered finite difference scheme after lumping the mass matrix. Therefore, the numerical tests considered here are directly related  to the normal-mode analysis conducted in Section~\ref{sec:analysis}.
In order to respectively study the effects of the divergence damping and the pressure boundary conditions, we consider the following four cases: (i) $\alpha=0$ with TN boundary condition; (ii)  $\alpha=1/h^2$ with TN boundary condition; (iii) $\alpha=0$ with WABE boundary condition; and (iv) $\alpha=1/h^2$ with WABE boundary condition.

 \subsubsection{Periodic in $x$ direction}
 
We begin the convergence  study  by assuming periodicity in $x$ direction and  no-slip boundary conditions on the other boundaries (i.e., $y=0$ and $y=1$).  This test is designed to match the assumptions of the analysis discussed in Section~\ref{sec:analysis}. 
Results  for the cases (i) TN boundary condition without divergence damping, (ii)  TN boundary condition with divergence damping and (iv) WABE boundary condition with divergence damping are collected in  in Figure~\ref{fig:errorPlotsPD}. The errors  for the velocity component $v$ and the pressure $p$ are plotted with  the first row of  images  for case (i), the second row for case (ii), and the third row for case (iv). The solutions plotted here are obtained using the uniform square mesh with grid spacing $h=1/160$. Here we observe that the errors in all cases  are well behaved in that the magnitudes are small and they are smooth throughout the domain except near the boundary layers for cases (ii) and (iv). We also observe that the interior accuracy is improved by including divergence damping, and the boundary accuracy is further improved by using the WABE condition, which is consistent with our analytical results.

{
\newcommand{\figWidth}{7cm}
\def\xa{14}
\def\ya{19.5}
\newcommand{\trimfig}[2]{\trimw{#1}{#2}{0.}{0.}{0.}{0.0}}
\begin{figure}[b]
\begin{center}
\begin{tikzpicture}[scale=1]
  \useasboundingbox (0.0,0.0) rectangle (\xa,\ya);  
\draw(-0.5,13.0) node[anchor=south west,xshift=0pt,yshift=0pt] {\trimfig{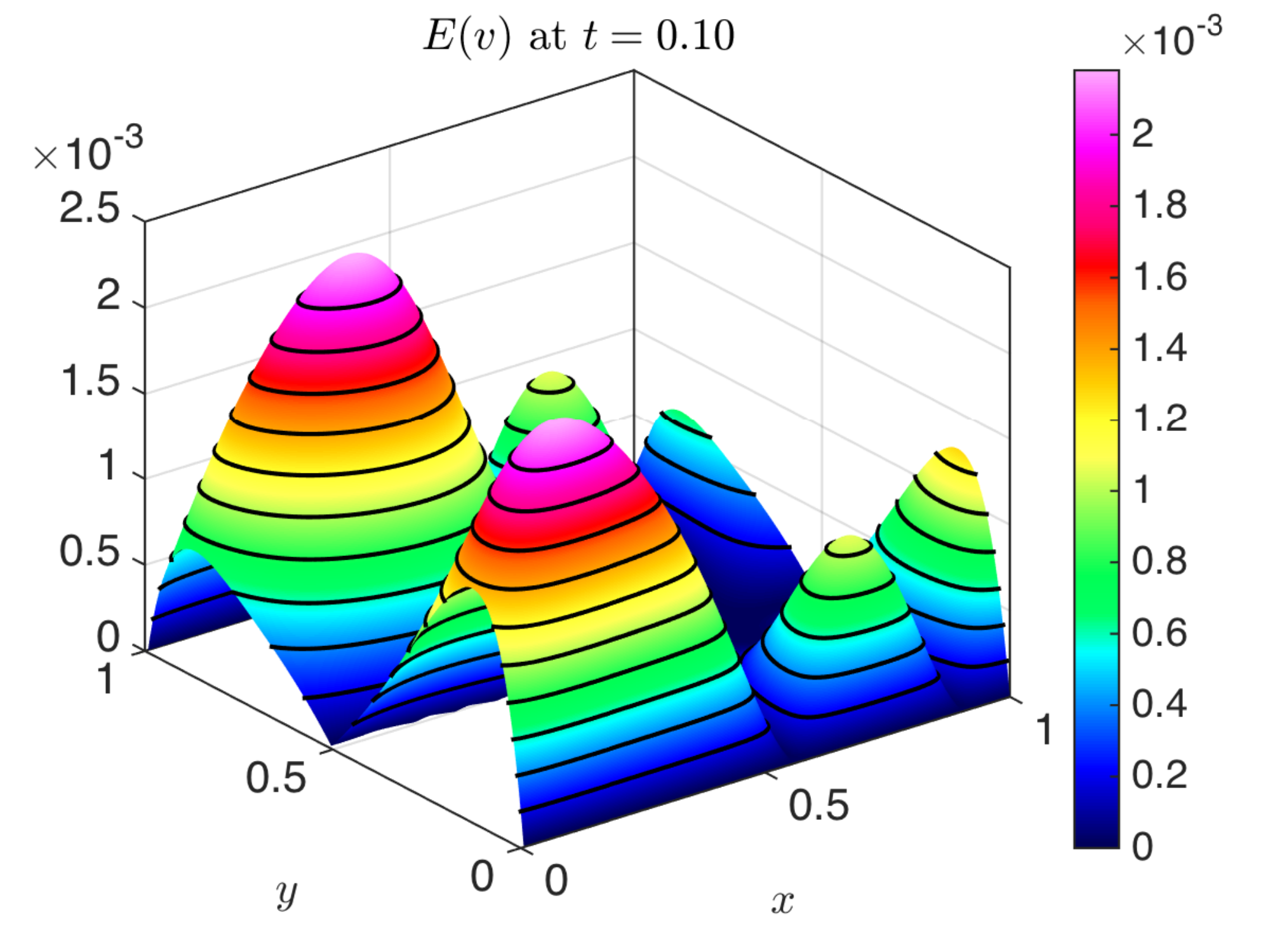}{\figWidth}};
\draw(6.8,13.0) node[anchor=south west,xshift=0pt,yshift=0pt] {\trimfig{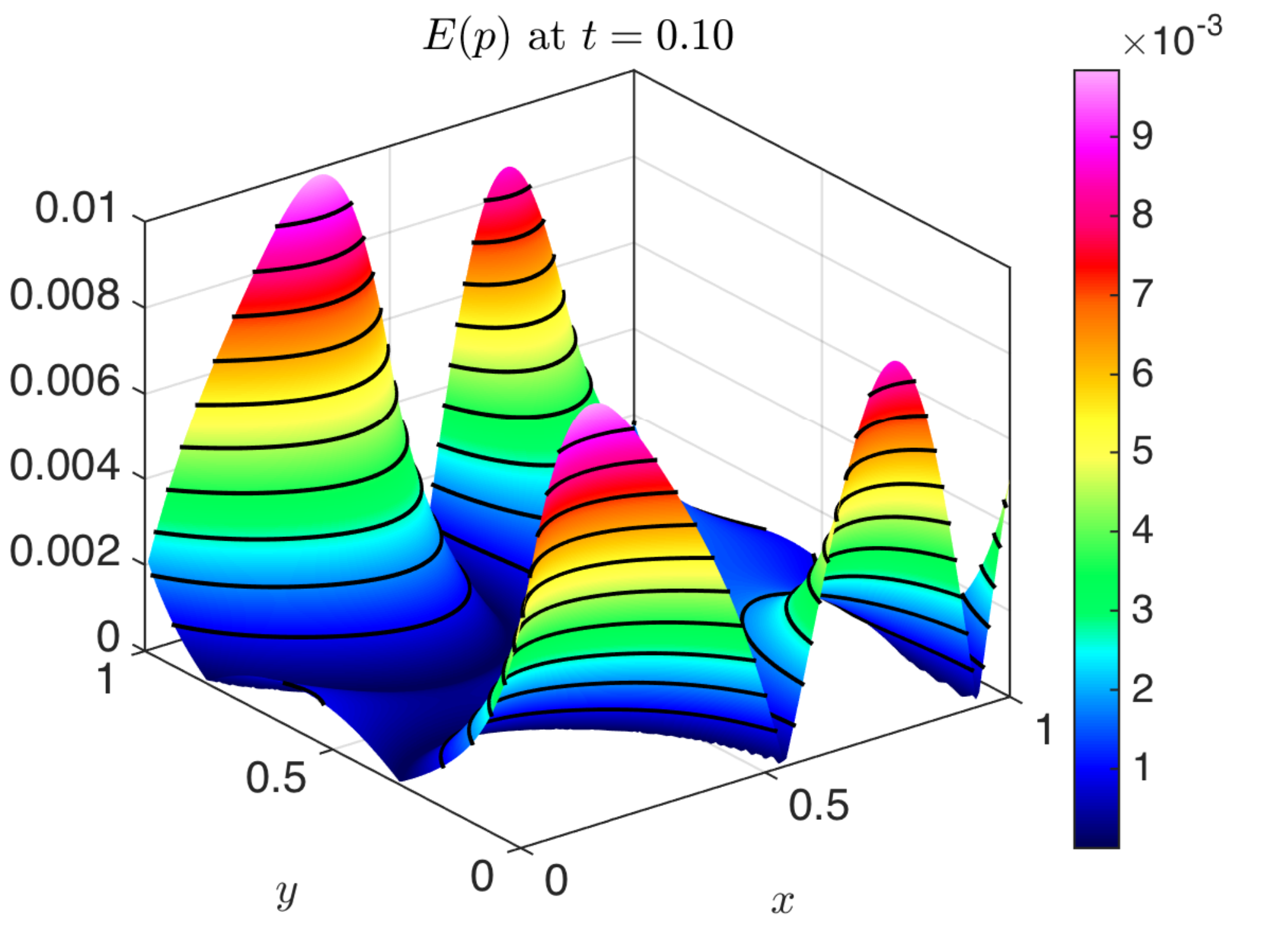}{\figWidth}};

\draw(-0.5,6.5) node[anchor=south west,xshift=0pt,yshift=0pt] {\trimfig{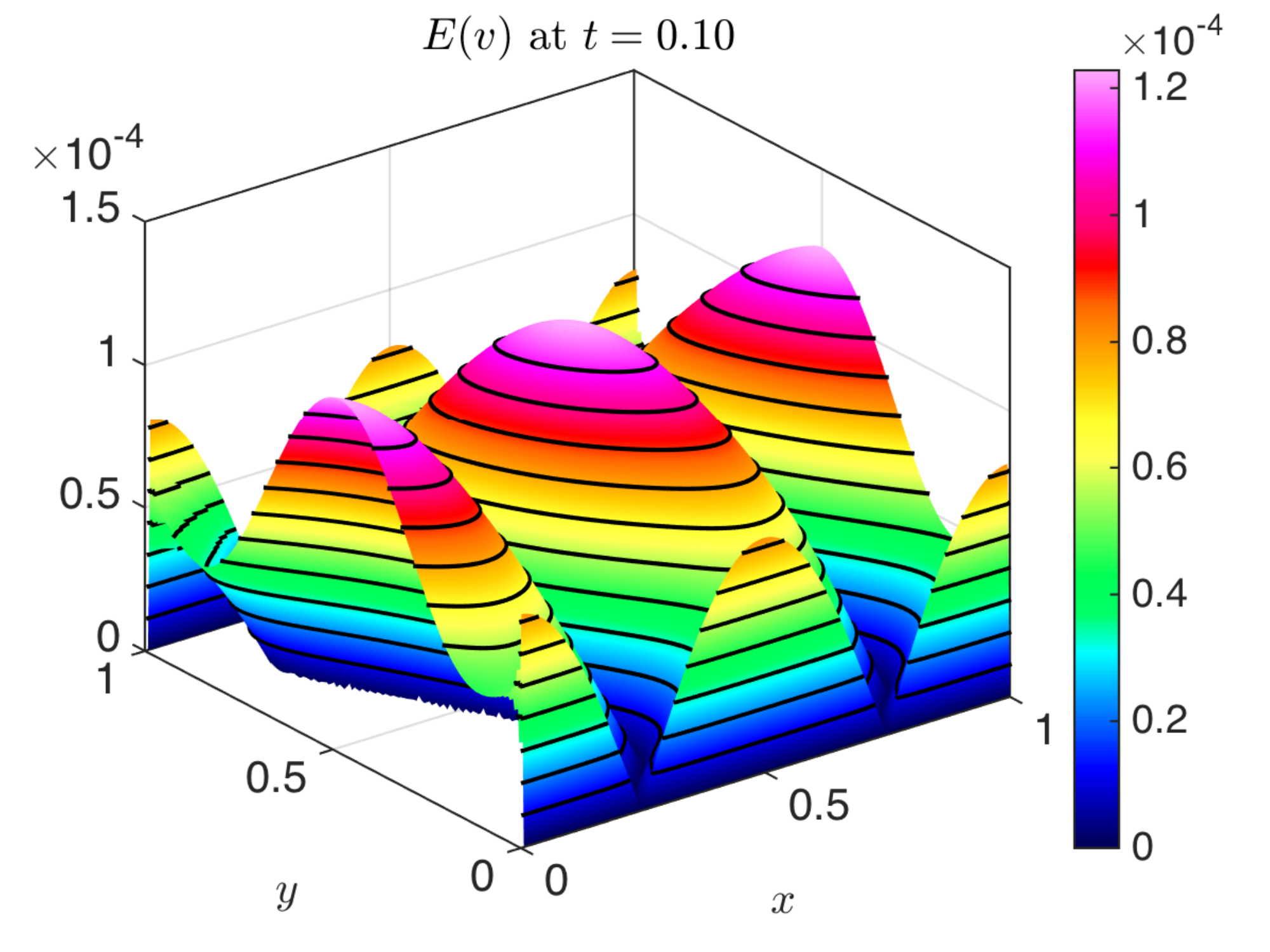}{\figWidth}};
\draw(6.8,6.5) node[anchor=south west,xshift=0pt,yshift=0pt] {\trimfig{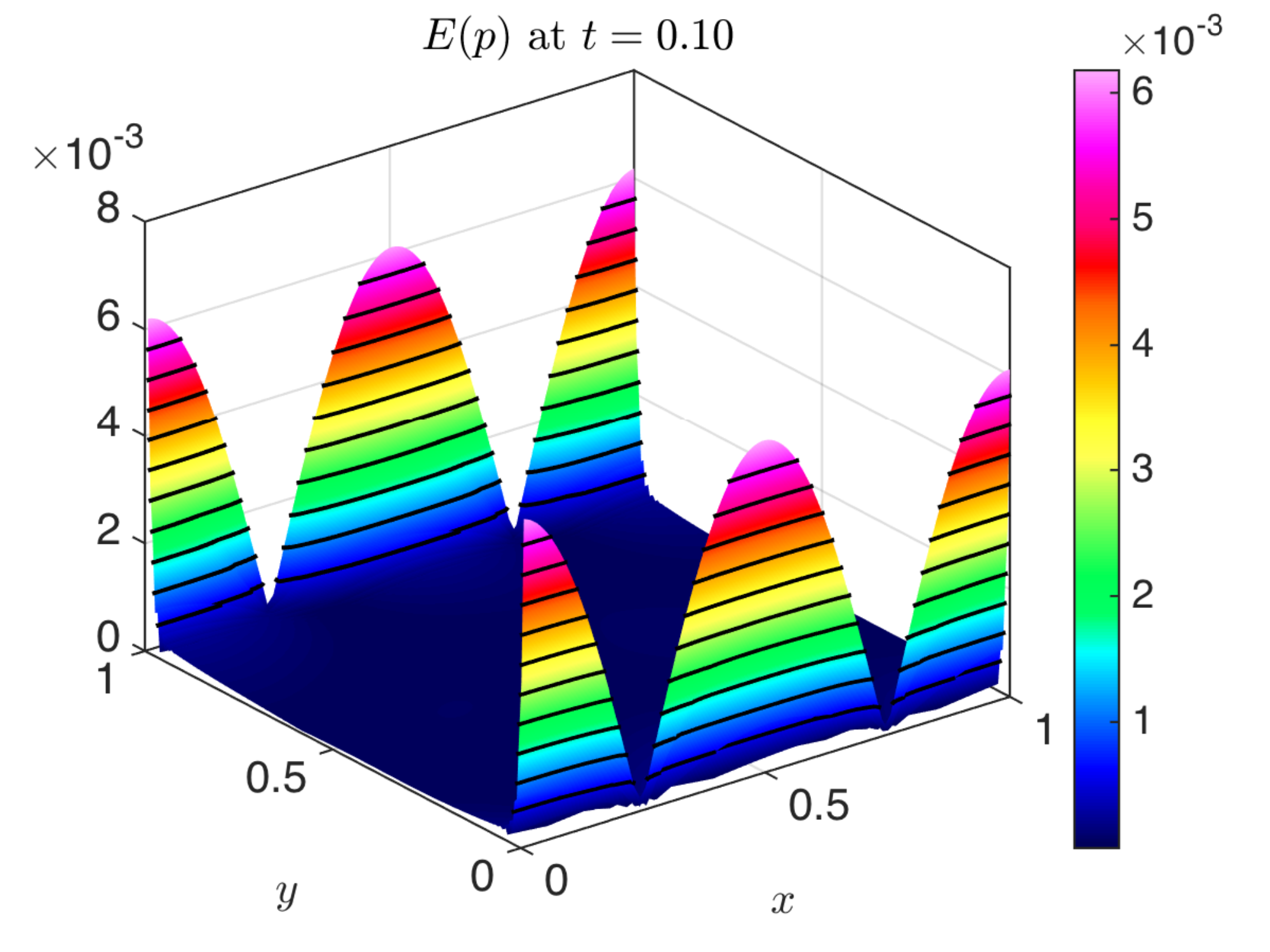}{\figWidth}};

\draw(-0.5,0.0) node[anchor=south west,xshift=0pt,yshift=0pt] {\trimfig{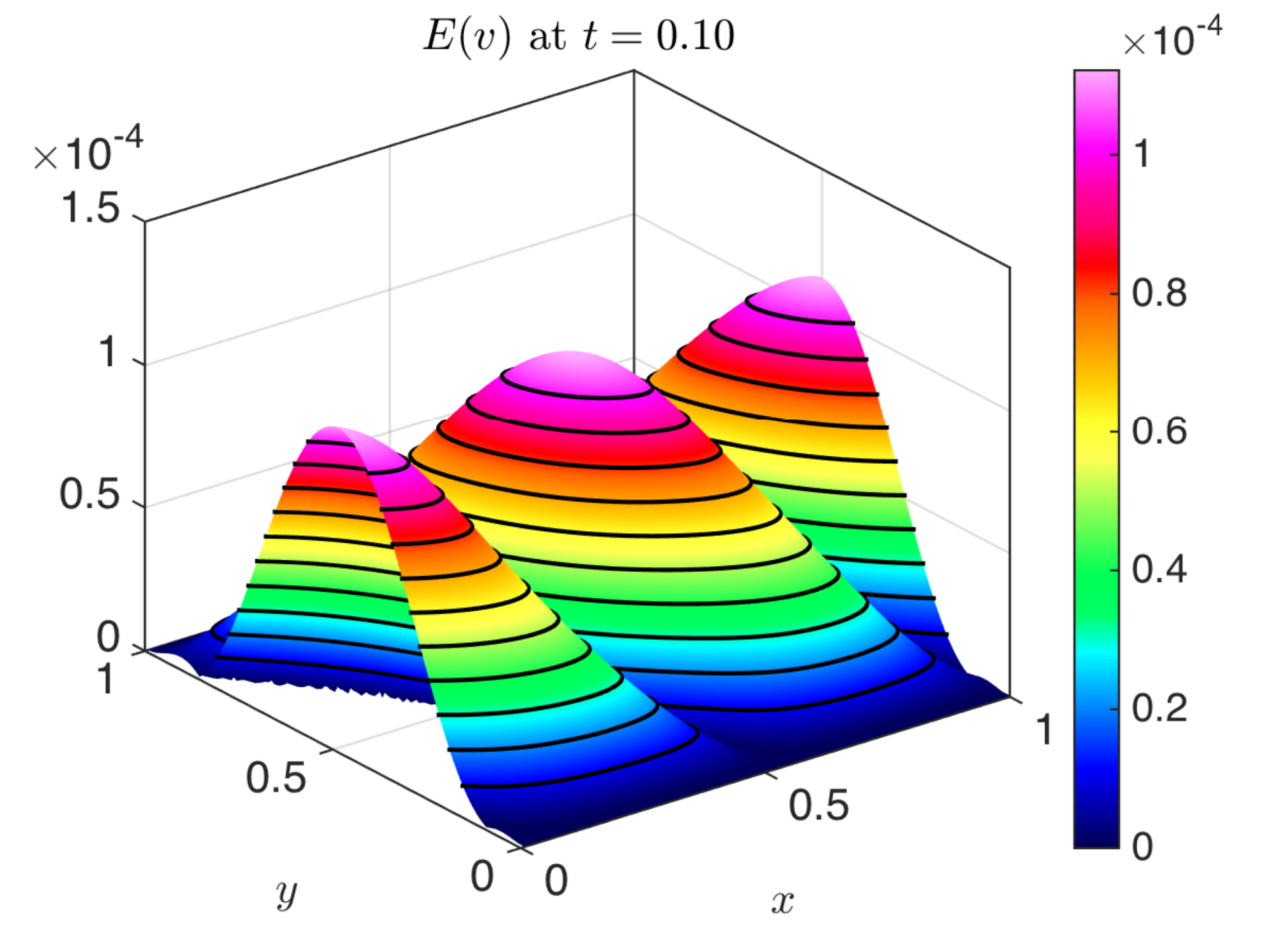}{\figWidth}};
\draw(6.8,0.0) node[anchor=south west,xshift=0pt,yshift=0pt] {\trimfig{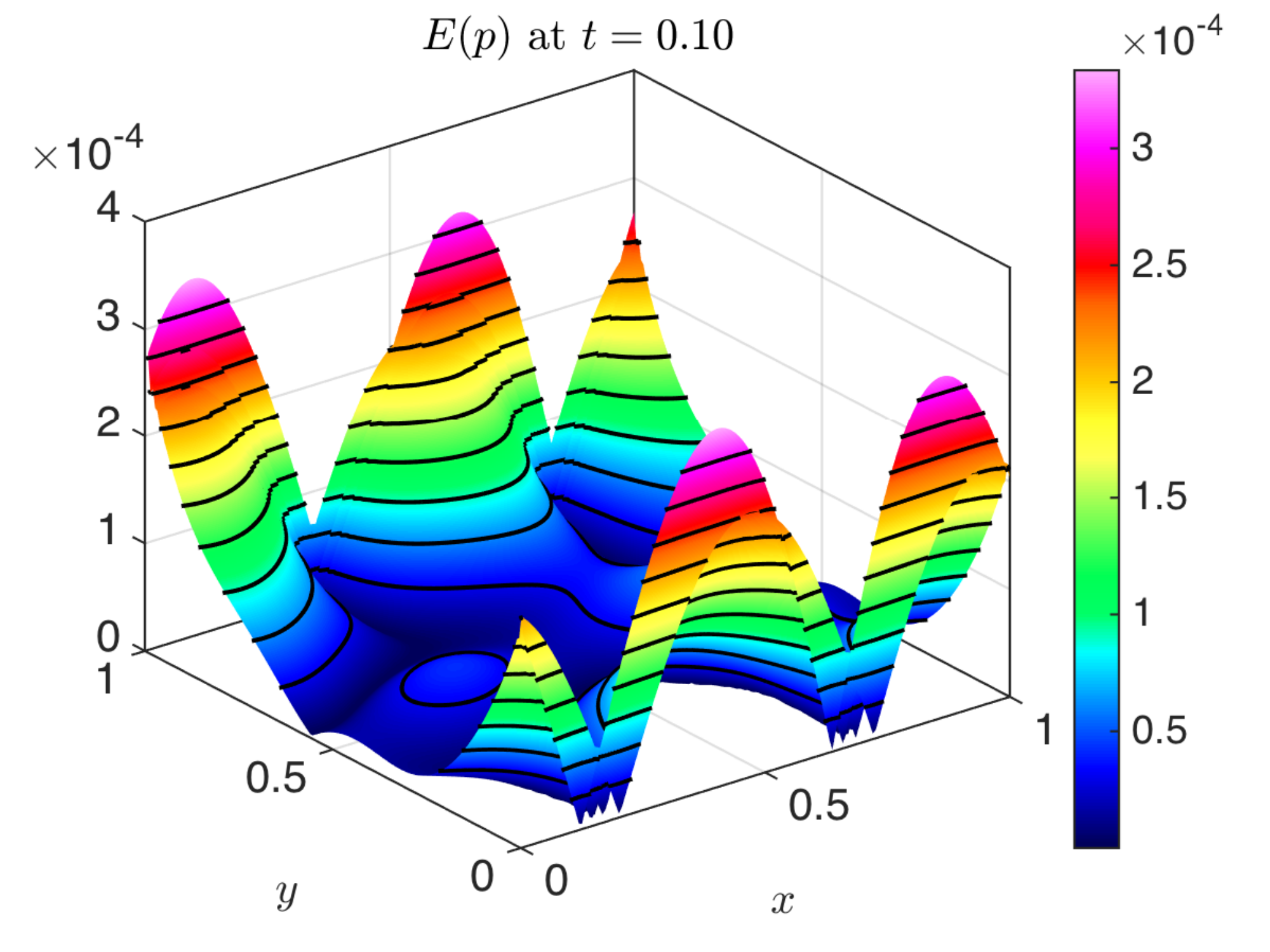}{\figWidth}};

\draw(7,18.5) node[anchor=south]{case (i): $\alpha=0$ with $\oldBC$ boundary condition};

\draw(7,12.) node[anchor=south]{case (ii): $\alpha=1/h^2$ with $\oldBC$ boundary condition};

\draw(7,5.5) node[anchor=south]{case (iv): $\alpha=1/h^2$ with $\newBC$ boundary condition};

%
\end{tikzpicture}

\end{center}
\caption{Absolute values of the errors for the velocity component $v$ and the pressure $p$ are plotted, i.e., $E(v)=|v-v_e|$ and $E(p)=|p-p_e|$. Periodicity in x direction is assumed, while   no-slip boundary conditions are enforced  on the boundaries in y direction. The grid spacing is $h=1/160$ and the time for  plot is $t=0.1$.}\label{fig:errorPlotsPD}
\end{figure}
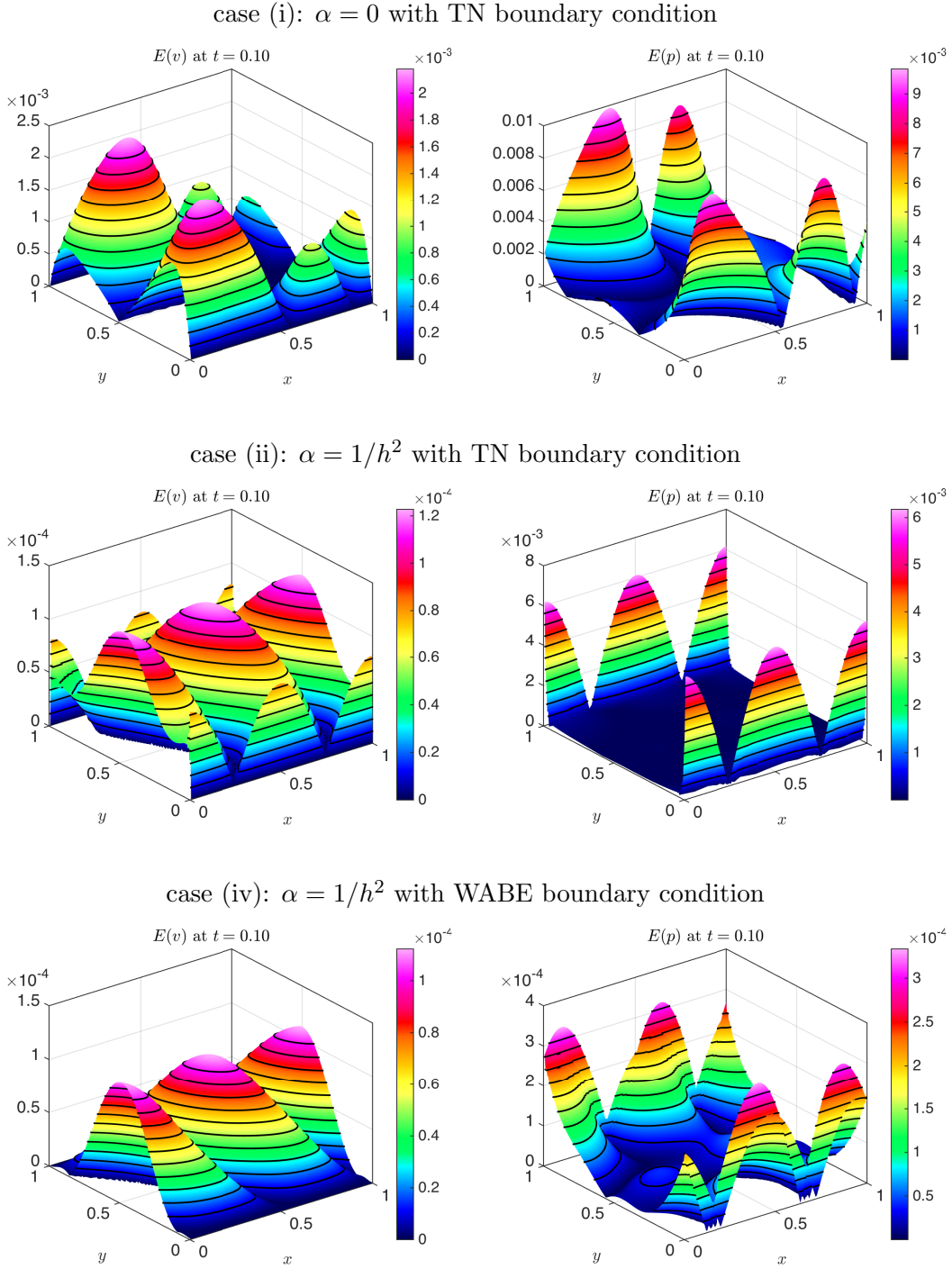
}

{
  \newcommand{\figWidth}{6cm}
  \def\xa{13}
  \def\ya{11.5}
\newcommand{\trimfig}[2]{\trimw{#1}{#2}{0.}{0.}{0.}{0.0}}
\begin{figure}[h]
\begin{center}
\begin{tikzpicture}[scale=1]
\useasboundingbox (0.0,0.0) rectangle (\xa,\ya);  
\draw(-0.5,6) node[anchor=south west,xshift=0pt,yshift=0pt] {\trimfig{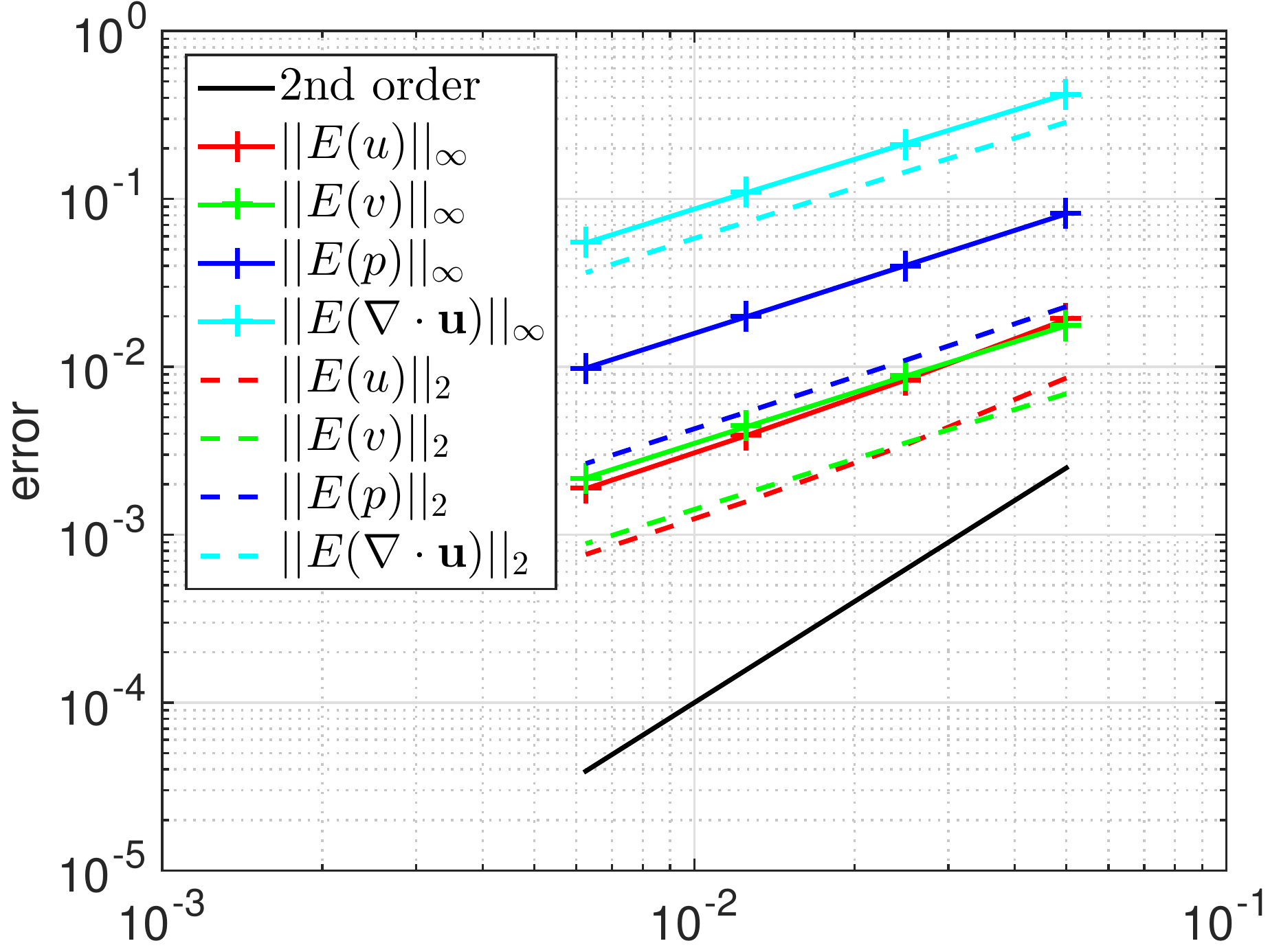}{\figWidth}};
\draw(6.5,6) node[anchor=south west,xshift=0pt,yshift=-0pt] {\trimfig{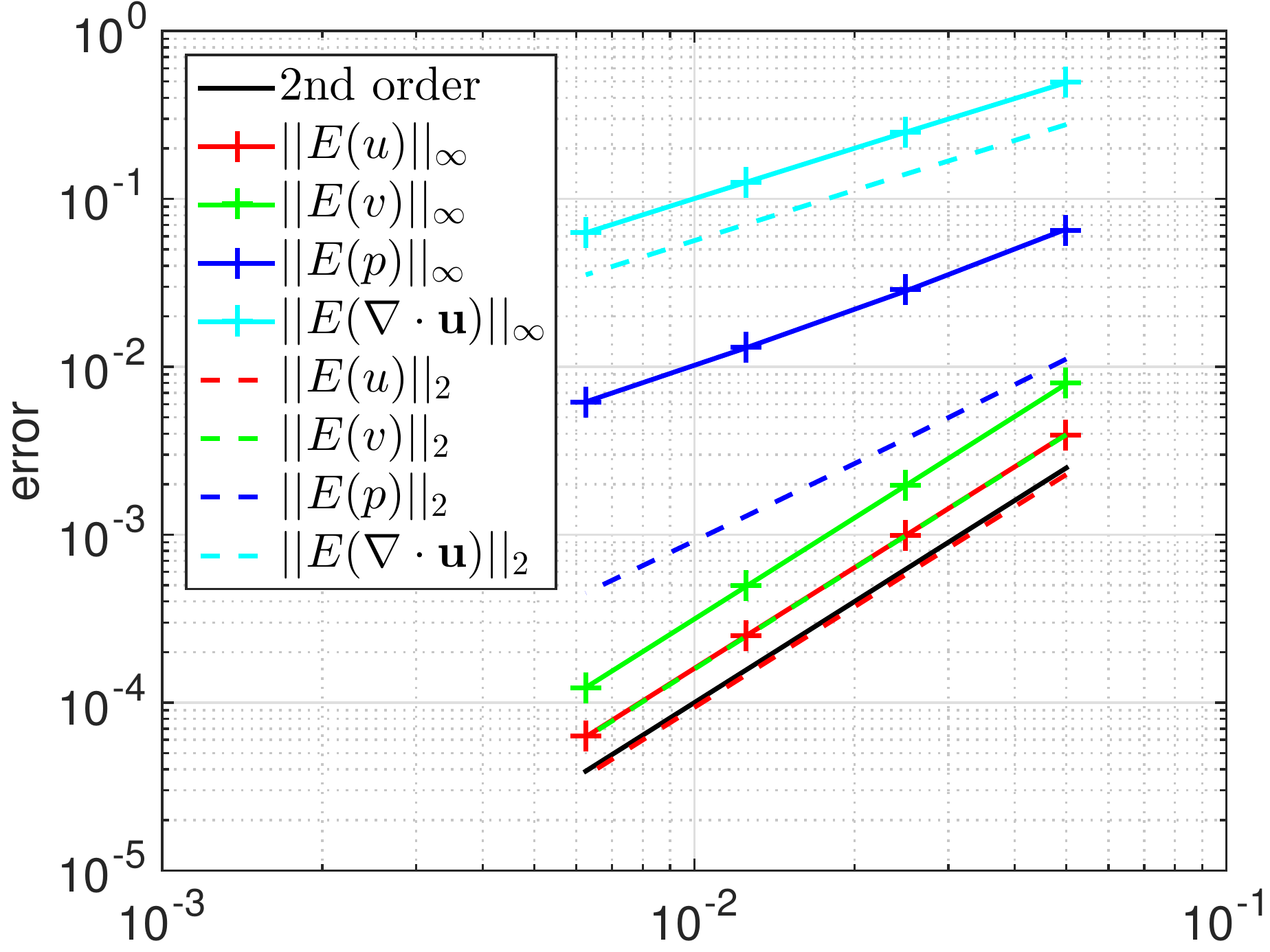}{\figWidth}};
\draw(-0.5,0) node[anchor=south west,xshift=0pt,yshift=0pt] {\trimfig{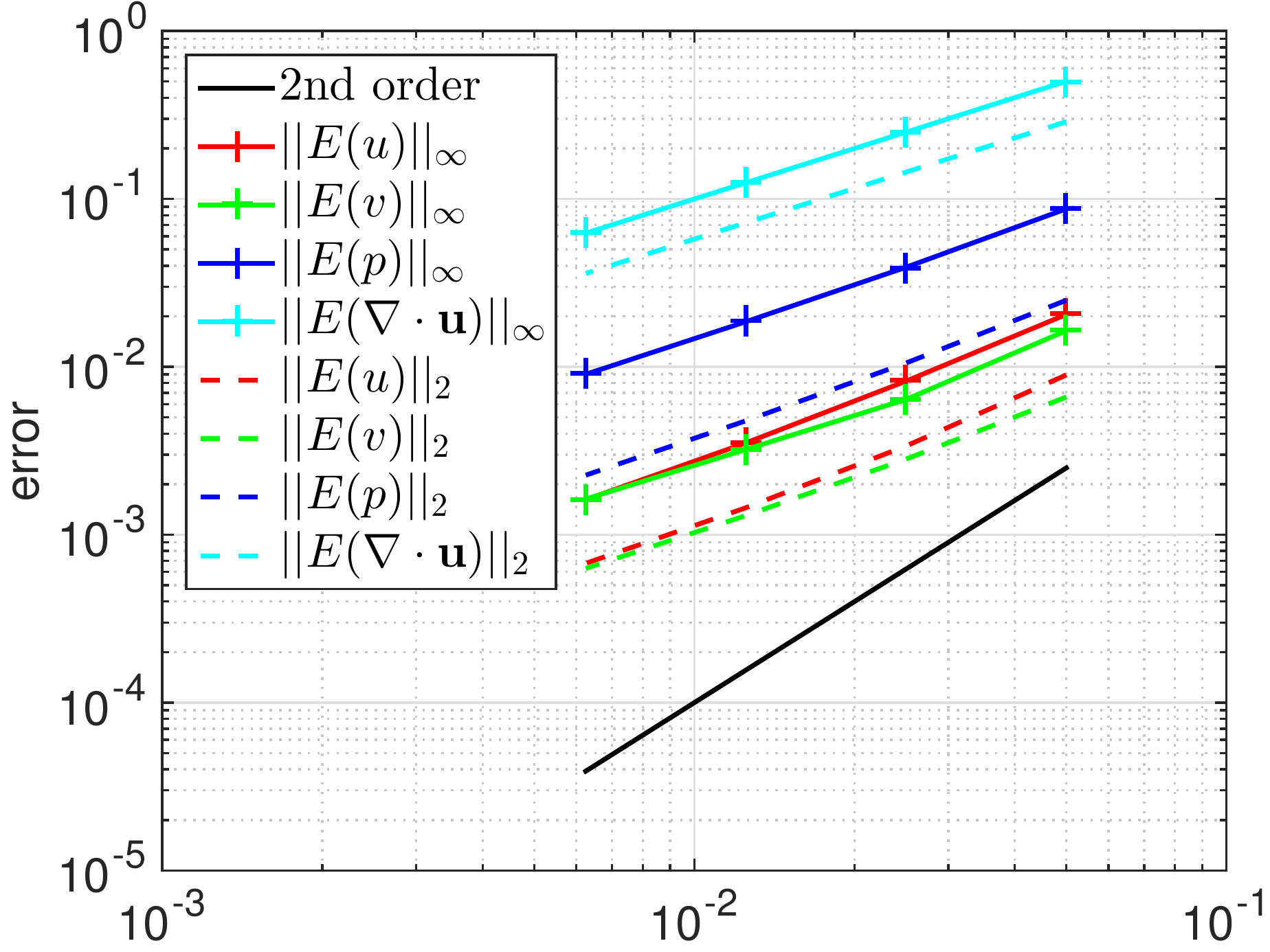}{\figWidth}};
\draw(6.5,0.0) node[anchor=south west,xshift=0pt,yshift=0pt] {\trimfig{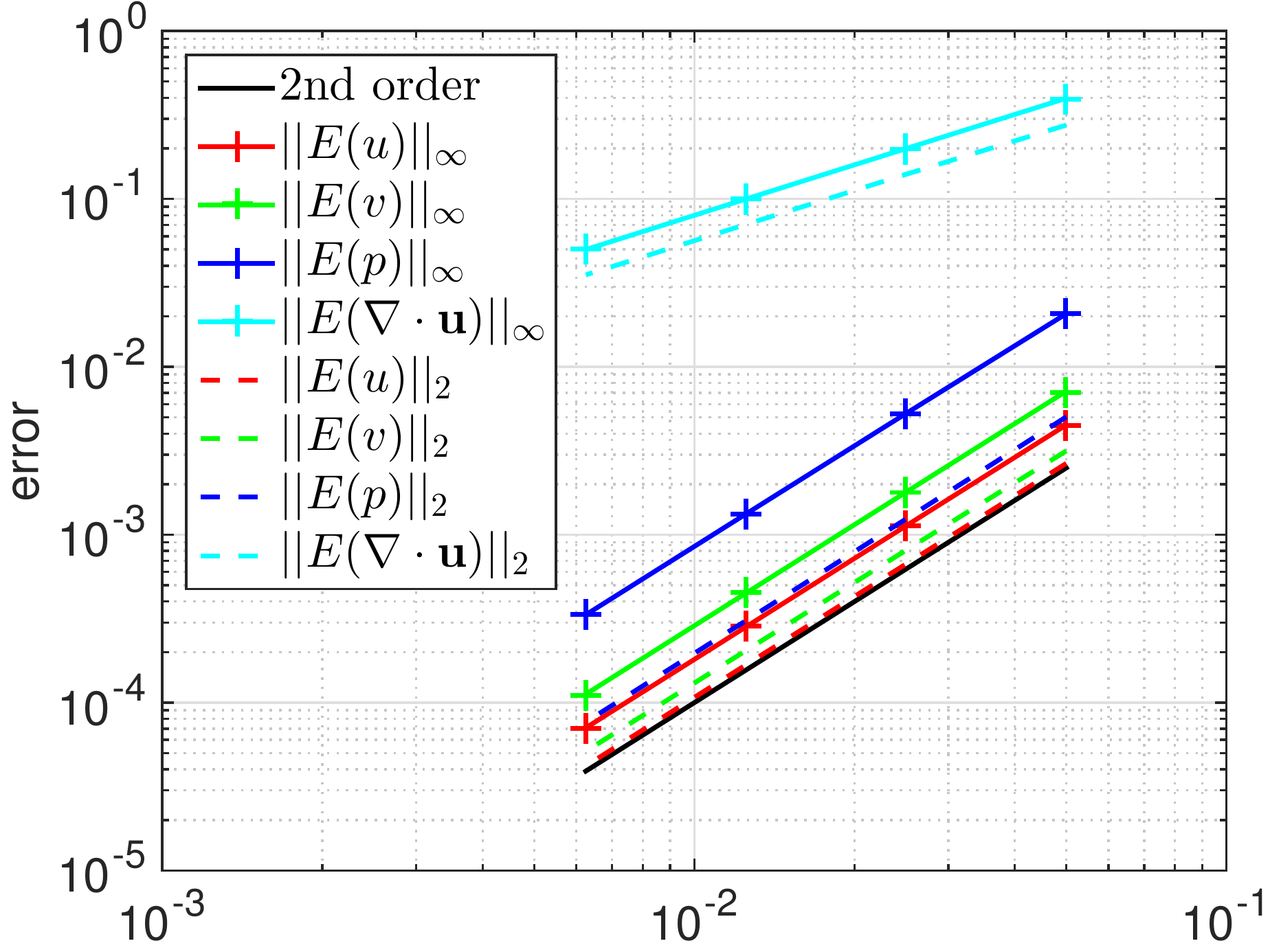}{\figWidth}};
%
\draw(3.,11) node[anchor=center]{case (i): $\alpha=0$ with $\oldBC$ BC};
\draw(10.25,11) node[anchor=center]{case (ii): $\alpha=1/h^2$ with $\oldBC$ BC};
\draw(3.,5) node[anchor=center]{case (iii): $\alpha=0$ with $\newBC$ BC};
\draw(10.25,5) node[anchor=center]{case (iv): $\alpha=1/h^2$ with $\newBC$ BC};

%
\end{tikzpicture}

\end{center}
\caption{Convergence rates of $u,v,p$ and  $\nabla\cdot\uv$ are plotted  at $t=0.1$. The errors are measured using both $L_\infty$ and $L_2$ norms with the former represented by solid lines and the latter by dashed lines. Periodicity in x direction is assumed, while   no-slip boundary conditions are enforced  on the boundaries in y direction.}\label{fig:convRatePD}
\end{figure}
}

A convergence study  for all the four cases is shown in Figure~\ref{fig:convRatePD}.  We see that, without divergence damping (i.e., $\alpha=0$), the errors for all the components  are  about first order accurate regardless   of  the pressure boundary conditions used. For the cases with the  divergence damping  turned on (i.e., $\alpha=1/h^2$), we observe that the  accuracies for the velocity components are second order for both TN and WABE boundary conditions. However, the pressure accuracy is first order in $L_\infty$ norm and a little bit inferior to second order in $L_2$ norm if TN boundary condition is implemented for pressure;  in contrast, the pressure accuracy is second order in both norms if WABE boundary condition is used. The performance of the TN boundary condition can be explained by looking at the second row of the error plots in Figure~\ref{fig:errorPlotsPD}. We see that it   improves the accuracy in the interior of the domain by adding divergence damping, so boundary layers are observed in both $E(v)$ and $E(p)$.  However, since the velocity error is still dominated by the interior, we observe second order accuracy in both norms, while the error for the pressure is dominated by the boundary-layer errors  so we observe first order in $L_\infty$ norm and almost second order in  $L_2$ norm.

We note that the divergence $\nabla\cdot\uv$ is always first order.  This is because we simply evaluate $\nabla\cdot\uv$ from the finite element solution for the velocity, and it is well-known that the derivative of a function represented in $\Pe_1$ elements is of first order accuracy. This can be improved by using some other post-processing techniques that compute  $\nabla\cdot\uv$ more accurately.  However, since we only want to keep track of the magnitude of the  divergence to make sure it remains small throughout the computation and it does not affect the accuracy of our scheme, it suffices for us to stick with this simple approach.

 \subsubsection{No-slip boundary conditions on all boundaries}
 We then consider the convergence study with no-slip boundary conditions  enforced on all boundaries with the rate of convergence for all four cases shown in Figure~\ref{fig:convRateNonPD}. Similar convergence properties are observed for the cases (i), (ii) and (iii). For  case (iv) (the lower right plot in Figure~\ref{fig:convRateNonPD}),  we still observe second order accuracy for the solutions $u,v$ and $p$  in $L_2$ norm; however, we see only first order accuracy for $p$ in $L_\infty$ norm.  By looking at the error plots   in Figure~\ref{fig:errorWABEddNonPD} for this case,  we see steep gradients near the corners for both $E(v)$ and $E(p)$.  For $E(v)$, the error is still dominated by the interior; however, for $E(p)$, the error is dominated by the corner spikes. Thus,  the max-norm error for the pressure is strongly affected by its  behavior in the corners.  We note that  similar corner behavior for the pressure solution is also reported in \cite{LiuLiuPego2010b} for numerical examples computed on domains with sharp corners.    In fact, there are two independent normals at corners where two edges meet; these independent normals cause trouble when implementing Neumann boundary conditions there.  Thus, we can also see the corner spikes in the error plot when solving Poisson equation with pure Neumann boundary condition on a square domain using finite element method.  Fortunately,  the corner issue only affects the accuracy of the scheme locally. Globally, the scheme is still well-behaved and   second-order accurate, which is confirmed by the $L_2$-norm errors.
 
 {
  \newcommand{\figWidth}{6cm}
  \def\xa{13}
  \def\ya{11.5}
\newcommand{\trimfig}[2]{\trimw{#1}{#2}{0.}{0.}{0.}{0.0}}
\begin{figure}[h]
\begin{center}
\begin{tikzpicture}[scale=1]
\useasboundingbox (0.0,0.0) rectangle (\xa,\ya);  
\draw(-0.5,6) node[anchor=south west,xshift=0pt,yshift=0pt] {\trimfig{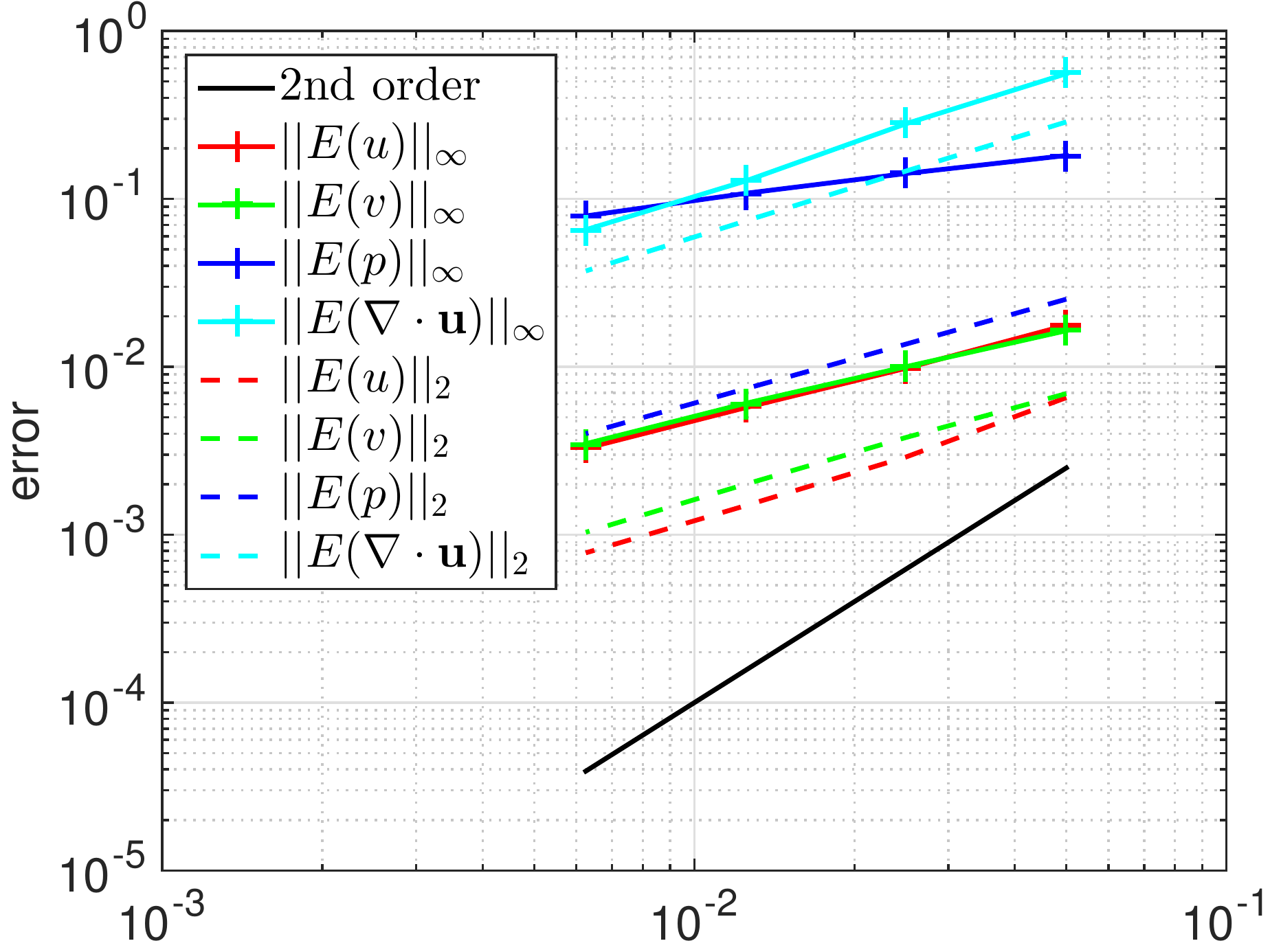}{\figWidth}};
\draw(6.5,6) node[anchor=south west,xshift=0pt,yshift=-0pt] {\trimfig{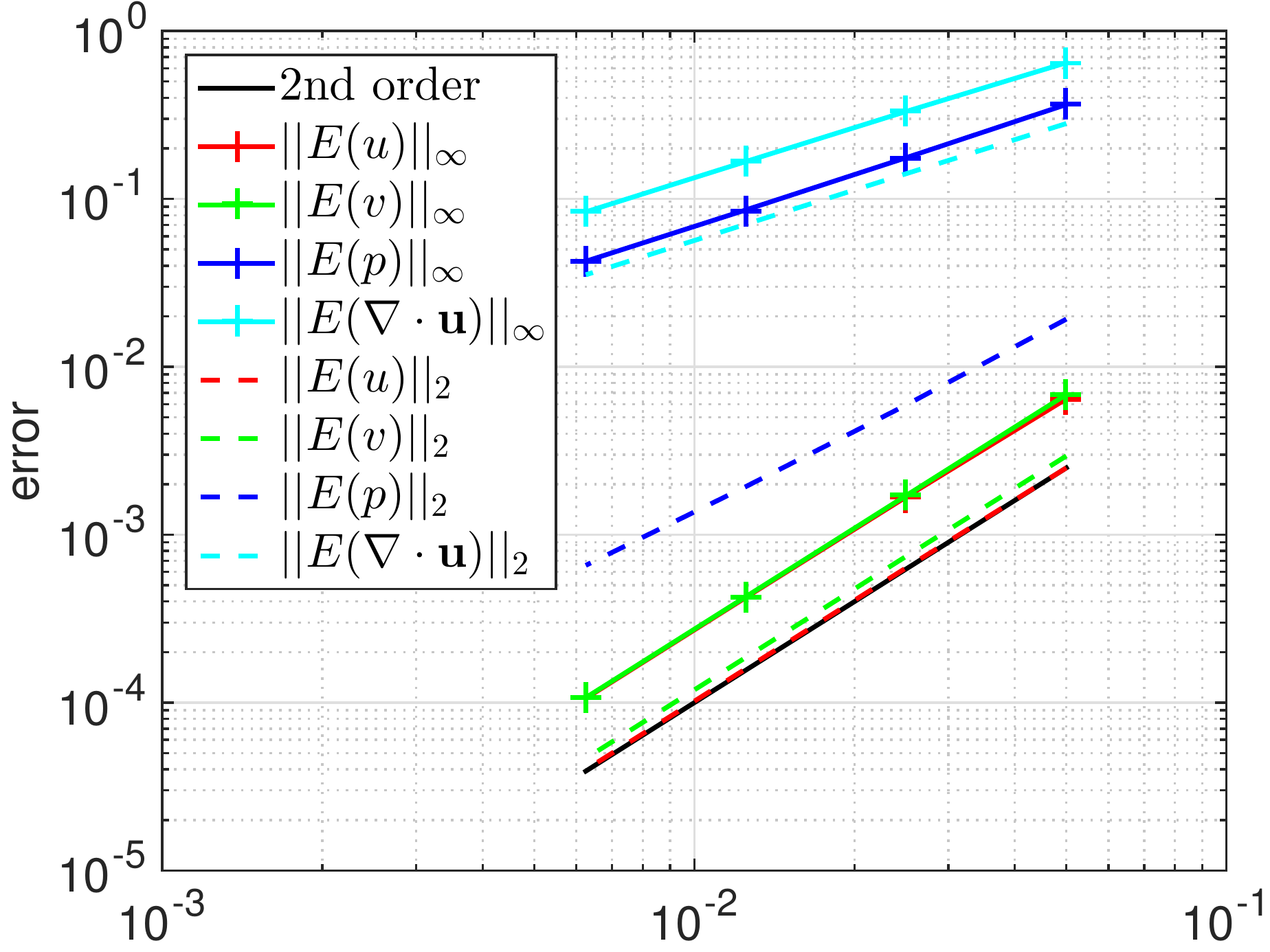}{\figWidth}};
\draw(-0.5,0) node[anchor=south west,xshift=0pt,yshift=0pt] {\trimfig{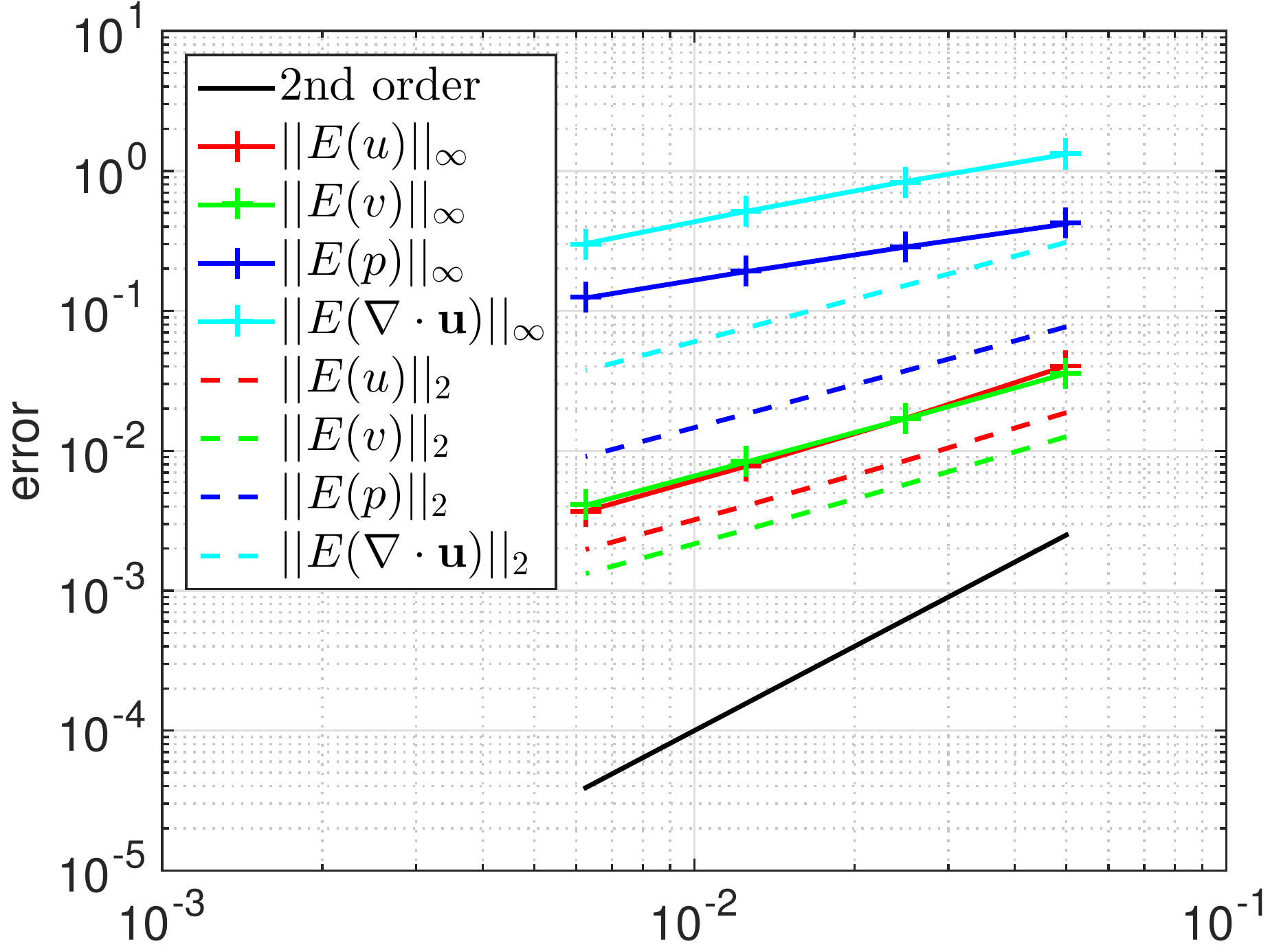}{\figWidth}};
\draw(6.5,0.0) node[anchor=south west,xshift=0pt,yshift=0pt] {\trimfig{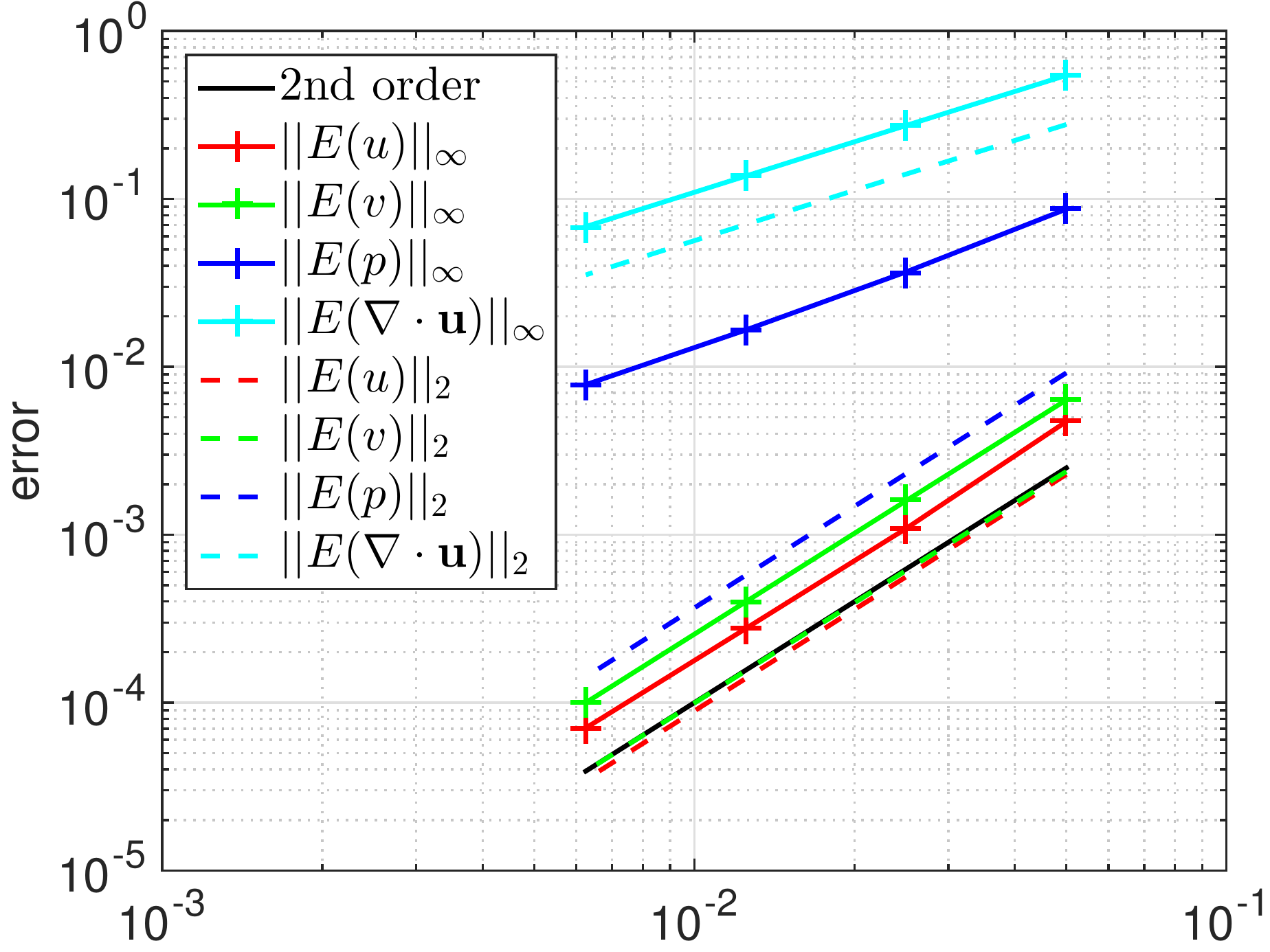}{\figWidth}};
%
\draw(3.,11) node[anchor=center]{case (i): $\alpha=0$ with $\oldBC$ BC};
\draw(10.25,11) node[anchor=center]{case (ii): $\alpha=1/h^2$ with $\oldBC$ BC};
\draw(3.,5) node[anchor=center]{case (iii): $\alpha=0$ with $\newBC$ BC};
\draw(10.25,5) node[anchor=center]{case (iv): $\alpha=1/h^2$ with $\newBC$ BC};

%
\end{tikzpicture}

\end{center}
\caption{Convergence rates of $u,v,p$ and  $\nabla\cdot\uv$ are plotted  at $t=0.1$. The errors are measured using both $L_\infty$ and $L_2$ norms with the former represented by solid lines and the latter by dashed lines. No-slip boundary conditions are enforced  on all the boundaries.}\label{fig:convRateNonPD}
\end{figure}
}
 
{
\newcommand{\figWidth}{7cm}
\newcommand{\trimfig}[2]{\trimw{#1}{#2}{0.}{0.}{0.}{0.0}}
\begin{figure}[h]
\begin{center}
\begin{tikzpicture}[scale=1]
\useasboundingbox (0.0,0.0) rectangle (14.,6);  
\draw(-0.5,0.0) node[anchor=south west,xshift=0pt,yshift=0pt] {\trimfig{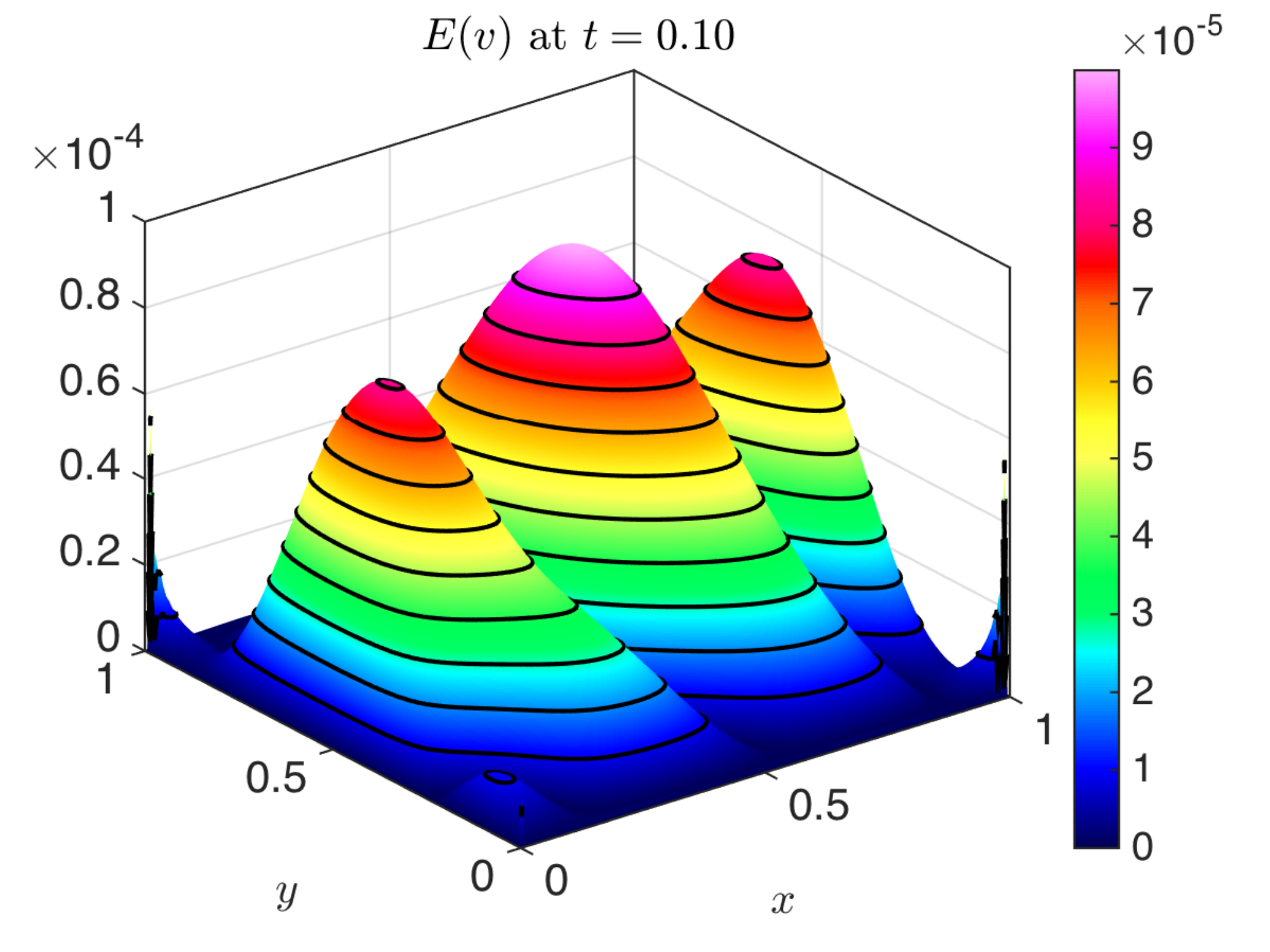}{\figWidth}};
\draw(6.8,0.0) node[anchor=south west,xshift=0pt,yshift=0pt] {\trimfig{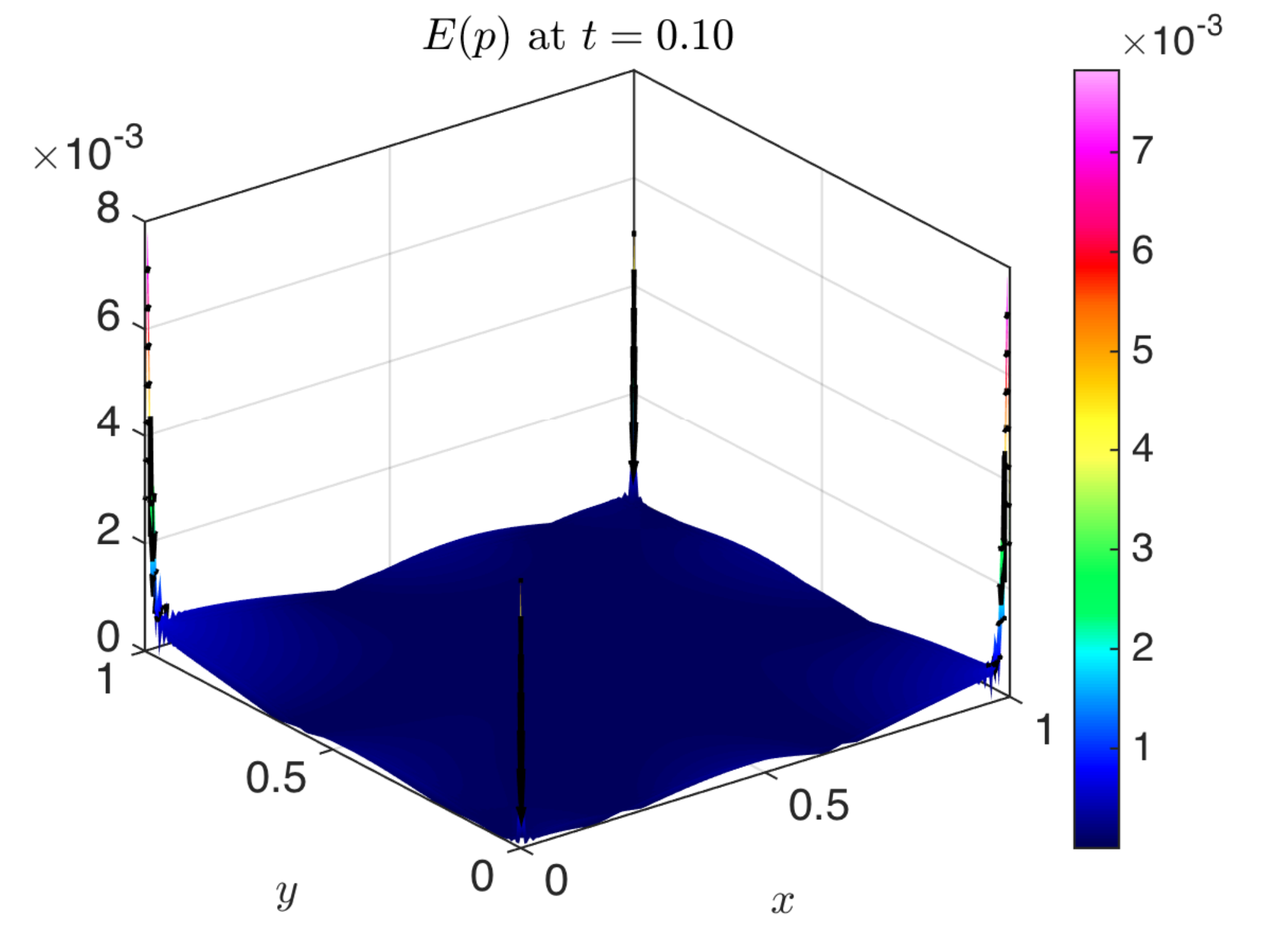}{\figWidth}};

\draw(7,5.5) node[anchor=south]{case (iv): $\alpha=1/h^2$ with $\newBC$ BC};

%
\end{tikzpicture}

\end{center}
\caption{Absolute values of the errors for the velocity component $v$ and the pressure $p$ are plotted, i.e., $E(v)=|v-v_e|$ and $E(p)=|p-p_e|$. We  enforce  no-slip boundary conditions at all the  boundaries. The grid spacing is $h=1/160$ and the time for  plot is $t=0.1$.}\label{fig:errorWABEddNonPD}
\end{figure}
}

\clearpage
\subsection{Modified Lid-Driven Cavity}
To further verify the accuracy and efficiency of the proposed scheme, we solve a modified  lid-driven cavity  problem. Specifically, we consider the flow in the  square  domain  $[0,1]\times[0,1]$. The associated  boundary conditions are $\uv= (u_0, 0)$ on the top of the domain (i.e., $y=1$) and $\uv= (0,0)$  on the other three sides. It is well known that the classical lid-driven cavity problem that specifies $u_0=1$  introduces  singularities at the corners since the horizontal  velocity at the top corners suddenly changes from $0$ to $1$. In spite of the singularities, this problem is popular for testing and evaluating numerical methods. We will  compare our results with that from  \cite{GiahEtal82, BotellaPeyret98, LiuLiuPego2010b}. For results reported in \cite{GiahEtal82, LiuLiuPego2010b}, the authors did nothing to  suppress the corner singularities.  However,  a spectral method has been used  in \cite{BotellaPeyret98} so the singularities at corners cannot be ignored  because   the ``spectral'' accuracy is   generally  associated with the smoothness of the solution.  To  minimize  the effect of the singularities,  the authors   extracted  analytically the corner singularities from the dependent variables of the problem.
Here we take a different approach to remove the singularities  by  modifying the boundary condition. Specifically, we define $u_0(x)$ such that its value smoothly transitions from $0$ to $1$ when $x$ is away from the ends, i.e.,
$$
u_0(x) =\frac{1}{2}\left[-\tanh\left(\frac{|x-0.5|-0.495}{0.01}\right)+1\right].
$$
In the right image of \Fig~\ref{fig:squareStretchedMesh}, we show the  horizontal velocity (u) on part of the top boundary grids near the top-left corner  (0,1) to illustrate how the boundary condition is smoothed over a couple of grid points.

The problem for  $\nu=1/1000$  is  solved using our numerical methods with both TN and WABE boundary conditions.  The  mesh used for computation consists of  4225 degrees of freedom (dof), and the  maximum  and minimum values of the grid spacings are $\max(h) = 0.027616$ and $\min(h) = 0.011561$, respectively. In Figure~\ref{fig:squareStretchedMesh},
a coarsened version of the computational mesh is  shown.  Note that  the  grids are stretched to cluster towards the boundaries for numerical purposes.

{
\newcommand{\figHeight}{5cm}
\newcommand{\trimfiga}[2]{\trimh{#1}{#2}{0.15}{0.15}{0.}{0.0}}
\newcommand{\trimfig}[2]{\trimh{#1}{#2}{0.}{0.}{0.}{0.0}}
\begin{figure}[h]
\begin{center}
\begin{tikzpicture}[scale=1]
  \useasboundingbox (0.0,0.0) rectangle (14.,5.5);  

\draw(0,0.0) node[anchor=south west,xshift=0pt,yshift=0pt] {\trimfiga{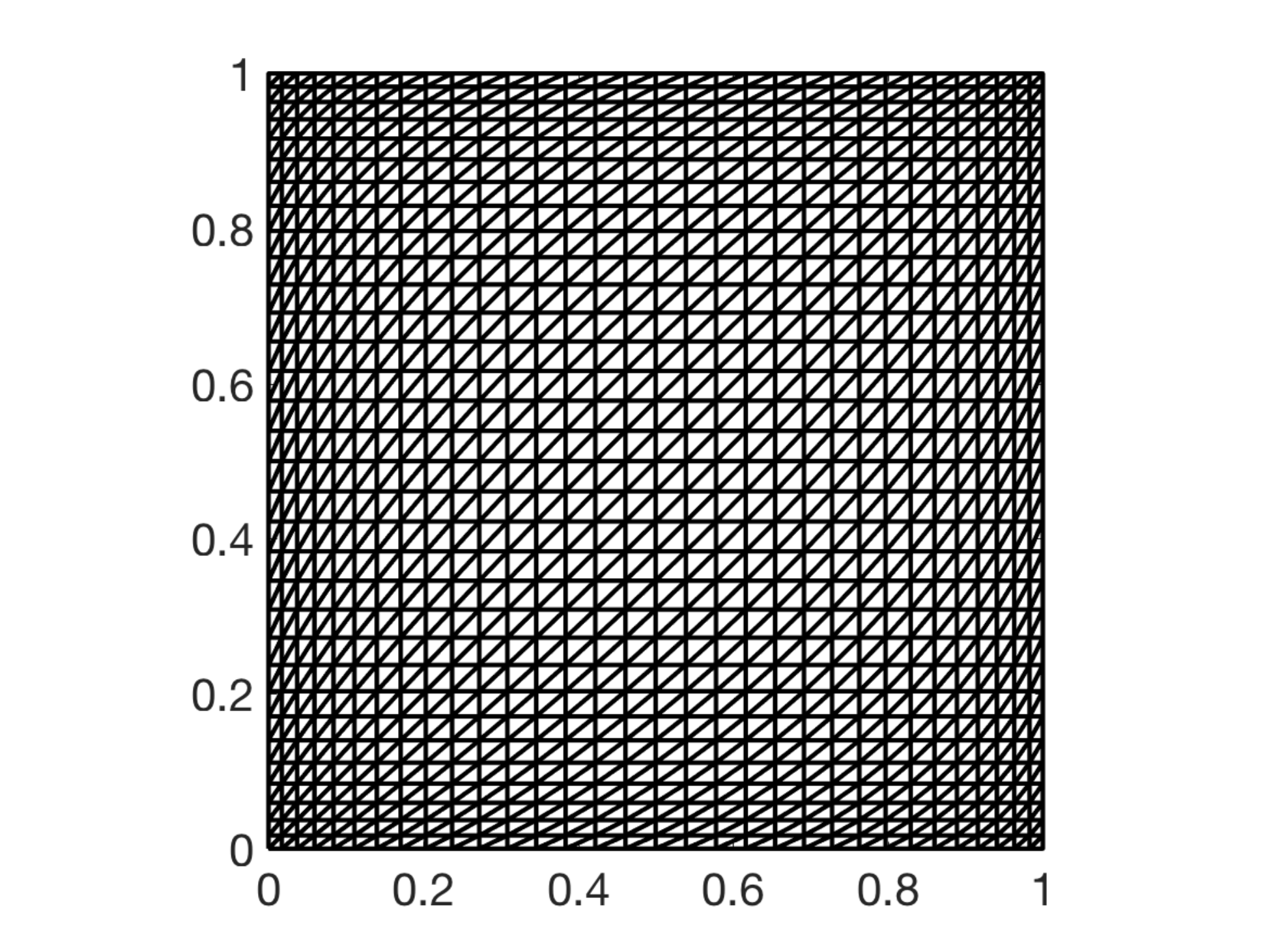}{\figHeight}};
\draw(6.8,0.0) node[anchor=south west,xshift=0pt,yshift=0pt] {\trimfig{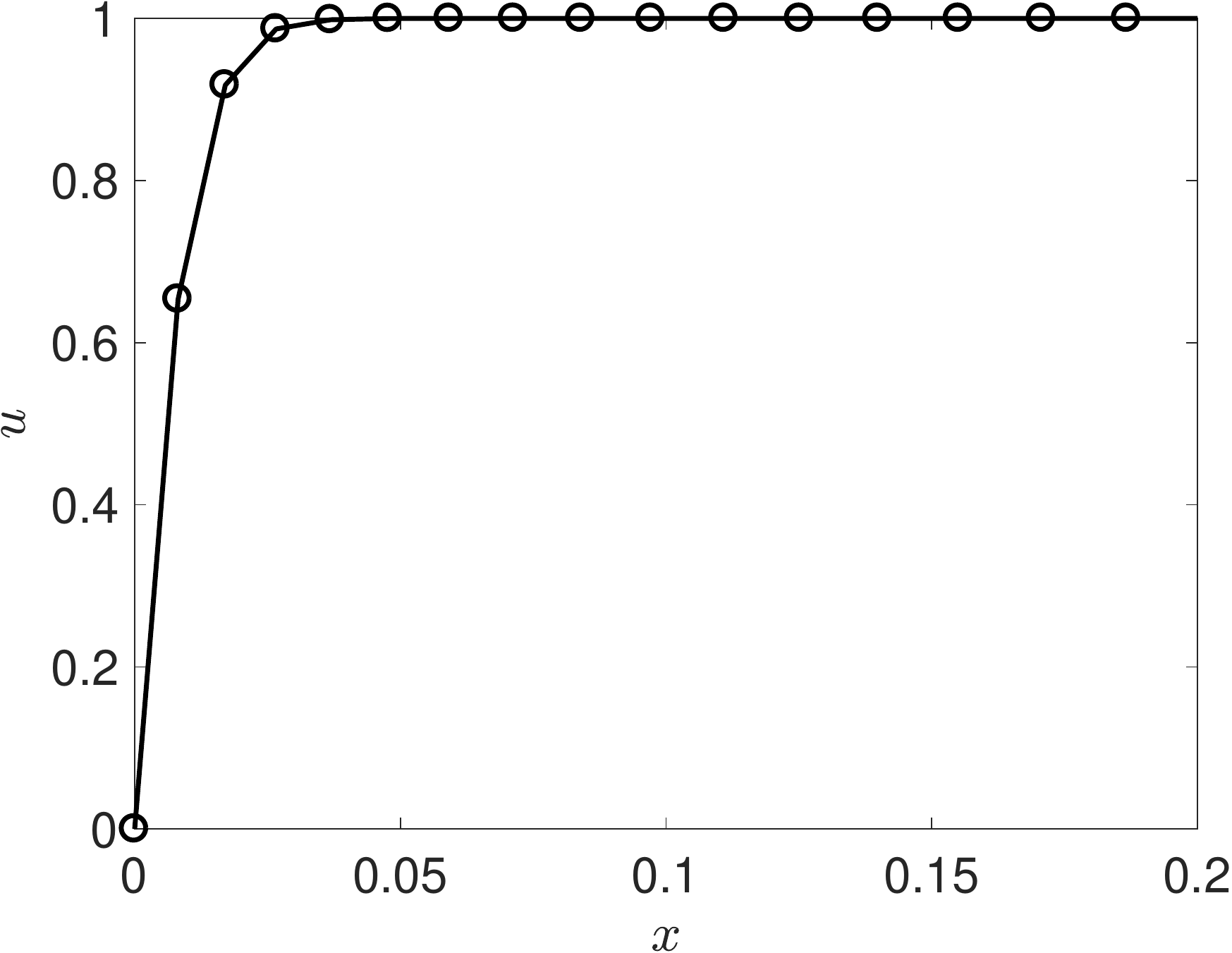}{\figHeight}};

%
\end{tikzpicture}

\end{center}
\caption{Left: a coarsened version of the computational mesh with  grids  stretched to cluster towards the boundaries. For the actual  mesh used for computation, we have  $\text{dof}=4225$,$\max(h)=0.027616$, $\min(h)=0.011561$. Right: horizontal velocity (u) on part of the top boundary grids near the top-left corner  (0,1).}\label{fig:squareStretchedMesh}
\end{figure}
}

In  Figure~\ref{fig:modifiedLDCSL}, we show the streamlines of the lid-driven cavity flow at $t=50$. And  in Figure~\ref{fig:modifiedLDCGhia}, we plot the velocity components $u$ and $v$ along the vertical and horizontal  lines through the geometric center; i.e., $u(0.5,y)$ and $v(x,0.5)$. Reference data from  \cite{GiahEtal82} are also plotted on top of our results for comparison. We see that our results using the proposed scheme with the divergence damping coefficient   $\alpha=1/h_{\min}^2$  match very well  with existing computations  reported in \cite{GiahEtal82, BotellaPeyret98, LiuLiuPego2010b}. 

{
\newcommand{\figWidth}{5cm}
\newcommand{\trimfig}[2]{\trimw{#1}{#2}{0.}{0.}{0.}{0.0}}
\begin{figure}[h]
\begin{center}
\begin{tikzpicture}[scale=1]
  \useasboundingbox (0.0,0.0) rectangle (11.,5);  


\draw(-1,0.0) node[anchor=south west,xshift=0pt,yshift=0pt] {\trimfig{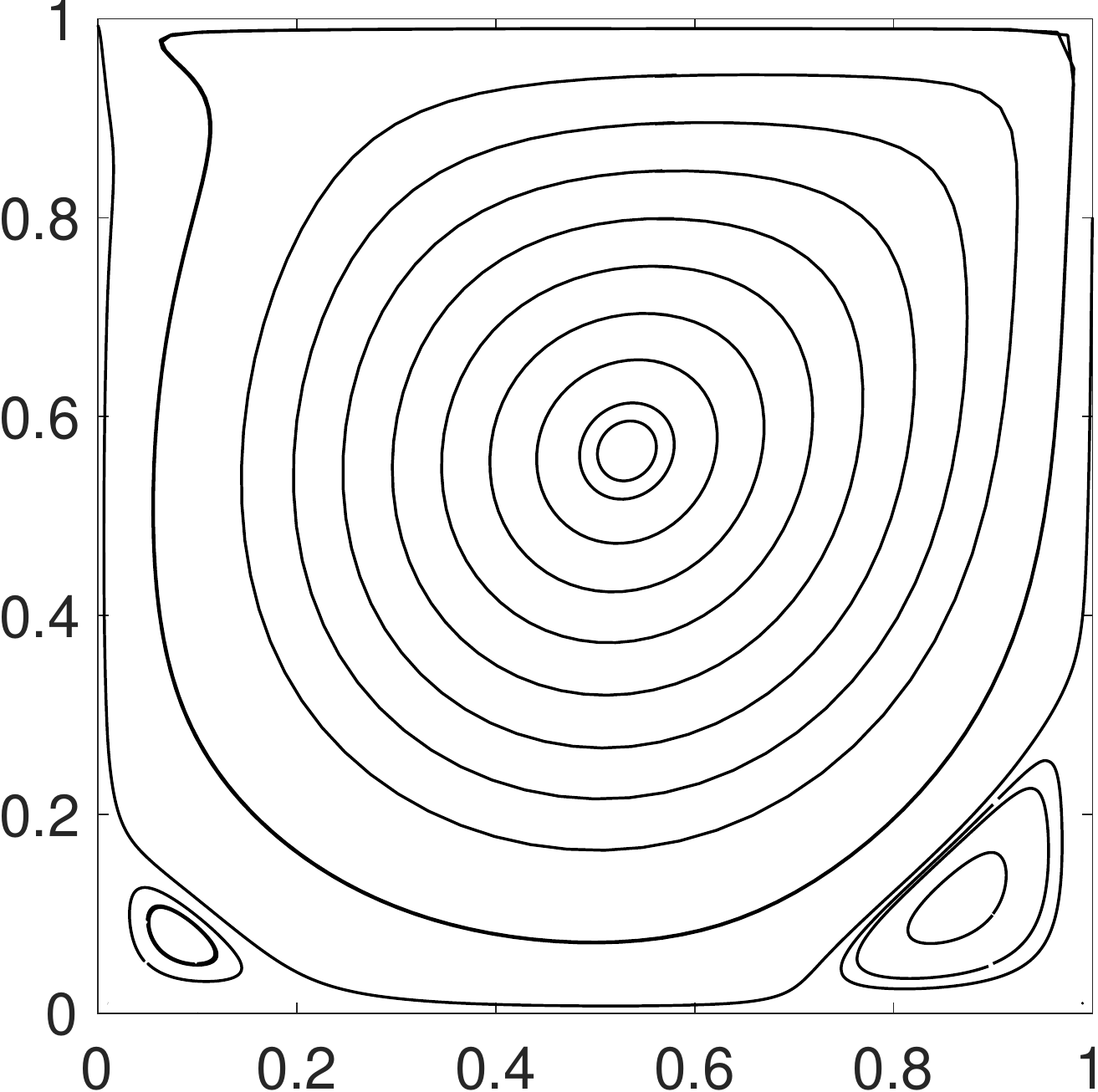}{\figWidth}};
\draw(5.5,0.0) node[anchor=south west,xshift=0pt,yshift=0pt] {\trimfig{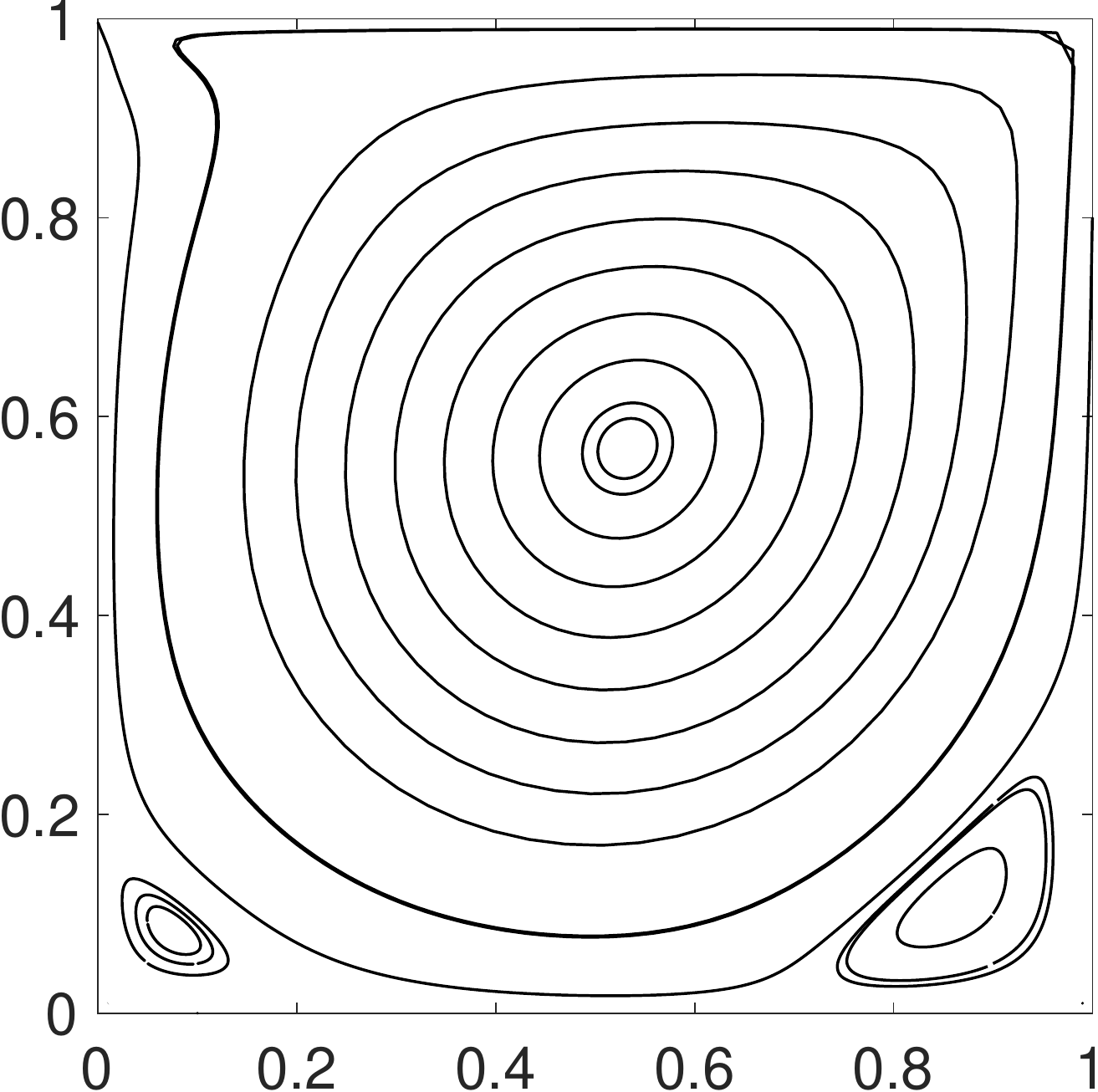}{\figWidth}};

%
\end{tikzpicture}

\end{center}
\caption{Streamline plots at $t=50$. Left: TN boundary condition. Right: WABE boundary condition.}\label{fig:modifiedLDCSL}
\end{figure}
}

Following \cite{BotellaPeyret98,LiuLiuPego2010b},  we also  show the vorticity contour at levels $[-5,-4,-3,-2,-1,-0.5,0,0.5,1,2,3]$  and the pressure contour at levels $[0.3, 0.17, 0.12, 0.11, 0.09, 0.07, 0.05, 0.02, 0, -0.002]$.  Plots from the solutions of both the TN and WABE boundary conditions  are collected in Figure~\ref{fig:modifiedLDCResults}.
In referring to the results presented in \cite{BotellaPeyret98,LiuLiuPego2010b},  we find that streamlines shown in Figure \ref{fig:modifiedLDCResults} are in excellent qualitative agreement with the ones presented there.

{
\newcommand{\figWidth}{6cm}
\newcommand{\trimfig}[2]{\trimw{#1}{#2}{0.}{0.}{0.}{0.0}}
\begin{figure}[h]
\begin{center}
\begin{tikzpicture}[scale=1]
  \useasboundingbox (0.0,0.0) rectangle (14.,6.5);  

  \draw(-0.5,0.) node[anchor=south west,xshift=0pt,yshift=0pt] {\trimfig{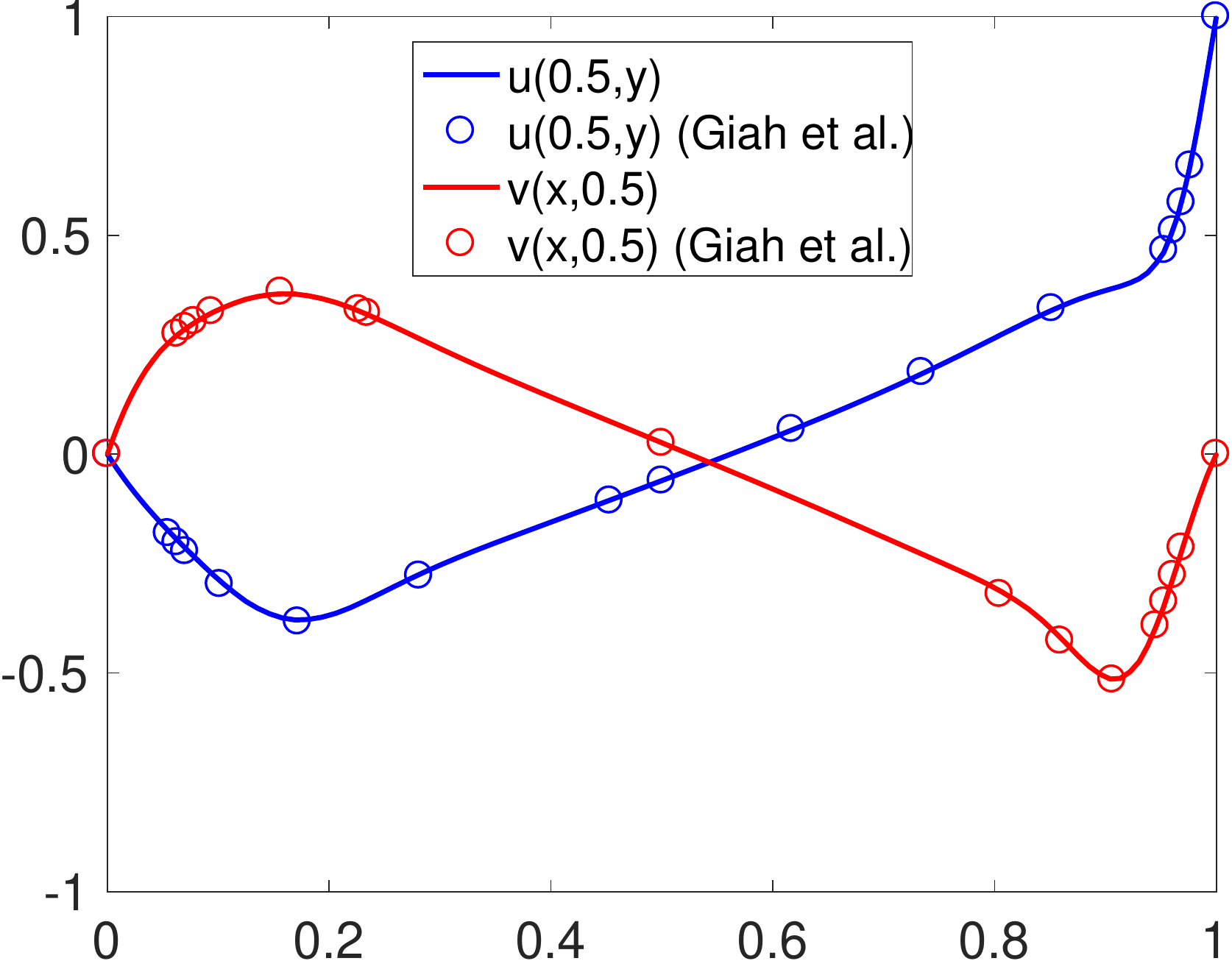}{\figWidth}};
\draw(6.8,0.) node[anchor=south west,xshift=0pt,yshift=0pt] {\trimfig{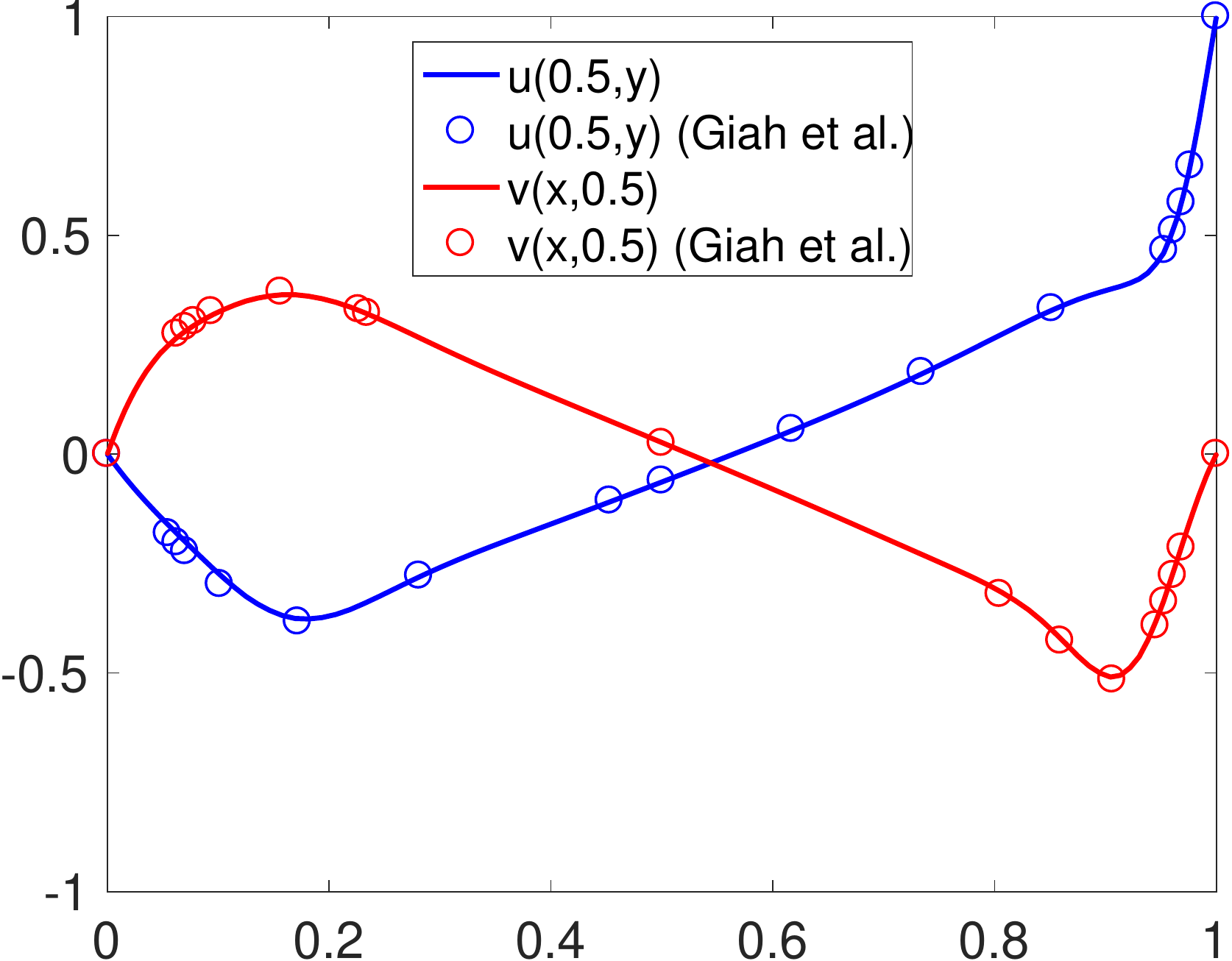}{\figWidth}};

%
\end{tikzpicture}

\end{center}
\caption{ Plots for $u(0.5,y)$ and $v(x,0.5)$ at $t=50$. Left: TN boundary condition. Right: WABE boundary condition.}\label{fig:modifiedLDCGhia}
\end{figure}
}

{
\newcommand{\figWidth}{6.cm}
\def\xa{14}
\def\ya{13.}
\newcommand{\trimfig}[2]{\trimw{#1}{#2}{0.}{0.}{0.}{0.0}}
\begin{figure}[b]
\begin{center}
\begin{tikzpicture}[scale=1]
  \useasboundingbox (0.0,0.0) rectangle (\xa,\ya);  

\draw(0.25,6.5) node[anchor=south west,xshift=0pt,yshift=0pt] {\trimfig{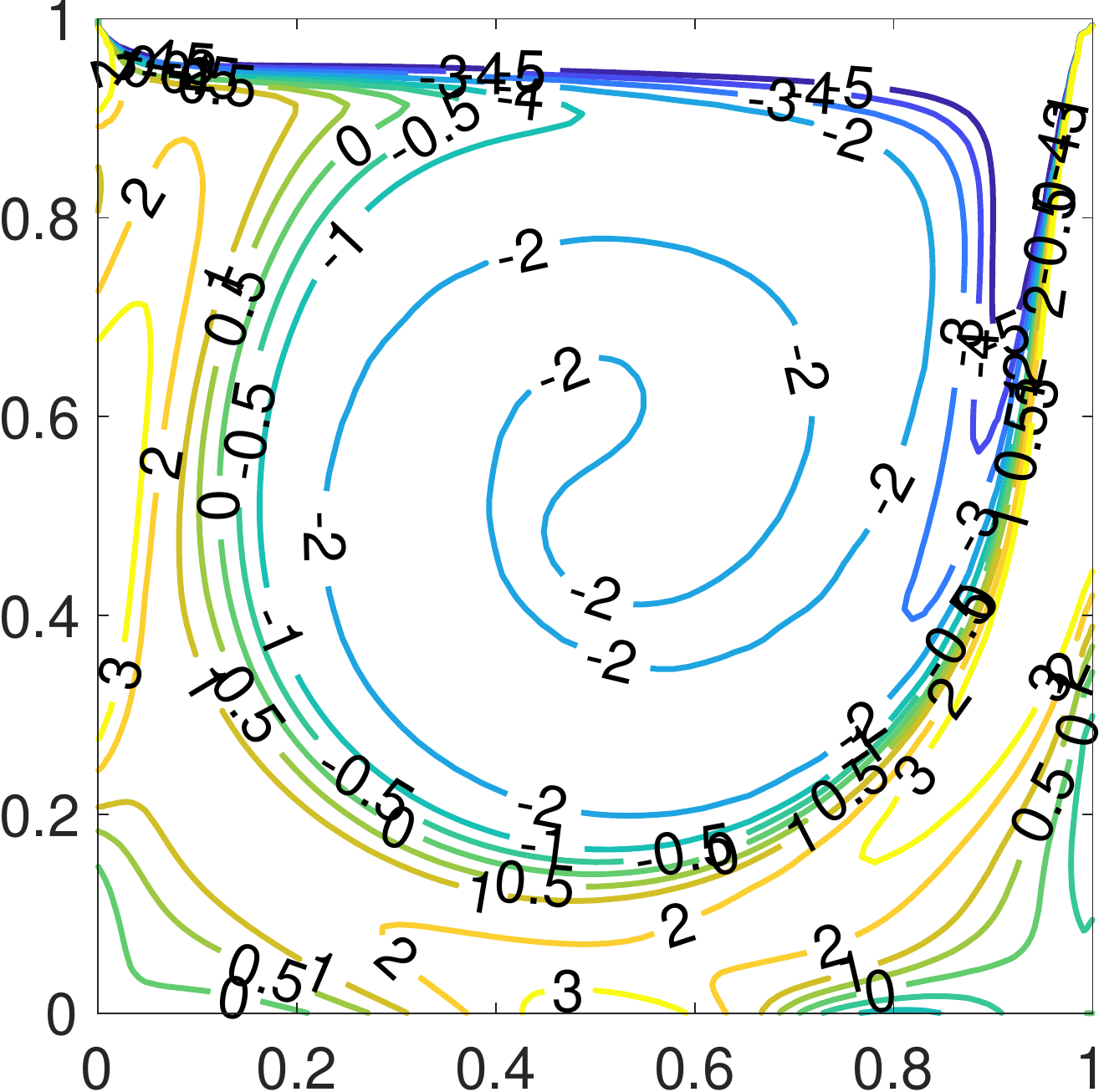}{\figWidth}};
\draw(6.8,6.5) node[anchor=south west,xshift=0pt,yshift=0pt] {\trimfig{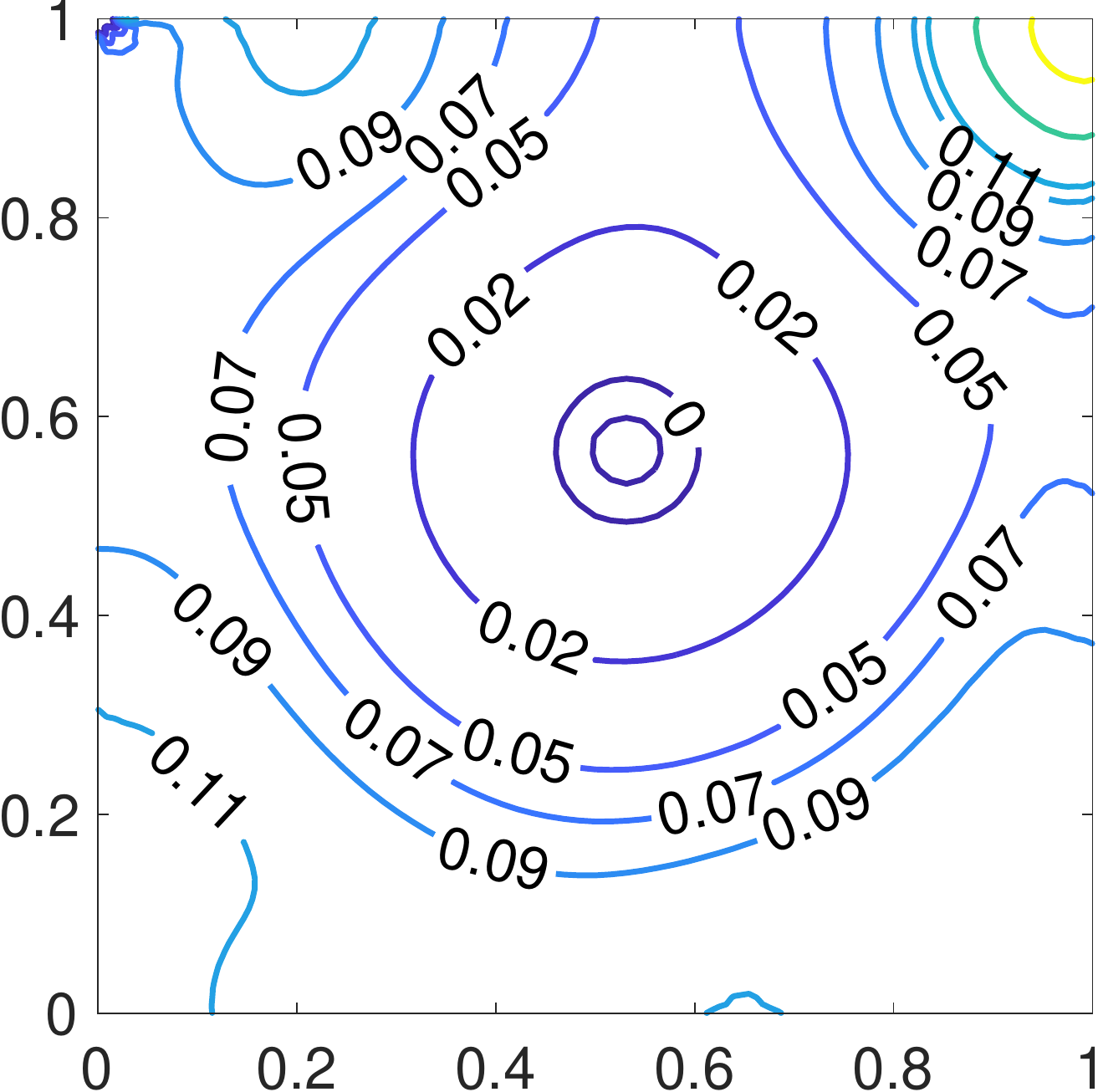}{\figWidth}};

\draw(0.25,0.0) node[anchor=south west,xshift=0pt,yshift=0pt] {\trimfig{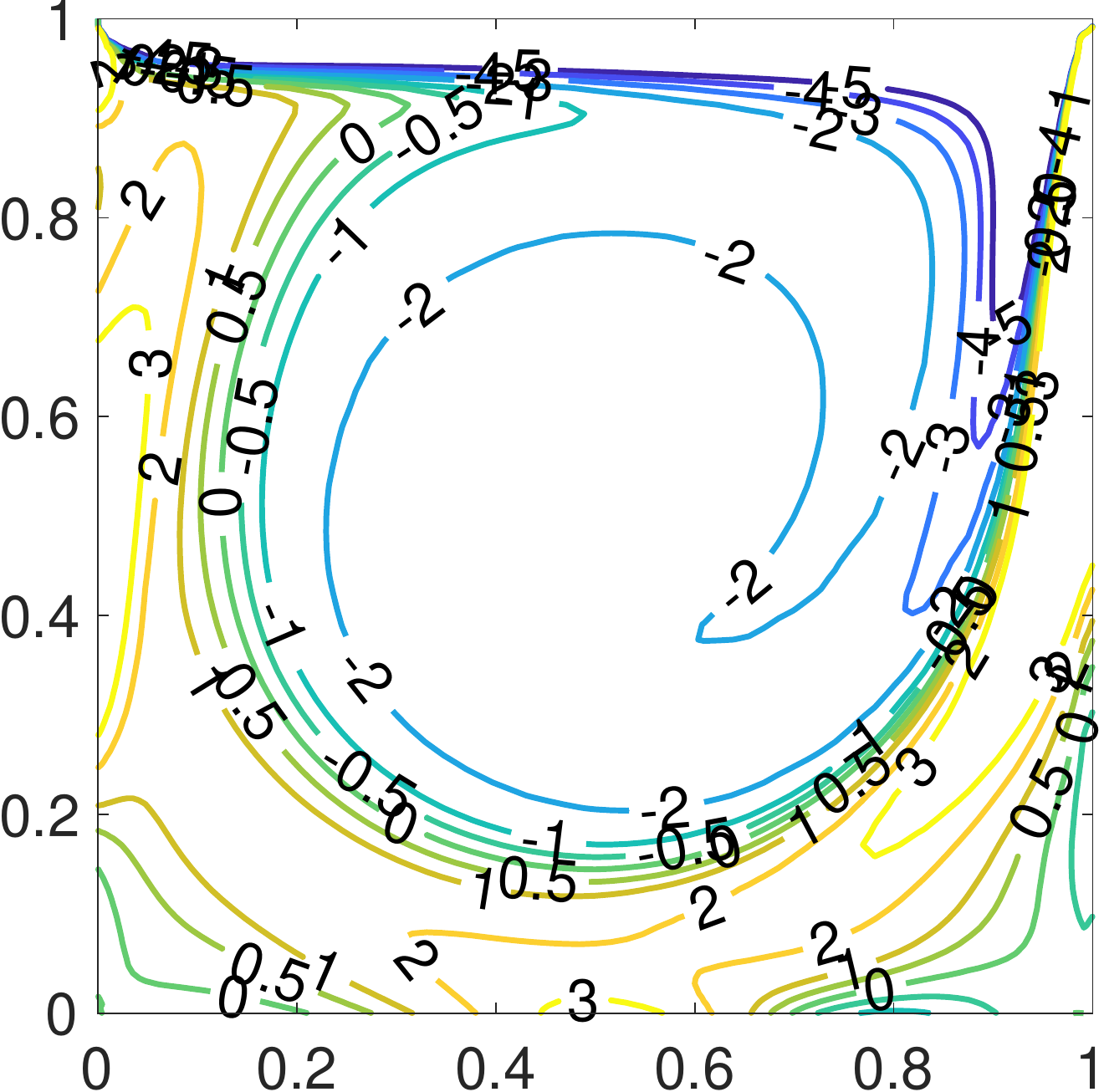}{\figWidth}};
\draw(6.8,0.0) node[anchor=south west,xshift=0pt,yshift=0pt] {\trimfig{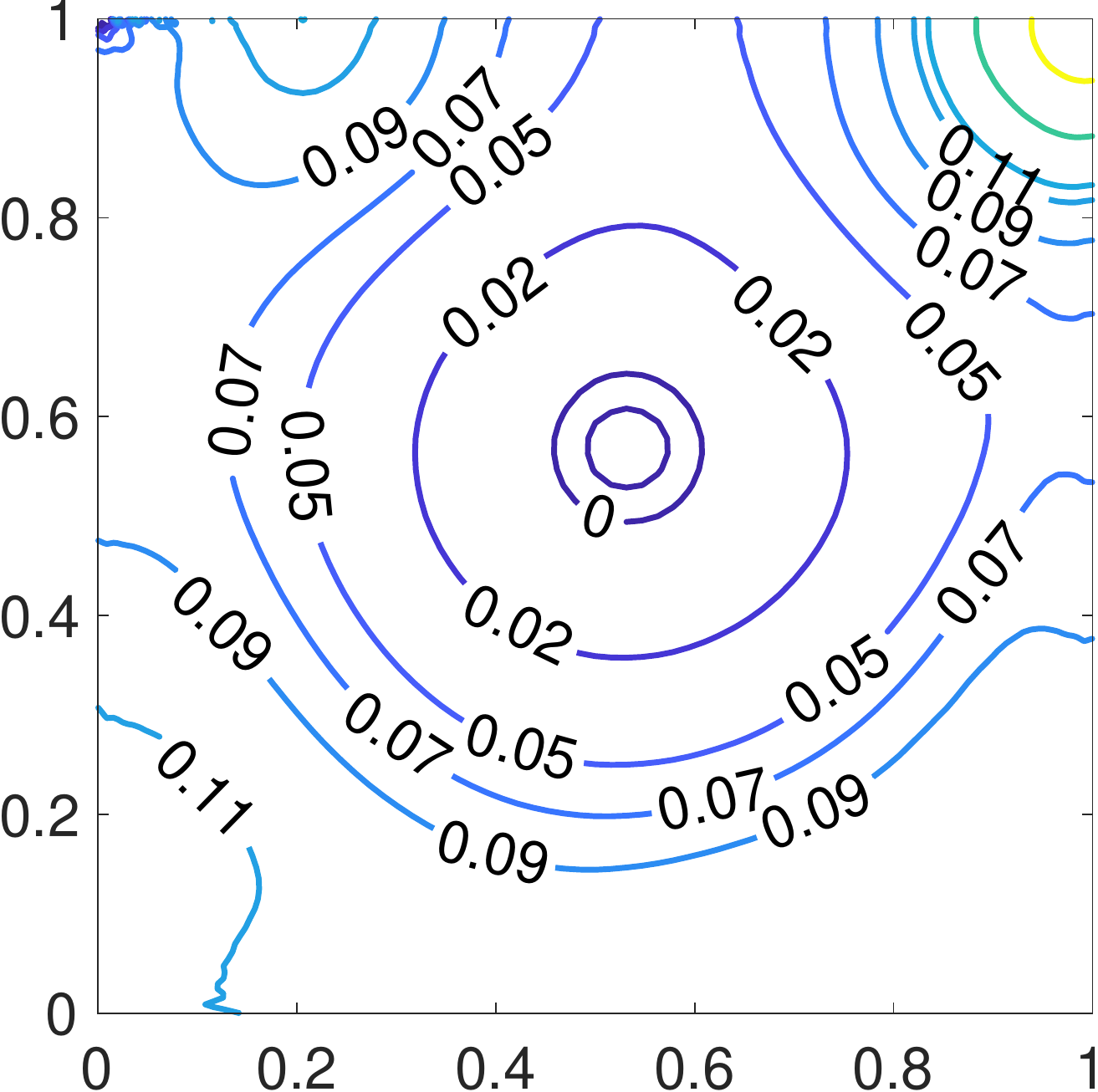}{\figWidth}};

\draw(0,3)  node[anchor=south ,xshift=0pt,yshift=0pt] {WABE:};
\draw(0,9.5)  node[anchor=south,xshift=0pt,yshift=0pt] {TN:};

%
\end{tikzpicture}

\end{center}
\caption{Vorticity  and  pressure contours at $t=50$.}\label{fig:modifiedLDCResults}
\end{figure}
}

\subsection{Flow Past a Cylinder}
As an example to illustrate the efficiency and   robustness of our approach when  used with  higher-order  finite elements,  we solve a classical flow-past-cylinder problem using $\Pe_n$ finite elements with $n=1,2,4$. The settings of the test problem follow the example in \cite{SchaferEtal96,John04,LiuLiuPego2010b}. To be specific, the domain of the problem is $\Omega=[0,2.2]\times [0,0.41] \backslash\left\{(x,y)|  (x-0.2)^2+(y-0.2)^2<0.05^2\right\}$. The inflow and outflow velocity profiles are prescribed as a time-dependent function,
$
\uv(0,y,t)=\uv(2.2,y,t)=[0.41^{-2}\sin(\pi t/8)(6y(0.41-y)),0]^T.
$
The top and bottom boundaries are enforced as no-slip walls. We  numerically solve the INS equations using our  algorithm and $\Pe_n$ finite elements in the most straightforward way; {that is, the numerical solutions are represented in \eqref{eq:spatialRep} with the basis function $\varphi_i\in \Pe_n$.}  Both the TN and WABE conditions are considered.

\begin{figure}[h]
{
  \newcommand{\figWidth}{7cm}
  \newcommand{\trimfig}[2]{\trimw{#1}{#2}{0.125}{0.095}{0.415}{0.385}}
  \newcommand{\trimfigmesh}[2]{\trimw{#1}{#2}{0.0675}{0.}{0.222}{0.03}}

\def\xa{15}
\def\ya{5}
\def\xtext{3.2}
\def\ytext{1.5}
%
\begin{center}
\begin{tikzpicture}[scale=1,line width= 1pt]
  \useasboundingbox (0.,0) rectangle (\xa,\ya);  

  \begin{scope}[xshift=0cm,yshift=2.5cm]
  \draw (-0.5, 0) node[anchor=south west] {\trimfigmesh{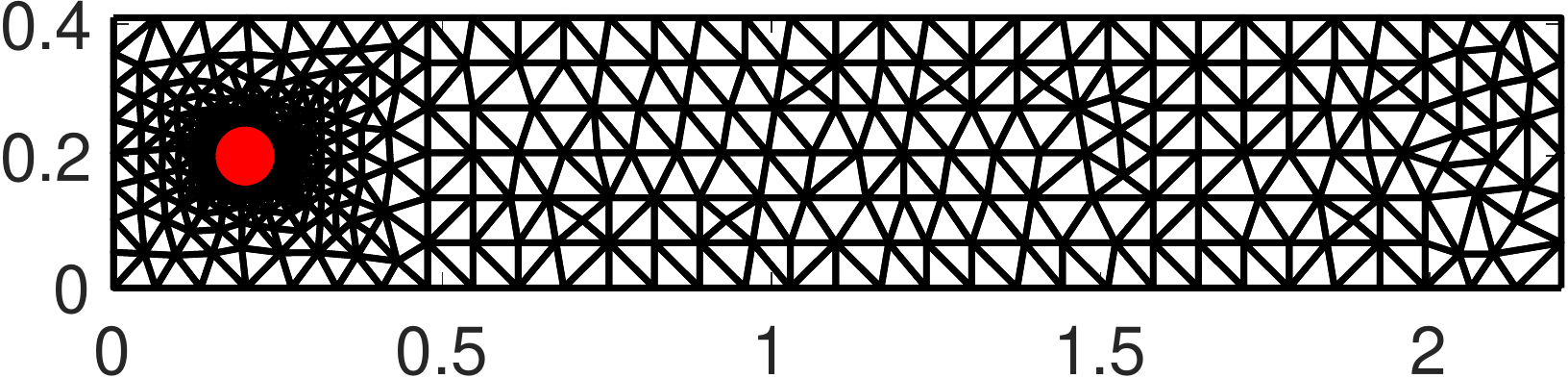}{\figWidth}};
  \draw (\xtext, \ytext) node[anchor=south west] { {$\mathcal{G}_1$:} };
\end{scope}
  
\begin{scope}[xshift=8cm,yshift=2.5cm]
  \draw (-0.5, 0) node[anchor=south west] {\trimfig{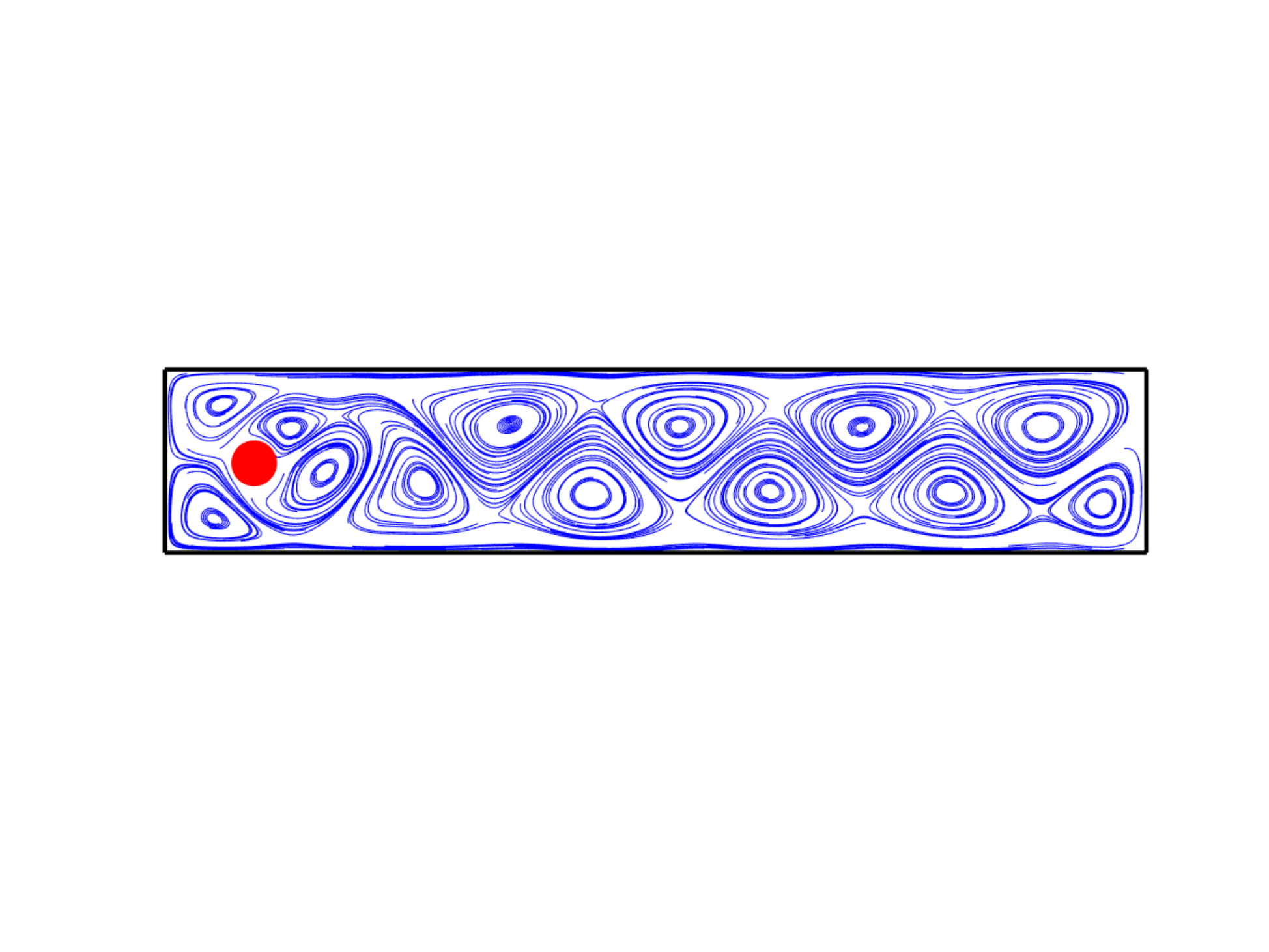}{\figWidth}};
  \draw (\xtext, \ytext) node[anchor=south west] { ($\mathbb{P}_1,\mathcal{G}_4 $): };
\end{scope}

\begin{scope}[xshift=0cm,yshift=-.2cm]
  \draw (-0.5, 0) node[anchor=south west] {\trimfig{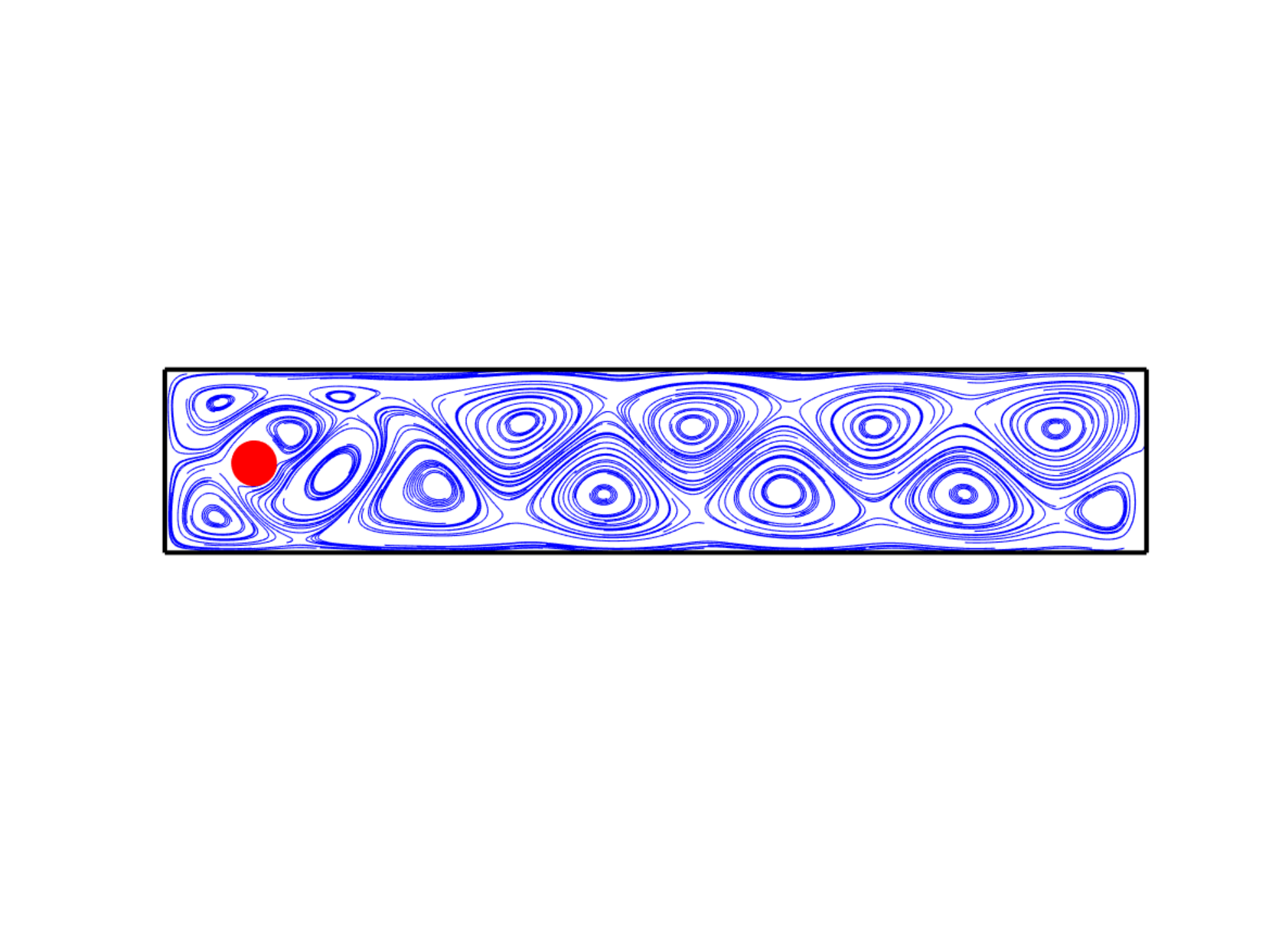}{\figWidth}};
  \draw (\xtext, \ytext)  node[anchor=south west] { ($\mathbb{P}_2,\mathcal{G}_2$): };
\end{scope}
\begin{scope}[xshift=8cm,yshift=-0.2cm]
  \draw (-0.5, 0.0) node[anchor=south west] {\trimfig{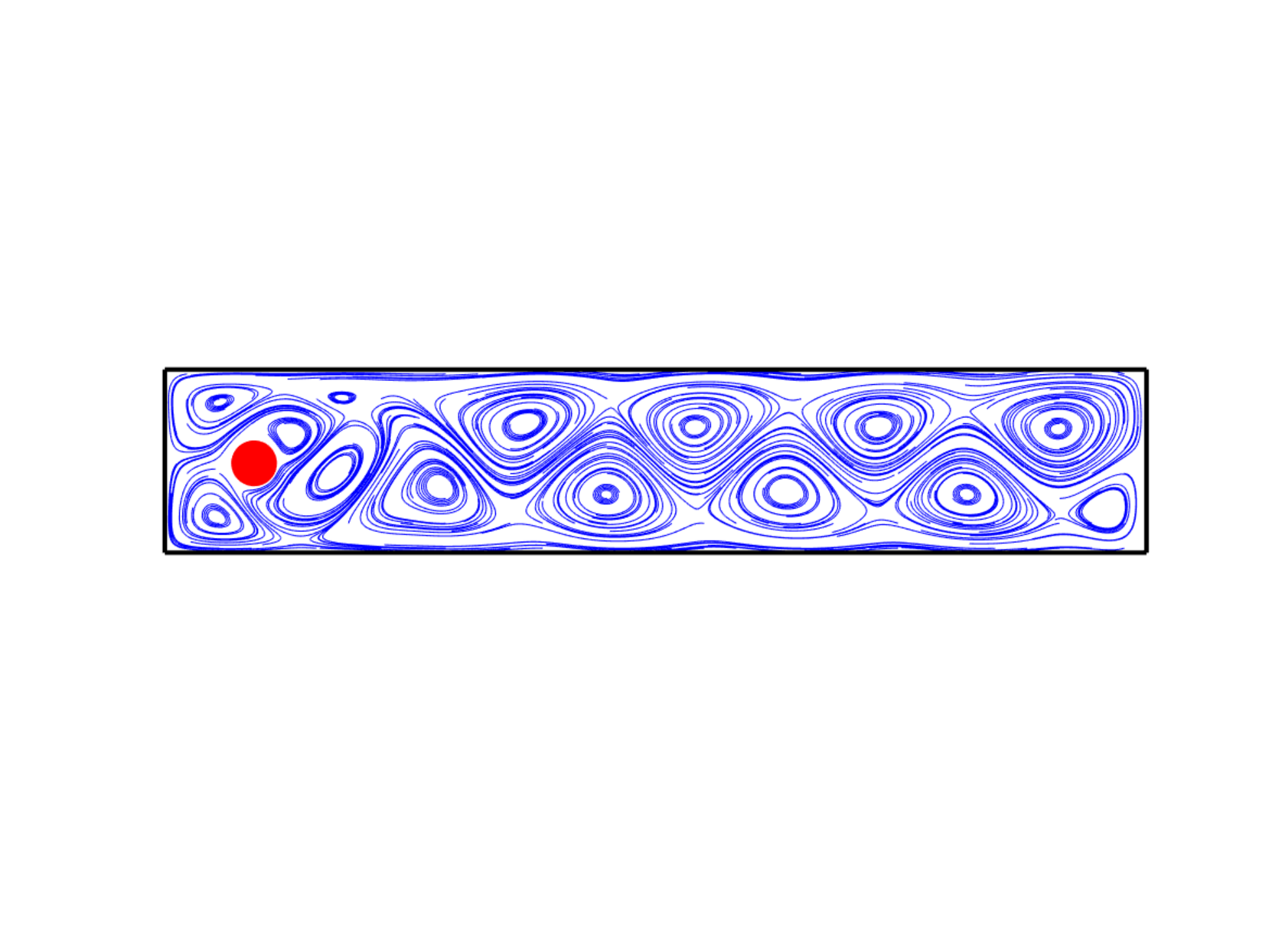}{\figWidth}};
  \draw (\xtext, \ytext) node[anchor=south west] { ($\mathbb{P}_4,\mathcal{G}_1$): };
\end{scope}
%
\end{tikzpicture}
\end{center}
}
\caption{The computational mesh (top left) and streamlines of the solutions obtained using ($\mathbb{P}_1,\mathcal{G}_4$)  (top right), ($\mathbb{P}_2,\mathcal{G}_2$) (bottom left) and($\mathbb{P}_4,\mathcal{G}_1$)  (bottom right) elements  with $\nu=1\times 10^{-3}$  at $t=8$.}
\label{fig:FPCMeshAndVorticiy}
\end{figure}
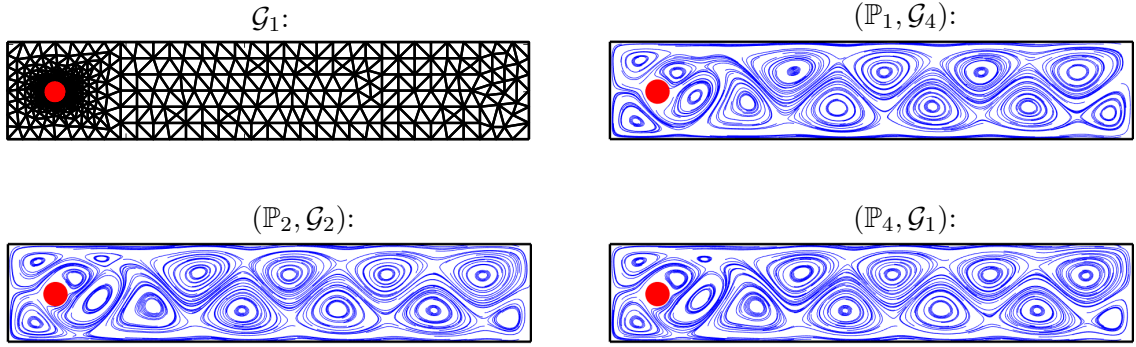

{The coarsest   computational mesh ($\mathcal{G}_1$) used to solve this problem  is shown  in the top-left image of Figure~\ref{fig:FPCMeshAndVorticiy}, which consists of  814 triangles  and  486 vertices  with the largest and smallest grid spacings  being $\max(h)=0.1055$ and $\min(h)=0.0078$, respectively. For convergence study, we refine  $\mathcal{G}_1$ by splitting  each side of the triangles into $n$ equal segments, and denote the refined mesh as $\mathcal{G}_n$.  We study the convergence with mesh refinement by solving the problem using ${\Pe}_1$ elements on $\mathcal{G}_1$,   $\mathcal{G}_2$ and  $\mathcal{G}_4$. We also study the convergence  related to the order of the elements  by   solving the problem on the same mesh ($\mathcal{G}_1$)  using   finite  elements with increasing orders ($\Pe_1$, $\Pe_2$ and $\Pe_4$).  Note that   $(\Pe_n, \mathcal{G}_k)$   indicates  the  element and mesh used  for a particular simulation. 

}

{We  show the streamlines  of the flow  from  the simulations using $\Pe_n (n=1,2,4)$ elements in  Figure~\ref{fig:FPCMeshAndVorticiy}; for comparison purposes, the results shown here maintain  the same number  (6828) of    dofs  and the same damping coefficient $\alpha=5521.08$.      We can see that these solutions are consistent with each other, and our solutions are comparable with  reference results given in  \cite{LiuLiuPego2010b}.}

{
{
\def\xa{16}
\def\ya{4.5}
\newcommand{\figWidth}{5cm}
\newcommand{\trimfig}[2]{\trimw{#1}{#2}{0.}{0.}{0.}{0.0}}
\begin{figure}[h]
\begin{center}
\begin{tikzpicture}[scale=1]
\useasboundingbox (0.0,0.0) rectangle (\xa.,\ya);  
\draw(-0.5,0.0) node[anchor=south west,xshift=0pt,yshift=0pt] {\trimfig{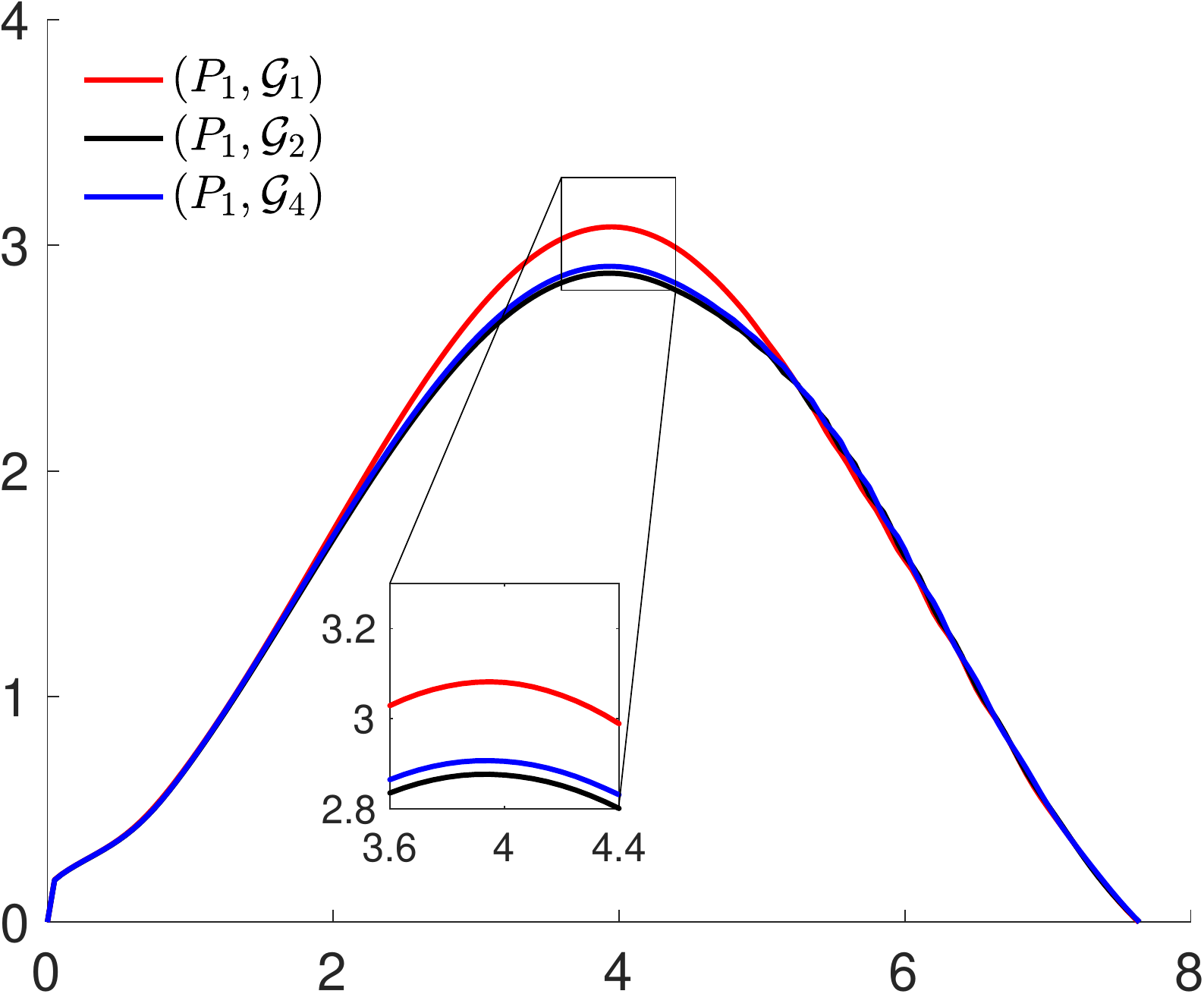}{\figWidth}};
\draw(4.5,0.0) node[anchor=south west,xshift=0.5cm,yshift=0pt] {\trimfig{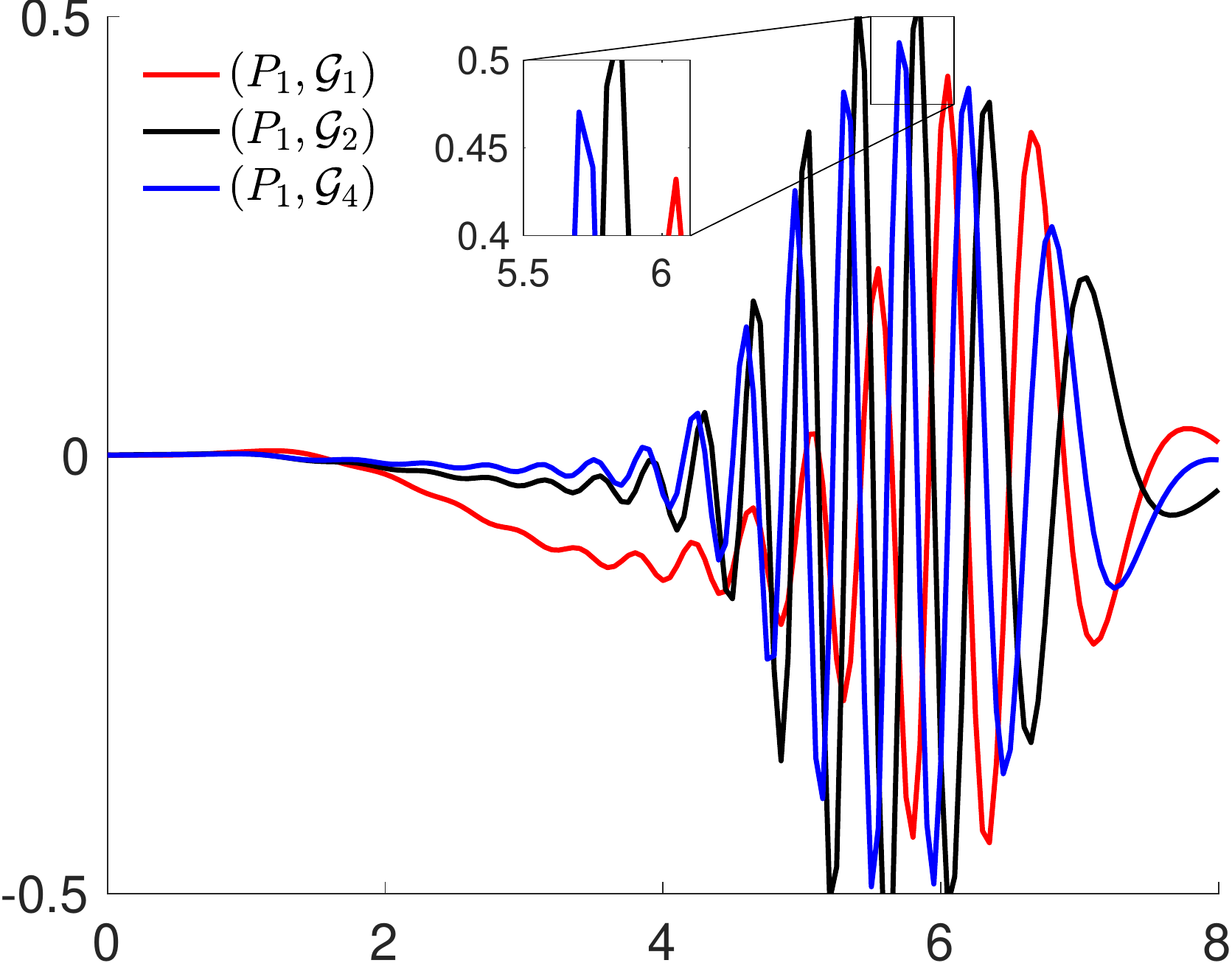}{\figWidth}};
\draw(9.5,0.0) node[anchor=south west,xshift=1.cm,yshift=0pt] {\trimfig{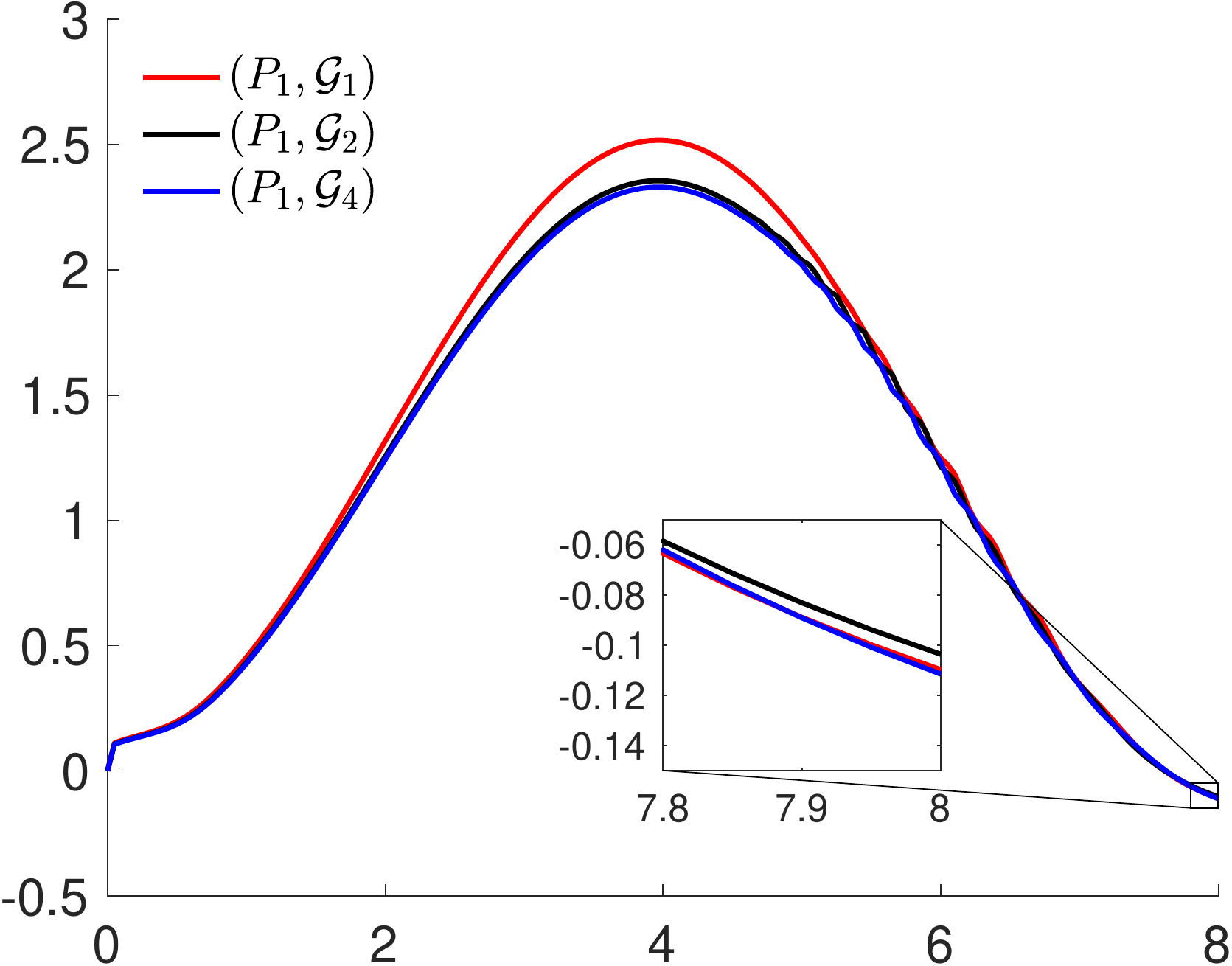}{\figWidth}};

  \draw (2.5, 4) node[anchor=south ,xshift=0pt,yshift=0pt] { $C_d(t)$};
  \draw (8.5, 4) node[anchor=south ,xshift=0.5cm,yshift=0pt] { $C_l(t)$ };
  \draw (13, 4) node[anchor=south ,xshift=1cm,yshift=0pt] { $\Delta p(t)$ };
%
\end{tikzpicture}

\end{center}
\caption{Simulations obtained using $\mathbb{P}_1$ finite elements with TN boundary condition on a sequence of refined meshes: $\mathcal{G}_1, \mathcal{G}_2$, and $\mathcal{G}_4$.}\label{fig:dragLiftDpConvP1GnTN}
\end{figure}
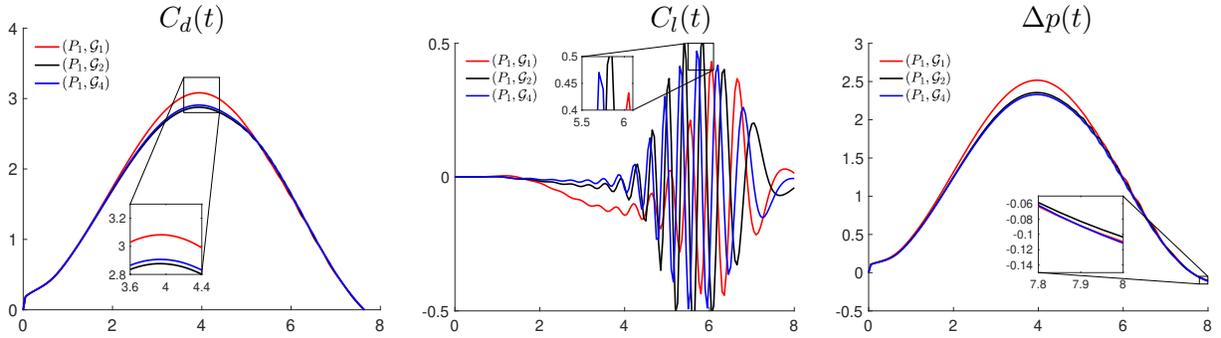
}

{
\def\xa{16}
\def\ya{4.5}
\newcommand{\figWidth}{5cm}
\newcommand{\trimfig}[2]{\trimw{#1}{#2}{0.}{0.}{0.}{0.0}}
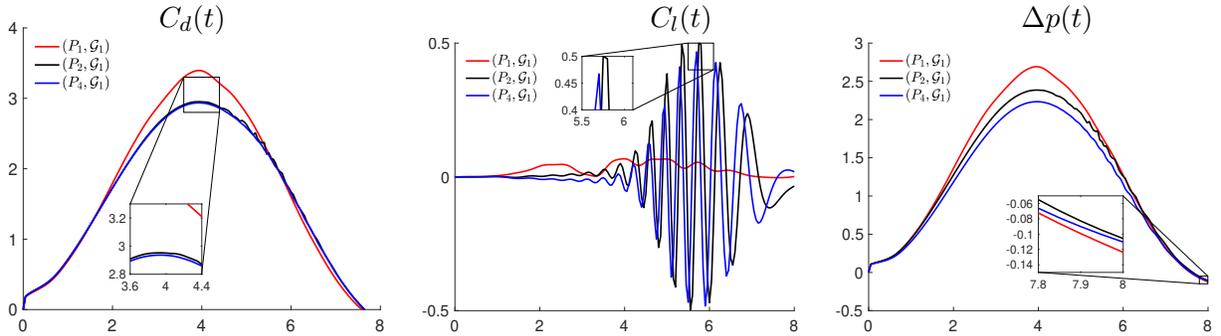
\begin{figure}[h]
\begin{center}
\begin{tikzpicture}[scale=1]
\useasboundingbox (0.0,0.0) rectangle (\xa.,\ya);  
\draw(-0.5,0.0) node[anchor=south west,xshift=0pt,yshift=0pt] {\trimfig{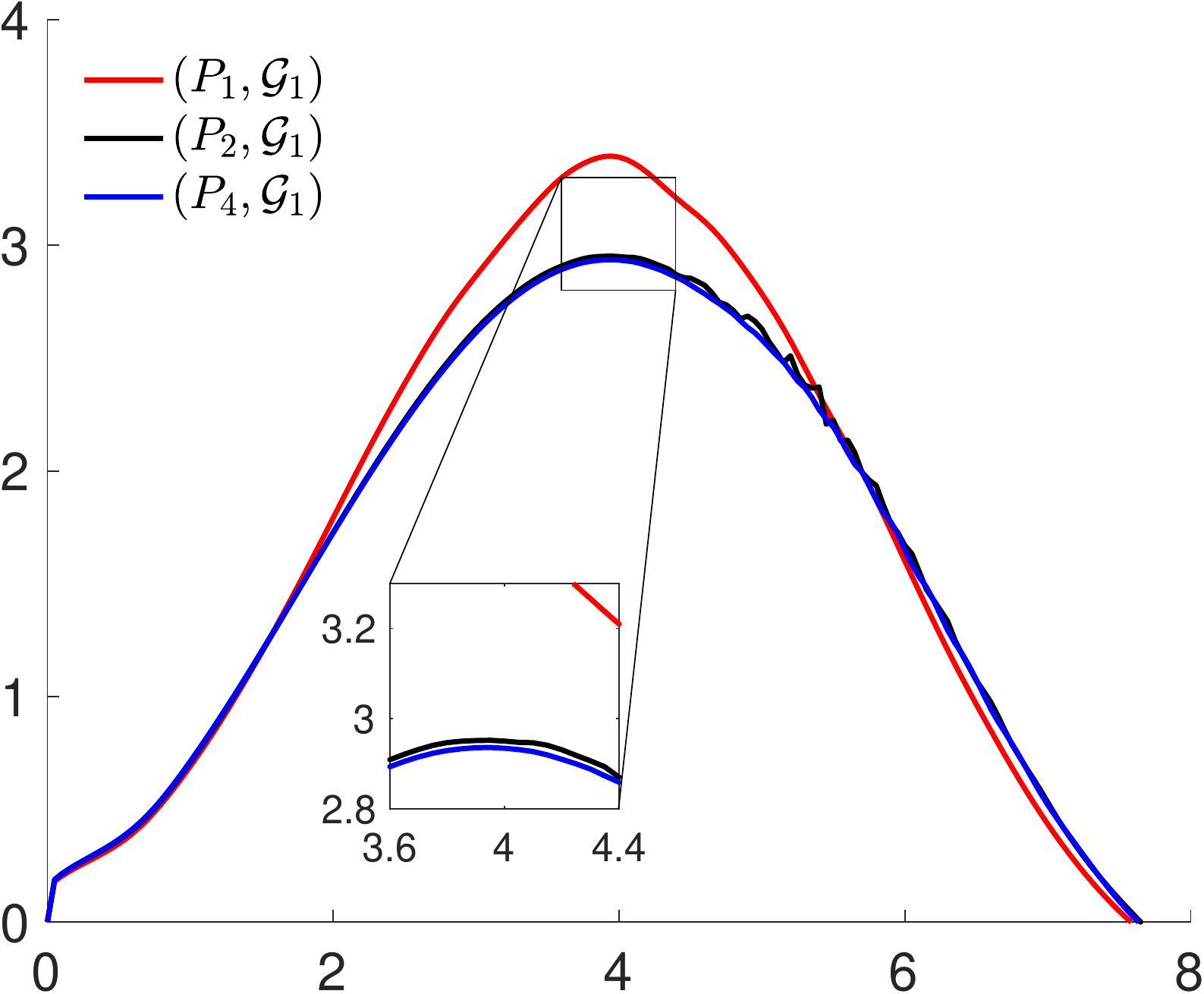}{\figWidth}};
\draw(4.5,0.0) node[anchor=south west,xshift=0.5cm,yshift=0pt] {\trimfig{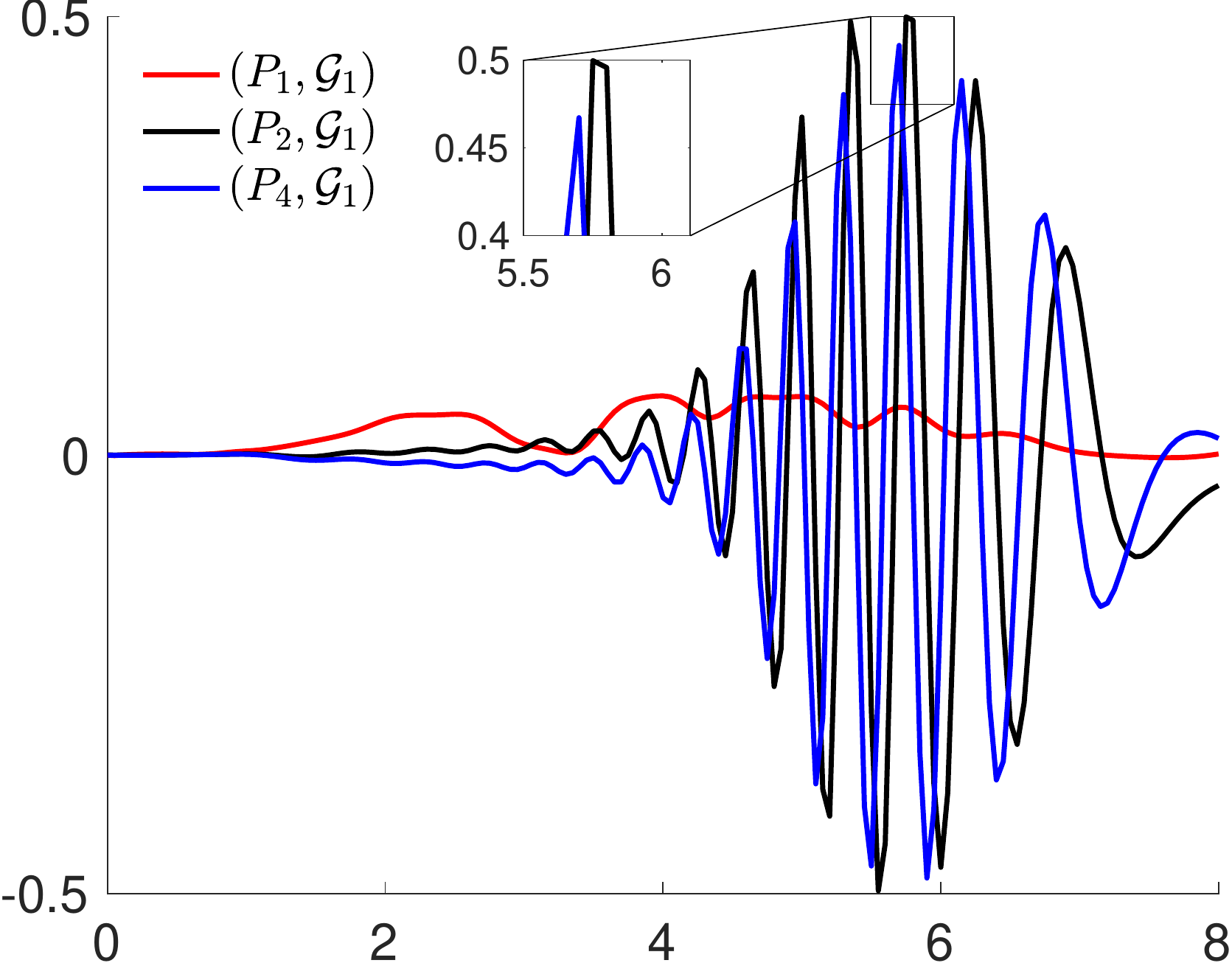}{\figWidth}};
\draw(9.5,0.0) node[anchor=south west,xshift=1.cm,yshift=0pt] {\trimfig{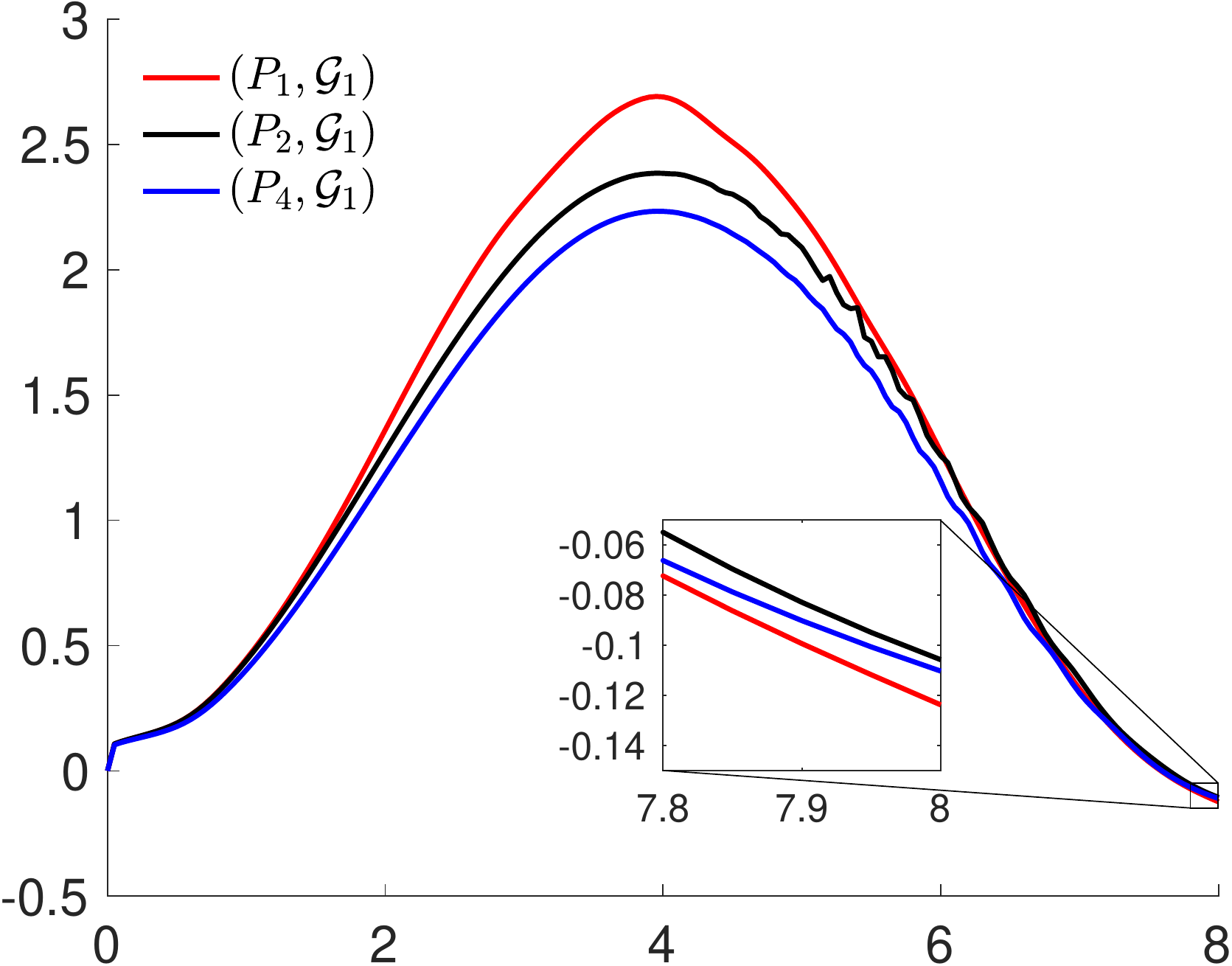}{\figWidth}};

  \draw (2.5, 4) node[anchor=south ,xshift=0pt,yshift=0pt] { $C_d(t)$};
  \draw (8.5, 4) node[anchor=south ,xshift=0.5cm,yshift=0pt] { $C_l(t)$ };
  \draw (13, 4) node[anchor=south ,xshift=1cm,yshift=0pt] { $\Delta p(t)$ };
%
\end{tikzpicture}

\end{center}
\caption{Simulations obtained using $\mathbb{P}_1$,  $\mathbb{P}_2$ and  $\mathbb{P}_4$   finite elements with WABE boundary condition on the  mesh $\mathcal{G}_1$.}\label{fig:dragLiftDpConvPnG1WABE}
\end{figure}
}
}

{
To further validate our algorithm, we compute the drag and lift coefficients at the cylinder, denoted by $C_d(t)$ and $C_l(t)$, and  the pressure difference between the front and the back of the cylinder, $\Delta p(t)$. 
The evolution of these variables  are shown in Figure~\ref{fig:dragLiftDpConvP1GnTN} and Figure~\ref{fig:dragLiftDpConvPnG1WABE}. Here  Figure~\ref{fig:dragLiftDpConvP1GnTN} collects the results  obtained using $\mathbb{P}_1$ finite elements with TN boundary condition on a sequence of refined meshes ( $\mathcal{G}_1, \mathcal{G}_2$, and $\mathcal{G}_4$), while Figure~\ref{fig:dragLiftDpConvPnG1WABE} collects the results  obtained  using $\mathbb{P}_1$,  $\mathbb{P}_2$ and  $\mathbb{P}_4$   finite elements with WABE boundary condition on the same mesh $\mathcal{G}_1$. 
We also calculate the maximum values of $C_d(t)$ and $C_l(t)$ and the times when they occur. All the  results obtained with 6828 dofs   are tabulated together with reference values in the literature in Table~\ref{tab:dragLiftDpTable}.  The sense of convergence for these quantities can be clearly seen in  Figure~\ref{fig:dragLiftDpConvP1GnTN} and Figure~\ref{fig:dragLiftDpConvPnG1WABE}. Although we use a rather coarse triangular mesh,   our results  agree quite well with  the  range of the reference values.
}


\begin{table}[h] 
\begin{center}
\begin{tabular}{rccccc} \hline 
\hline 
 & $C_{d,\max} $ & $ t(C_{d,\max}) $ & $ C_{l,\max} $ & $ t(C_{l,\max}) $ & $ \Delta p(8) $ \\ \hline 
TN  ($\mathbb{P}_1,   \mathcal{G}_4$) & $2.9070 $ & $ 3.9343 $ & $ 0.4901 $ & $ 5.7189 $ & $ -0.1115 $ \\ \hline 
TN  ($\mathbb{P}_2,   \mathcal{G}_2$) & $2.9527 $ & $ 3.9376 $ & $ 0.4871 $ & $ 5.6933 $ & $ -0.1115 $ \\ \hline 
TN  ($\mathbb{P}_4,   \mathcal{G}_1$) & $2.9371 $ & $ 3.9365 $ & $ 0.4797 $ & $ 5.6885 $ & $ -0.1111 $ \\ \hline 
WABE ($\mathbb{P}_1,   \mathcal{G}_4$) & $2.9036 $ & $ 3.9341 $ & $ 0.4663 $ & $ 5.7274 $ & $ -0.1111 $ \\ \hline 
WABE ($\mathbb{P}_2,   \mathcal{G}_2$) & $2.9415 $ & $ 3.9339 $ & $ 0.4663 $ & $ 5.6989 $ & $ -0.1126 $ \\ \hline 
WABE ($\mathbb{P}_4,   \mathcal{G}_1$) & $2.9363 $ & $ 3.9365 $ & $ 0.4738 $ & $ 5.6890 $ & $ -0.1103 $ \\ \hline 
\hline 
 \multicolumn{6}{c}{ Reference Values}   \\ \hline\hline
V. John \cite{John04} & $2.9509 $ & $ 3.9362 $ & $ 0.4779 $ & $ 5.6931 $ & $ -0.1116 $ \\ \hline 
Liu {\em et al.} \cite{LiuLiuPego2010b} & $2.9541 $ & $ 3.9364 $ & $ 0.4791 $ & $ 5.6928 $ & $ -0.1116 $ \\ \hline 
Sch\"afer {\em et al.}  \cite{SchaferEtal96} & $[2.930,2.970] $ &   & $ [0.470,0.490] $ &   & $[-0.115,-0.105] $ \\ \hline 
\end{tabular}
\caption{Maximum values of the drag and lift coefficients and the pressure difference at final time $t=8$.} \label{tab:dragLiftDpTable}
\end{center}
\end{table}

We have to note that we can not claim optimal order of accuracy is achieved for $\Pe_n$ elements with $n>1$ at this point even though we have obtained  accurate and comparable results in the  flow-past-a-cylinder example. This is because, for all the  shown computations, triangular meshes are used, which means that the cylinder is not a real cylinder but a polygon instead.  Therefore, the computational mesh contributes an $\order{h^2}$ error.   To achieve higher than second-order accuracy, we need to adapt our algorithms for isoparametric finite elements.
And this will be left for future work.

\section{Conclusion}\label{sec:conclusion}

In this paper, we propose an algorithm that solves the incompressible Navier-Stokes equations in the velocity-pressure formulation using a split-step method that separates the updates for  velocity and pressure at each time step. The separation  of the pressure solution is the key to avoid solving a saddle point problem whose  solution depends on the choice of finite-element spaces for velocity and pressure that is subject to the  LBB condition. Therefore, our algorithm has more   flexibility of choosing finite-element spaces; for efficiency and robustness,   Lagrange (piecewise-polynomial) finite elements of equal order for both velocity and pressure are used in the algorithm.

We also include a divergence damping term  into our formulation, this linear damping term plays no role at the PDE level, but helps suppress the numerical divergence in the discretized equations, and more importantly it improves the accuracy of the scheme. Motivated by a post-processing technique that produces super-convergent derivatives from finite-elememt solutions, we formulate an alternative compatibility boundary condition at the discrete level for the pressure equation. The new pressure boundary condition, referred to as WABE boundary condition, are shown to help the pressure solution achieve better accuracy near the boundary.

An important feature of the paper is that we use the normal-mode analysis, a technique that is often used for the analysis of finite difference schemes, to reveal the stability and accuracy properties of our finite element scheme via a simplified model problem. The model problem is discretized on a uniform mesh using $\Pe_1$ finite elements, so that we can rewrite the scheme as a finite difference method and then perform the  normal-mode analysis to the resulted discrete system.  The analysis shows that the scheme for the model problem is locally stable with the presence of a large divergence damping term  for both TN and WABE pressure boundary conditions. Further, by obtaining the leading order solutions of the error equations, we find that the  error introduced by the boundary condition is rapidly damped out by the divergence damping term producing boundary-layer errors. Since the
WABE boundary condition is  more accurate than the TN boundary condition in terms of truncation error,  it is expected to help  alleviate the boundary-layer errors for the pressure solution. 

Moreover, we  conduct careful numerical tests to verify the stability and accuracy  of our scheme. Through  convergence studies  using the method of manufactured
solutions, we find that the interior accuracy is
improved by including divergence damping, and the boundary accuracy is further improved by using the WABE pressure boundary  condition.  Mesh refinement study  using  $\Pe_1$ finite elements confirms that, with both divergence
damping and WABE condition, our scheme is 2nd order accurate up-to the boundary, which is optimal for the elements used.  The numerical results agree with the analysis. To further validate our scheme, benchmark problems such as lid-drive cavity and flow-past-a-cylinder are also considered; we have shown that solutions obtained using our algorithm are in excellent agreement with those reported in {the} literature.

 
%
%

{The split-step finite-element method developed in this paper exhibits 
  good numerical properties. In particular,  the  flexibility of using the standard Lagrange (piecewise-polynomial) finite elements of equal order for both velocity and pressure, which violates the classic LBB stability condition, makes it much easier to couple our fluid solver with a structure solver for solving FSI problems. Furthermore, the ability to suppress the boundary-layer errors in the pressure solution using WABE boundary condition ensures that information can be accurately transferred across the fluid-structure interface, which could be crucial for an FSI algorithm  to maintain high-order accuracy. In the future, we will investigate extending our AMP FSI schemes using  this finite element INS algorithm  to develop an  accurate and efficient partitioned FSI algorithm within the finite element framework.

In terms of improving the INS algorithm itself,} we will use isoparametric elements to achieve higher order accuracy.
In addition to the WABE boundary condition, we will also investigate the possibility of  addressing  the pressure boundary-layer issue using $p$-refinement by increasing the polynomial degree for the basis functions on the boundary nodes, thus we could obtain more accurate approximation of $\nabla\times\uv_h$ on the boundary. This $p$-refinement strategy is also computationally efficient since the degrees of freedom
 are only increased on the boundary  with all the other basis functions in the interior being unchanged. {Analysis for this new finite element framework is also  under investigation; we are interested in deriving some energy estimates for our scheme under more general assumptions. }

\section{Acknowledgment}
The author  would like to acknowledge Professors  W. D. Henshaw, Fengyan Li, J. W. Banks and   D. W. Schewendeman  of Rensselaer Polytechnic Institute (RPI)  for helpful conversations.  Portions of this research were conducted with high performance computational resources provided by the Louisiana Optical Network Infrastructure (http://www.loni.org).

\appendix
\section{}\label{sec:appendix}
In this appendix, we show technical details  for some of the results presented in Section~\ref{sec:analysis}. 

\begin{itemize}
\item Proposition~\ref{prop:errorEqns}
  \begin{proof}
Inserting the errors  $U_j=u_j-\uhat_j$, $V_j=v_j-\vhat_j$ and $P_j=p_j-\phat_j$ into the numerical scheme   \eqref{eq:fdEqn} and expanding the exact solutions at grid $j$, we arrive at the following error equations
\begin{align*}
  \text{for}~j>0:~
  \begin{cases}
 M\dot{U}_j = - ik MP_j-\nu k^2MU_j+\nu D_+D_-U_j+\order{h^2},\\
  M\dot{V}_j = -D_0P_j- \nu k^2MV_j+\nu D_+D_-V_j+\order{h^2}, \\
     -k^2MP_j+D_+D_-P_j=  \alpha ik  MU_j+ \alpha D_0V_j +\order{h^2}.
    \end{cases}
\end{align*}
Similarly, for the the TN condition, we have
\begin{equation*}
    D_+P_0 +\nu ik D_+U_0 =k^2\left(\frac{1}{3}P_0+\frac{1}{6}P_1\right)h+   \alpha ik \left(\frac{1}{3}U_0+\frac{1}{6}U_1\right)h+\alpha\frac{V_1-V_0}{2}+ \order{h},
  \end{equation*}
and, for the WABE condition, we have
\begin{equation*}
  D_+P_0+\nu ikD_+U_0=
 -\left(\frac{2}{3}\dot{V}_0+\frac{1}{3}\dot{V}_1\right)-\nu k^2\left(\frac{2}{3}{V}_0+\frac{1}{3}{V}_1\right)
   +
  \order{h^2}.
  \end{equation*}
Realizing that
$
Mf_j = f_j+\order{h^2}
$, and representing $\order{h^r}$ as $h^rF$ with some $\order{1}$ function $F$, we derive the error equations \eqref{eq:errorEqnLumped} and their boundary conditions \eqref{eq:errorBCLumped}, which completes the proof.
\end{proof}

Note that replacing $Mf_j$ with $f_j$ is in the same spirit
of mass lumping, a technique frequently used in finite-element methods. In addition, to see the second order accuracy of the WABE condition, it is most convenient  to  expand the exact solutions about the point $y=1/3$.

\item Lemma~\ref{lemma:eigenSolutionNodd}
\begin{proof}
If $\alpha=0$, the  general solution to the  pressure equation  in the eigenvalue problem \eqref{eq:eigenValueProblem} is found to be
$$
\ptilde_j=C_p e^{-\xi y_j},
$$
where $\xi$ satisfies
\begin{equation}\label{eq:noddAnalysisXiRelation}
\frac{4}{h^2}\sinh^2\left(\frac{\xi h}{2}\right)=k^2 ~~\text{and}~~\xi>0. 
\end{equation}
Substituting the pressure solution into the velocity equations, we have
\begin{align}
 & \left(s+\nu k^2\right) {\utilde}_j -\nu D_+D_-\utilde_j= - ik C_p e^{-\xi y_j}, \label{eq:noddAnalysisU}\\
 &  \left(s+\nu k^2\right) {\vtilde}_j -\nu D_+D_-\vtilde_j= \frac{1}{h}\sinh(\xi h)C_p e^{-\xi y_j}. \label{eq:noddAnalysisV}
\end{align}
We note that both velocity equations have the same homogeneous part, which is solved by
$$
{\utilde}^h_j ={\vtilde}^h_j= e^{-\gamma y_j}.
$$
Here $\gamma$ satisfies
\begin{equation}\label{eq:noddAnalysisGammaRelation}
\frac{4}{h^2}\sinh^2\left(\frac{\gamma h}{2}\right)=\frac{s}{\nu}+ k^2~~\text{and}~~\Re(\gamma)>0.
\end{equation}
Note that $\gamma=\gamma(s)$ depends on $s$.
The  particular solutions of  the velocity equations have the forms
$$
\utilde^p_j=A_ue^{-\xi y_j}~~\text{and}~~ \vtilde^p_j=A_ve^{-\xi y_j}.
$$
Substituting them into \eqref{eq:noddAnalysisU} and \eqref{eq:noddAnalysisV}, respectively,  we get
\begin{align*}
 & A_u\left[\left(s+\nu k^2\right)  -\nu\frac{4}{h^2}\sinh^2\left(\frac{\xi h}{2}\right) \right]= - ik C_p,\\
 &A_v\left[  \left(s+\nu k^2\right) -\nu\frac{4}{h^2}\sinh^2\left(\frac{\xi h}{2}\right) \right]= \frac{1}{h}\sinh(\xi h)C_p.
\end{align*}
Using \eqref{eq:noddAnalysisXiRelation}, we have
$$
A_u = -\frac{ ik C_p}{s} ~~\text{and}~~ A_v=\frac{1}{hs}\sinh(\xi h)C_p.
$$
The general solutions of \eqref{eq:noddAnalysisU} and \eqref{eq:noddAnalysisV} are then given by
\begin{align*}
  &\utilde_j= \utilde^p_j+C_u {\utilde^h}_j=A_ue^{-\xi y_j}+C_u e^{-\gamma y_j},\\
  &\vtilde_j= \vtilde^p_j+C_v {\vtilde^h}_j=A_ve^{-\xi y_j}+C_v e^{-\gamma y_j}.
\end{align*}
Implementing the no-slip boundary condition, $\utilde_0=\vtilde_0=0$, we have
$$
C_u=-A_u~~\text{and}~~C_v=-A_v.
$$
Therefore, we have found the solution given in \eqref{eq:noddAnalysisSolutions}

\end{proof}

\item Lemma~\ref{lemma:noddLemma_q1}
\begin{proof}
  If $q_1(s)$ is real, there is a real number $c$ such that
  $$
  q_1(s)=\frac{1}{s}\left(e^{-\xi h}- e^{-\gamma h}\right)=c.
  $$
Then we have
  \begin{equation} \label{eq:noddLemma1ProofEqGamma1}
e^{-\gamma h}= e^{-\xi h}-cs,
  \end{equation}
and squaring \eqref{eq:noddLemma1ProofEqGamma1} implies
    \begin{equation} \label{eq:noddLemma1ProofEqnMain}
e^{-2\gamma h} = c^2s^2-2cse^{-\xi h}+e^{-2\xi h}
    \end{equation}
    From \eqref{eq:noddAnalysisXiRelation} and \eqref{eq:noddAnalysisGammaRelation}, we have
    \begin{align}
      e^{-2\xi h}-(2+h^2k^2)e^{-\xi h}+1=0 \label{eq:noddLemma1ProofEqXi}\\
      e^{-2\gamma h}-[2+h^2(s/\nu+k^2)]e^{-\gamma h}+1=0 \label{eq:noddLemma1ProofEqGamma2}
     \end{align}
    After inserting \eqref{eq:noddLemma1ProofEqGamma1} and  \eqref{eq:noddLemma1ProofEqGamma2} into \eqref{eq:noddLemma1ProofEqnMain} to eliminate $e^{-\gamma h}$ and $e^{-2\gamma h}$, we have
    $$
   [2+h^2(s/\nu+k^2)](e^{-\xi h}-cs)-1=   c^2s^2-2cse^{-\xi h}+e^{-2\xi h}.
    $$
   Simplifying the above equation using \eqref{eq:noddLemma1ProofEqXi}, we arrive at
   $$
   c^2s^2-2cse^{-\xi h}+(2+h^2k^2)cs-h^2s/\nu e^{-\xi h}+ch^2s^2/\nu=0.
   $$
This is a quadratic equation for $s$ with real coefficients. Obviously, $s=0$ is a root; it follows that the other root must also be real, which proves the lemma.
   \end{proof}

\item Lemma~\ref{lemma:noddLemma_q}

\begin{proof}
  If $s$ is a root for $q(s)$, then $0=q(s)=\left(e^{-\xi h}-1\right)+\nu k^2q_1(s)$   implies $q_1(s)$ is real.
  According to Lemma~\ref{lemma:noddLemma_q1}, $s$ must be real. So it suffices to consider $s>0$ to prove this lemma.

  For $s>0$, we solve  \eqref{eq:noddLemma1ProofEqXi} and \eqref{eq:noddLemma1ProofEqGamma2} and obtain
  \begin{align*}
  e^{-\xi h}&=\frac{1}{2}\left[(2+h^2k^2)-\sqrt{4h^2k^2+h^4k^4} \right],\\
   e^{-\gamma h}&=\frac{1}{2}\left[(2+h^2(s/\nu+k^2))-\sqrt{4h^2(s/\nu+k^2)+h^4(s/\nu+k^2)^2} \right]. 
    \end{align*}
We note that both  \eqref{eq:noddLemma1ProofEqXi} and \eqref{eq:noddLemma1ProofEqGamma2} 
have two roots that are reciprocal of each other.  Because of the regularity conditions at infinity, the roots with magnitude less than one are kept in the above expressions. Therefore, we have
$$
q_1(s)=\frac{1}{s}\left(e^{-\xi h}- e^{-\gamma h}\right)=\frac{1}{2s}\left(- h^2s/\nu-\sqrt{4h^2k^2+h^4k^4} +\sqrt{4h^2(s/\nu+k^2)+h^4(s/\nu+k^2)^2}\right).
$$
Then the derivative of $q_1(s)$ can  be obtained explicitly as following:
$$q'_1(s)=
-\frac{N_1-N_2}{N_3}
$$
where
\begin{align*}
  &N_1=2h^2s + 4h^2k^2\nu + h^4k^4\nu + h^4k^2s, \\
  &N_2= \nu\sqrt{4h^2k^2+h^4k^4 }\sqrt{4h^2( s/\nu+k^2 )+h^4(s/\nu+k^2 )^2 },\\
  &N_3=2\nu s^2\sqrt{4h^2(s/\nu + k^2)+h^4(s/\nu + k^2)^2 }.
  \end{align*}
It is easily seen that
$$
N_1^2-N_2^2=4h^4s^2 >0. 
$$
Notice that $N_1>0$ and $N_2>0$, so we have $N_1-N_2>0$. Since $N_3>0$, we conclude that $q_1'(s)<0$. So $q_1(s)$ is decreasing, and we have
$
q_1(s)<q_1(0),
$
where
$$
q_1(0):=\lim_{s\rightarrow 0^+}q_1(s)=\frac{h^4k^2 + 2h^2}{2\nu\sqrt{4h^2k^2+h^4k^4 }} -\frac{1}{2} h^2/\nu.
$$
The  limit is evaluated using L'Hospital's rule.
Furthermore, $
q'(s)=\nu k^2q'_1(s)<0
$
implies that
$
q(s) < q(0).
$
Here
$$
 q(0)=\left(e^{-\xi h}-1\right)+\nu k^2q_1(0)=-\frac{h^2k^2}{\sqrt{4h^2k^2+h^4k^4} }<0.
 $$
 Therefore, we have shown that $q(s)<0$ for $s>0$, and 
this proves the lemma.
\end{proof}

\end{itemize}

\clearpage
\bibliographystyle{elsart-num}
\bibliography{journal-ISI,henshaw,henshawPapers,fsi,LongfeiRef}

\end{document}